\documentclass[a4paper, 10pt]{article}

\DeclareMathAlphabet{\mathpzc}{OT1}{pzc}{m}{it}

\usepackage{latexsym} 
\usepackage{amsfonts} 
\usepackage{amsmath} 
\usepackage{amsthm} 
\usepackage{mathrsfs} 
\usepackage{dsfont} 
\usepackage{bbold} 
\usepackage[english]{babel} 
\usepackage{caption} 
\usepackage{epsfig} 
\usepackage{float} 
\usepackage{subfigure} 
\usepackage{psfrag} 
\usepackage{graphicx} 
\usepackage{epsfig} 
\usepackage{amsbsy} 
\usepackage{hyperref}
\usepackage{xcolor}
\usepackage{color} 
\usepackage{enumitem}

\DeclareMathOperator{\Val}{\matV}

\newtheorem{theorem}{Theorem} 
\newtheorem{thm}{Theorem}

\newtheorem{lemma}[theorem]{Lemma} 
\newtheorem*{lemma*}{Lemma} 

\newtheorem{hyp}{Hypothesis} 
\newtheorem{rmk}[theorem]{Remark}

\newcommand{\zerarcounters}{\setcounter{equation}{0}\setcounter{theorem}{0}}

\newcommand{\ZZZ}{\mathds{Z}} 
 
\newcommand{\NNN}{\mathds{N}} 
 
\newcommand{\RRR}{\mathds{R}} 
\newcommand{\TTT}{\mathds{T}}

\newcommand{\calA}{{\mathcal A}}

\newcommand{\calE}{{\mathcal E}} 
 
\newcommand{\calG}{{\mathcal G}} 
 
\newcommand{\calI}{{\mathcal I}}

\newcommand{\NN}{{\mathcal N}} 
 
\newcommand{\calP}{{\mathcal P}}

\newcommand{\TT}{{\mathcal T}} 
 
\newcommand{\VV}{{\mathcal V}} 
\newcommand{\calW}{{\mathcal W}}

\newcommand{\FFF}{{\mathscr F}}

\newcommand{\calmP}{{\mathscr P}}

\newcommand{\gotm}{{\mathfrak m}}

\newcommand{\gotB}{{\mathfrak B}} 
\newcommand{\gotC}{{\mathfrak C}} 
\newcommand{\gotD}{{\mathfrak D}}

\newcommand{\gotK}{{\mathfrak K}} 
 
\newcommand{\gotM}{{\mathfrak M}}

\newcommand{\gotR}{{\mathfrak R}} 
 
\newcommand{\gotT}{{\mathfrak T}}

\newcommand{\matV}{{\mathscr V}}

\newcommand{\io}{\infty} 
\newcommand{\e}{\varepsilon} 
\newcommand{\al}{\alpha} 
 
\newcommand{\be}{\beta}

\newcommand{\ga}{\gamma} 
\newcommand{\om}{\omega}

\newcommand{\f}{\varphi}

\newcommand{\nn}{\nu} 

\def\ins#1#2#3{\vbox to0pt{\kern-#2 \hbox{\kern#1 #3}\vss}\nointerlineskip} 

\begin{document}

\title{\bf Pseudo-synchronous solutions\\for dissipative non-autonomous systems}
\author 
{\bf  Michele Bartuccelli$^1$, Livia Corsi$^2$, Jonathan Deane$^1$ and Guido Gentile$^2$
\vspace{2mm} 
\\ \small
$^1$ Department of Mathematics, University of Surrey, Guildford, GU2 7XH,
UK 
\\ \small
$^2$ Dipartimento di Matematica e Fisica, Universit\`a Roma Tre, Roma, I-00146,
Italy
\\ \footnotesize
E-mail: m.bartucelli@surrey.ac.uk, lcorsi@mat.uniroma3.it, j.deane@surrey.ac.uk, gentile@mat.uniroma3.it
}

\date{} 
 
\maketitle 
 
\begin{abstract}
In the framework of KAM theory, 
the persistence of invariant tori in quasi-integrable systems is proved by assuming a non-resonance condition
on the frequencies, such as the standard Diophantine condition or the milder Bryuno condition.
In the presence of dissipation, most of the quasi-periodic solutions disappear and one expects, at most,
only a few of them to survive together with the periodic attractors. However, to prove that a quasi-periodic solution really exists,
usually one assumes that the frequencies still satisfy a Diophantine condition and, furthermore,
that some external parameters of the system are suitably tuned with them.
In this paper we consider a class of systems on the one-dimensional torus, subject to a periodic perturbation
and in the presence of dissipation, and show that, however small the dissipation,
if the perturbation is a trigonometric polynomial in the angles and the unperturbed frequencies satisfy a non-resonance condition of finite order,
depending on the size of the dissipation, then a quasi-periodic solution exists with slightly perturbed frequencies
provided the size of the perturbation is small enough.
If on the one hand the maximal size of the perturbation is not uniform in the degree of the trigonometric polynomial,
on the other hand all but finitely many frequencies are allowed and there is no restriction arising from the tuning of the external parameters.
A physically relevant case, where the result applies, is the spin-orbit model, which describes the rotation
of a satellite around its own axis, while revolving on a Keplerian orbit around a planet, in the case in which
the dissipation is taken into account through the MacDonald torque.
\end{abstract} 

\zerarcounters 
\section{Introduction}
\label{sec:1} 

Almost all satellites are locked in an orbital resonance with their planets: the period of revolution is commensurable with the
period of rotation. In many cases the two periods are equal to each other (resonance 1:1), with the remarkable exception
of Mercury, considered as a satellite of the Sun (Mercury is entrapped in a 3:2 resonance).

A simple model widely used to study the problem is the \emph{spin-orbit model},
in which the satellite is described as an ellipsoidal rigid body orbiting around its primary in a Keplerian orbit \cite{D,GP,MD}.
Thus, the only variable is the sidereal angle $\theta$ that the longest axis of the satellite forms with respect
to the line of the apsides of the orbit. If $C$ denotes the maximal moment of inertia of the satellite,
the time evolution for $\theta$ is described
by a second order non-autonomous ordinary differential equation of the form
%
\begin{equation} \label{eq:1.1}
C \ddot\theta = \mathpzc{T}(\theta,\dot\theta,t) ,
\end{equation}
in which the dots represent derivatives with respect to time $t$
and the total torque $\mathpzc{T}(\theta,\dot\theta,t)$ takes into account
both the gravitational force -- through the \emph{triaxial torque}, depending only on $\theta$ and $t$ --
and the dissipation -- through the \emph{tidal torque}.

The appearance of attractors is ultimately related to the presence of dissipation.
While there is a universally accepted expression for the triaxial torque, 
the expression for the tidal torque is much more tricky.
In the literature, up to very recent times, the MacDonald model \cite{Mac,GP,CL,P,CC4,BDG2}
was the model mainly used for the tidal torque.
However, at least in some cases, the applicability of the MacDonald model appears doubtable. 
From a physical ponit of view, the main drawback of the MacDonald model is that it predicts the existence of
a pseudo-synchronous orbit, that is a solution in which the spin rate is not exactly a rational multiple of the mean motion
of the satellite. Another disappointing aspect of the model is that, 
in some instances, such as that of the Mercury,
the probability of capture of the 3:2 resonance is very low with respect to those of the coexisting
pseudo-synchronous orbit \cite{GP,CC4,BDG1}.
This runs counter to physical expectations;
only assuming a chaotic evolution of Mercury can the probability of capture become higher \cite{CL}.

Recently, in the case of Mercury, the validity of the MacDonald model -- as well as of any model based on constant time lag -- was
strongly questioned and a more realistic model was proposed \cite{EW,E,WE,EM,ME,NFME}.
The probability of capture of the attractor close to the resonance 3:2, as predicted by such a model, is much higher
with respect to the MacDonald model \cite{M}. There is strong numerical evidence that the orbit is still pseudo-synchronous,
but with spin rate much closer to the value $3/2$, and that it corresponds not to a KAM torus, but
to a Lagrangian torus inside a resonant gap \cite{BDG4}. By contrast, astronomical observations suggest that the orbit
is purely periodic and with smaller amplitude with respect to the realistic model.
A possible explanation of this apparent discrepancy is that the large-amplitude pseudo-synchronous orbit 
originally captured the satellite, which, much later in time, because of modifications in its internal structure,
evolved into the small-amplitude libration observed nowadays:
despite the fact that the probability of capture of the final resonant attractor was very low, the satellite
was already entrapped in its  (small) basin of attraction when the nearby pseudo-synchronous orbit disappeared.
We stress that, as widely pointed out in the literature \cite{E,WE,EM,ME,FM,NFME},
the model to use for the tidal torque strongly depends on the system one is interested in;
moreover, in order to detect the ultimate fate of a satellite, it is important to keep track of how the dissipation evolves in time,
considering the huge time scales involved \cite{BDG1,BDG2,WBDG1,WBDG2}.

Although the MacDonald model turns out not to be right, at least for Mercury,
the MacDonald model and the realistic model have one feature in common:
both predict the existence of a pseudo-synchronous orbit.
From a mathematical point of view, a rigorous proof of the existence of a solution which is not periodic is delicate.
In fact, one has to deal with a small divisor problem and, typically, some Diophantine condition must be assumed
on the two-dimensional frequency vector; in particular this means that one looks for a quasi-periodic solution.

The MacDonald model, because of its simplicity, is well suited for analytical computations
(the tidal torque of the realistic model, which is only $C^1$ and is characterised by narrow, large amplitudes
`spikes' at certain values of $\dot\theta$,
additionally make exhaustive numerical computations much harder \cite{BDG3,BDG4}).
The pseudo-synchronous attractor is expected to reduce to a KAM torus in the absence of dissipation,
so that one may use KAM techniques for maximal tori.
However, in order to prove that a quasi-periodic solution exists, one has to adjust an external parameter
(depending on the eccentricity), which appears in the MacDonald torque,
as a function of the frequency vector \cite{CC4,C};
for similar results see also refs.~\cite{BST,CD} for the dissipative standard map
and refs.~\cite{SL,M} for generalisations of the spin-orbit model to higher dimensions.
In principle one could make the parameter equal to the physical one by
using some implicit function argument, but this is hindered by the parameter not being a smooth
function of the frequency vector (in fact, it is defined on a set full of holes). So, with KAM arguments, one can prove that
there are quasi-periodic solutions provided the external parameter is suitably tuned;
whether the physical values of the parameters are compatible with the tuning in general remains undecidable.

In this paper we aim to study what happens when one looks for a pseudo-synchronous solution, 
requiring no condition at all on the external parameter and
assuming a non-resonance condition as weak as possible on the frequency vector.
The spin orbit model \eqref{eq:1.1}, with the MacDonald torque, is a particular case of the class of
dynamical systems described by the ordinary differential equation
\begin{equation} \label{eq:1.2}
\ddot\theta + \ga_0 \, ( \dot \theta - \al ) + \e \, f(\theta,t) = 0,
\end{equation}
where $\theta \in \TTT=\RRR / 2\pi\ZZZ$, $\e$ is the \emph{perturbation parameter} (measuring the size of the triaxial torque), 
$\ga_0$ is the \emph{dissipation parameter} and $\al$ is the aforementioned \emph{external parameter} \cite{GP,CL,CC4}.
The function $f$ depends periodically on both $\theta$ and $t$, and is real analytic.

Most physically relevant is the weakly dissipative case; so we allow $\ga_0$ to be arbitrarily small,
by setting $\ga_0=\e^{\gotm} \ga$, with $\ga$ independent of $\e$ and $\gotm\in\NNN$.
It is quite natural to assume $f$ to have zero mean, since it is, up to the sign, the derivative with respect to $\theta$
of a potential energy which depends periodically on $\theta$.
The main assumption we make, additionally, on $f$ is that it is a trigonometric polynomial of finite order $N$.
Usually the Fourier sum of $f$ is truncated to some finite order and one expects the remainder
to give negligible corrections. In the KAM approach, to deal with full Fourier series or
trigonometric polynomials makes no real difference \cite{CC4}. By contrast, for us this a crucial hypothesis.

Set, here and henceforth, $|a|:=|a|_1=|a_1|+|a_2|$ if $a=(a_1,a_2)\in\RRR^2$, and
let $\cdot$ denote the standard scalar product in $\RRR^2$, i.e.~$a\cdot b = a_1b_1+a_2b_2$
for any $a=(a_1,b_1)$ and $b=(b_1,b_2)$ in $\RRR^2$.

\begin{hyp} \label{hyp:1}
The function $f$ is a zero-mean trigonometric polynomial of degree $N$, i.e.~its Fourier components $f_{\nu}$
vanish for all $\nu\in\ZZZ^2$ such that $|\nu|>N$. Set $\Phi \! :=\! \max\{ |f_\nu| \! : \! 0 \! < \! |\nu| \! \le \! N\}$.
\end{hyp}

\begin{hyp} \label{hyp:2}
The vector $\om:=(\al,1)$ is non-resonant up to order $4\gotm N$, i.e.~one has $\om\cdot\nu \neq 0$
$\forall \nu\in \ZZZ^2$ such that $0 < |\nu| \le 4 \gotm N$, for some $\gotm\in\NNN$.
\end{hyp}

\begin{thm} \label{thm:1}
Consider \eqref{eq:1.2}, with $\ga_0=\e^{\gotm}\ga$ for some $\gotm\in\NNN$, and assume Hypotheses \ref{hyp:1} and \ref{hyp:2}.
Then there exists $\e_0>0$, depending on $\gotm$, $N$, $\ga$, $\al$ and $\Phi$,
such that for all $|\e|\le \e_0$ there exists an analytic multi-periodic
solution to \eqref{eq:1.2} with frequency vector $\om_0=(\al_0,1)$,
with $|\al-\al_0| \le C |\e|^2$, for a suitable constant $C$.
\end{thm}

\begin{rmk} \label{rmk:1.1}
\emph{
The value $\e_0$ 
goes to 0 when either $\gotm$ or $N$ tends to infinity.
That $\e_0$ tends to vanish with $\gamma$ is expected:
without requiring any stronger non-resonance condition on the frequency vector,
such as the standard Diophantine or the Bryuno condition, the invariant torus breaks up in the absence
of dissipation, as numerical evidence \cite{LL,T,O,CCV} and analytical results \cite{F,B} suggest. On the contrary,
that $\e_0$ becomes zero when $f$ is a generic analytical function
might be a technical issue of our approach and we leave as an open problem whether Theorem \ref{thm:1}
can be extended to analytical torques containing an infinite number of harmonics.
}
\end{rmk}

For fixed values of $\e$ and $\ga_0$, when writing $\ga_0=\e^{\gotm}\ga$ for some $\gotm\in\NNN$,
it may happen that the constant $\ga$ be rather large. However, 
Theorem \ref{thm:1} may be straightened by only requiring $\ga_0$
to be of the form $\ga_0=|\e|^{s} \bar\ga$, with $\bar\ga$ independent of $\e$ and $s$ an arbitrary positive constant.
To do this, one writes $\ga_0=\e^{\gotm}\ga$, where $\gotm:=\lceil s \rceil$, so that $s=\gotm - a$, for some $a\in[0,1)$.
Then the existence of a multi-periodic solution follows from Theorem \ref{thm:1}, proving that $\e_0$ does not depend
on $a$ as long as $0 \le a <1$. This leads to the following extension of Theorem \ref{thm:1}.

\begin{thm} \label{thm:2}
Consider \eqref{eq:1.2}, with $\ga_0=|\e|^{s}\bar\ga$ for some $s>0$, and assume Hypotheses \ref{hyp:1} and \ref{hyp:2},
with $\gotm=\lceil s \rceil$.
Then there exists $\e_0>0$, depending on $\gotm$, $N$, $\bar\ga$, $\al$ and $\Phi$,
such that for all $|\e|\le \e_0$ there exists an analytic multi-periodic
solution to \eqref{eq:1.2} with frequency vector $\om_0=(\al_0,1)$,
with $|\al-\al_0| \le C |\e|^{2}$, 
for a suitable constant $C$.
\end{thm}

\begin{rmk} \label{rmk:1.2}
\emph{
Theorem \ref{thm:2} shows that, assuming only that the dissipation parameter $\ga_0$ is not too large with respect to $\e$
(essentially it must be sublinear in $\e$),
for $\e$ small enough there exists a multi-periodic solution with frequency vector close to $(\al,1)$.
Typically, in physical applications, $\ga_0$ is (significantly) smaller than $\e$, so one may take $\bar\ga=1$ and fix $s$ consequently. 
}
\end{rmk}

\begin{rmk} \label{rmk:1.3}
\emph{
Theorems \ref{thm:1} and \ref{thm:2} leave out what happens when $\al$ is a rational number $p/q$, with 
the two natural numbers $p$ and $q$ not too large (more precisely such that one has $p+q\le 4\gotm N$).
Although this is a very unlikely situation -- since it corresponds to a set of values for the external parameter
which not only has zero measure but consists in a finite number of points --
however, the estimates provided by the proof of Theorem \ref{thm:1} are not uniform on $\al$:
in fact, the closer $\al$ is to such a rational, the smaller the value of $\e_0$.
In any case, when $\al$ is equal -- or even close -- to a rational number $p/q$, with $q$ small,
periodic solutions with frequency $p/q$ are known to exist \cite{BDG1,BC1}.
In such a case, a periodic and a quasi-periodic attractor may coexist:
if $\al$ is very close to a rational $p/q$, when $q$ is small, the periodic attractor is expected to be dominant,
while the opposite occurs when $q$ is large: compare the cases of the the Moon, with $\al$ very close to 1,
and Mercury, with $\al$ not quite close to $3/2$, with the continued fraction expansions giving
$\al=[1;55,1,1,4,\ldots]$ and $\al=[1;3,1,9,1,\ldots]$, respectively \cite{CC3,BDG2}.
}
\end{rmk}

\begin{rmk} \label{rmk:1.4}
\emph{
In refs.~\cite{CC4,SL,M}, for $\om_0=(\al_0,1)$ fixed to be Diophantine, a suitable function $\tilde\al(\al_0,\e)$
is proved to exist such that \eqref{eq:1.2}, with $\al:=\al_0+\e \tilde \al(\e,\al_0)$, admits
a quasi-periodic solution with frequency vector $\om_0$.
However, the map $\al_0 \mapsto \tilde \al(\e,\al_0)$ is not smooth, so, for given $\al$, there is no way to establish whether
the equation $\al=\al_0+\e \tilde \al(\e,\al_0)$ admits a solution $\al_0$ such that $\om_0$ is Diophantine.
On the contrary, under the Hypothesis \ref{hyp:2} on $\al$, the function $\tilde\al(\e,\al_0)$ is found to be smooth,
so that the implicit function theorem applies.
}
\end{rmk}

\begin{rmk} \label{rmk:1.5}
\emph{
Theorems \ref{thm:1} and \ref{thm:2} state the existence of a multi-periodic solution
(which describes an invariant torus if it is quasi-periodic), but says nothing about its possible attractiveness.
In refs~\cite{SL,M}, the quasi-periodic solution with Diophantine frequency vector is obtained by a KAM algorithm,
using the fact that the special form of the dissipation -- and its isotropy in higher dimensions --
allows to extend the construction of the normal form to the dissipative case,
and hence to deduce that the solution is normally attractive.
Assuming only Hypothesis \ref{hyp:2} on $\om$ does not guarantee $\om_0$ to satisfy any Diophantine condition,
so that the KAM scheme does not apply.
Certainly, the multi-periodic solution is not a global attractor, since the system admits many attracting periodic orbits: essentially
as many as the harmonics of the function $f$ \cite{BC1,CC1,CC3,BDG1}, even though to prove that
the analytical results hold for values of the parameters which fit the astronomical values is not quite trivial \cite{ABC}.
However, not all coexisting attractors are equally relevant \cite{CC3,BDG1}
and, in general, the basins of attraction in dissipative forced systems are strongly entwined
and are apparently fractal \cite{BDG1,WBDG1,WBDG2,WBDG3}, which suggest that
any analytical investigation might be very difficult. 
}
\end{rmk}

\begin{rmk} \label{rmk:1.6}
\emph{
It would be interesting to extend the results to the case of the Lagrangian tori which, in the absence of dissipation,
arise inside the small oscillatory islands which appear where the KAM tori are destroyed by the perturbation \cite{BC2,BC3,MNT}.
Indeed, in the case of Mercury, as already stressed, both numerical and analytical computations for the
the spin-orbit model with realistic tidal torque suggest that the main attractor is a torus of this kind.
}
\end{rmk}

The paper is organised as follows. In Section \ref{sec:2} we introduce an approximate solution to \eqref{eq:1.2},
 up to order $\gotm$ in $\e$,
and write the equation for the correction to the approximate solution.
In Section \ref{sec:3} we first bound the approximate solution, and then
we write the correction as a series expansion and provide a graphical representation for its Fourier coefficients.
In particular, we show that 
the series expansion is well defined to any order, provided suitable counterterms are added to the external parameter.
In Section \ref{sec:4} we show that the series expansion converges.
In Section \ref{sec:5} we show that the counterterm can be
absorbed in the original value of the external parameter, by slightly changing the frequency vector
without modifying the properties of convergence of the series expansion. This concludes the proof of Theorem \ref{thm:1}.
Finally, in Section \ref{sec:6} we show that, by writing $\ga:=\bar\ga |\e|^{-a}$,
the radius of convergence of the series expansion can be proved to be independent of $a$ for all $a\in[0,1)$:
this proves Theorem \ref{thm:2}, as soon as one takes $a=\gotm -s$.

\zerarcounters 
\section{Set-up for the proof of Theorem \ref{thm:1}}
\label{sec:2} 

Consider the non-autonomous ordinary differential equation in $\TTT$
\vspace{-.1cm}
\begin{equation} \label{eq:2.1}
\ddot \theta + \e^\gotm \ga ( \dot \theta - \al) + \e f(\theta,t) = 0 , 
\vspace{-.1cm}
\end{equation}
where $\e\in\RRR$ is a small parameter, $\gotm\in\NNN$ and (cf.~Hypothesis \ref{hyp:1})
%
%
%
%
\begin{equation} \label{eq:2.2}
f(\psi) = \sum_{\substack{\nu \in \ZZZ^2 \\ 0 < |\nu | \le N}} {\rm e}^{i\nu \cdot \psi} f_{\nu} ,
\qquad \psi = (\theta ,t) , \qquad f_{-\nu}^* = f_{\nu} . 
\vspace{-.1cm}
\end{equation}
Since Theorems \ref{thm:1} and \ref{thm:2} hold trivially if $f$ identically zero,
in the following, we assume $\Phi>0$.

\begin{rmk} \label{rmk:2.1}
\emph{
In the case of the spin-orbit model the sum is restricted to integer vectors $\nu\in\ZZZ^2$ of the form $\nu=(2,-k)$,
with $k\neq0$ \cite{GP,MD}. Specifically for the purpose of numerical simulations, the sum is usually truncated by
 requiring $N_1 \le k \le N_2$, for some $N_1,N_2\in \ZZZ$; usually one takes 
$N_1\ge -3$ and $ N_2\le 8$ \cite{CC3,BDG2,NFME,BDG3}. For the Moon and Mercury,
astronomical observations suggest the sizes of the triaxial torque and the tidal torque to be of order $10^{-4}$
and $10^{-8}$, respectively \cite{C,CL,P,CC3}: this means that, if one aims to describe the physical problem
through the differential equation \eqref{eq:2.1}, one has to select $\gotm=2$ in \eqref{eq:2.1}.
}
\end{rmk}

The system described by \eqref{eq:2.1} has a fixed frequency, equal to $1$. We look for a multi-periodic
solution with frequency vector $\om_{0}=(\al_0,1) \in \RRR^2$, with $\al_0$ close to $\al$ to be determined.

Set, for notational convenience,
\vspace{-.1cm}
\begin{equation} \nonumber
\ZZZ^2_* :=\ZZZ^2\setminus\!\{0\}, \qquad \ZZZ^{2}(N) = \{ \nu \in \ZZZ^2: |\nu| \le N\} , \qquad
\ZZZ_*^{2}(N) = \{ \nu \in \ZZZ^2_* : |\nu| \le N\},
\end{equation}
and define
\begin{equation} \label{eq:2.3}
\be_0: = \min\{ |\om_{0}\cdot\nu| : \nu \in \ZZZ^2_{*}(4\gotm N)  \} . 
\end{equation}
%

\begin{rmk} \label{rmk:2.2}
\emph{
If $\om_{0}$ is non-resonant, i.e.~if its components are rationally independent,
one has $\om_{0}\cdot\nu \neq 0$ for all $\nu\in \ZZZ_*^2$ and hence $\be_0 >0$.
A vector $\om_0$ is said to be non-resonant up to order $N$ if $\om_0\cdot\nu \neq 0$ for all $\nu\in \ZZZ_*^2(N)$.
Of course it is sufficient that $\om_0$ be non-resonant up to order $4\gotm N$ for $\be_0$ to be positive:
thus, assuming Hypothesis \ref{hyp:2} implies $\be_0>0$.
Furthermore, if $\om=(\al,1)$ and $|\om\cdot\nu| \ge \be$ for all $\nu\in \ZZZ_*^2(4\gotm N)$,
then one has $\be_0 \ge \be/2$ provided $|\al-\al_0|$ is small enough. 
}
\end{rmk}

For $\e=0$ any solution to \eqref{eq:2.1} is of the form $\theta_0(t) = \bar \theta + \al_0 t$,
where both $\bar\theta\in\TTT$ and $\al_0\in\RRR$ can be arbitrarily chosen; without loss of generality
we can fix $\bar\theta=0$.
Set $\al=\al_0 + \e \tilde\al$, where $\tilde\al$ will be referred to as the \emph{counterterm},
and define the $\gotm$\emph{-th order approximate solution} to \eqref{eq:2.1},
with $f$ given by \eqref{eq:2.2}, as
\begin{equation} \nonumber
\theta_0(t) + h(\om_0 t; \e),\qquad h(\om_0 t;\e):=\e h^{(1)}(\om_{0} t) +\ldots +\e^\gotm h^{(\gotm)}(\om_0 t),
\end{equation}
such that
\begin{equation} \label{eq:2.4}
\ddot h^{(k)} 
+ f^{(k)}(\om_0 t) 
=0, 
\qquad k =1,\ldots,\gotm,
\end{equation}
with $f^{(1)}(\om_0 t)=f(\om_0 t)$ and
\begin{equation} \label{eq:2.5}
f^{(k)}(\om_0 t) := \sum_{p=1}^{k-1} \frac{1}{p!} \frac{\partial^p}{\partial \theta^p} f(\om_0 t)
\!\!\!\!\!\!\!\!\!
\sum_{\substack{k_{1},\ldots,k_{p}\ge 1 \\ k_{1}+\ldots+k_p=k-1}} \!\!\!\!\!\!\!\!
h^{(k_{1})}(\om_0 t) \ldots h^{(k_p)}(\om_0 t) , \qquad k =2,\ldots,\gotm .
\end{equation}
Writing
\begin{equation} \nonumber
f^{(k)}(\psi) = \sum_{\nu \in \ZZZ^2(pN)} {\rm e}^{i\nu \cdot \psi} f^{(k)}_{\nu} , \qquad k=1,\ldots,\gotm,
\end{equation}
one can immediately check by induction that $f^{(k)}_{\nu}$ -- and hence $h^{(k)}_{\nu}$ -- vanishes whenever $|\nu|>kN$.
If $\be_0 \neq 0$, the system \eqref{eq:2.4} can be solved recursively and gives, for all $k=1,\ldots,\gotm$,
\begin{equation} \label{eq:2.6}
h^{(k)}(\psi) = \sum_{\nu \in \ZZZ^2(kN)} {\rm e}^{i\nu \cdot \psi} h^{(k)}_{\nu} , \qquad \qquad
h^{(k)}_{\nu} := \frac{f^{(k)}_{\nu}}{(\om_{0}\cdot\nu)^2} \qquad \forall \nu \in \ZZZ^2_*(kN) ,
\end{equation}
%
%
%
while $f^{(k)}_0=0$ (cf.~Lemma \ref{lem:3.4} below for details) and $h^{(k)}_{0}$ is arbitrary.
The coefficients $h^{(k)}_{\nu}$ are well defined if $\om_0$ is non-resonant up to order $\gotm N$;
hence $h^{(k)}$ is a trigonometric polynomial of degree $k N$, for $k=1,\ldots,\gotm$, and
$h$ is a trigonometric polynomial of degree $\gotm N$.

Once the function $h(\psi;\e)$ has been determined, write
\begin{equation} \label{eq:2.7}
\theta(t) = \al_0 t +  h(\om_{0} t;\e) + \e^{\gotm+1} H(\om_{0} t;\e) , \qquad
H(\psi;\e) = \sum_{\nu \in \ZZZ^2} {\rm e}^{i\nu \cdot \psi} H_{\nu} , 
\end{equation}
which, inserted into \eqref{eq:2.1}, provides a differential equation for $H$:
\begin{equation} \label{eq:2.8}
\ddot H + \e^\gotm \ga \dot H - \ga \tilde \al  +\frac{1}{\e^{\gotm }}  F(\om_0 t, H ;\e)
 = 0 ,
 \vspace{-.1cm}
\end{equation}
with
%
\begin{equation} \label{eq:2.9} 
\begin{aligned}
F(\om_0 t,H;\e) 
:= & \sum_{k=\gotm + 1}^{\io} \e^{k-1} \sum_{p=1}^{k-1} \frac{1}{p!} \frac{\partial^p}{\partial \theta^p} f(\om_0 t)
\!\!\!\!\!\!\!\!\!
\sum_{\substack{1 \le k_{1},\ldots,k_p \le \gotm \\ k_{1}+\ldots+k_p=k-1}} \!\!\!\!\!\!\!\!
h^{(k_{1})}(\om_0 t) \ldots h^{(k_p)}(\om_0 t)
\\ + & 
\sum_{k=1}^{\io} \frac{\partial^k}{\partial \theta^k} f(\om_0 t)
\!\!\!
\sum_{\substack{ p \ge 0, q \ge 1 \\ p+q=k}}
\frac{1}{p!q!} (h(\om_0 t;\e) )^{p} (\e^{\gotm+1} H)^{q} - \e^{\gotm-1} \ga \dot h (\om_0 t;\e). 
\end{aligned}
\end{equation}
Note that, due to the constraints on the exponents, $F$ is at least of order $\e^\gotm$.
Moreover, since $h_0+ \e^{\gotm+1} H_0$ can be absorbed into $\bar\theta$, we can and shall 
assume that both $h_0$ and $H_0$ vanish.

In Fourier space, if we set
\begin{subequations} \label{eq:2.10}
\begin{align}
F^0_{\nu_0} (q) & := \left( i\nu_{0,1} \right)^{q} f_{\nu_0} +
\sum_{p=1}^{\io} \frac{1}{p!} 
\sum_{1 \le k_{1},\ldots,k_p \le \gotm} 
\!\!\!\!\!\ \e^{k_{1}+\ldots+k_p} \!\!\!\!\!\!\!\!\!\!\!
\sum_{\substack{\tilde\nu_0 \in \ZZZ_*^2(N) \\
\tilde\nu_{1},\ldots, \tilde\nu_{p} \in \ZZZ_*^2(\gotm N) \\ \tilde\nu_0+\tilde\nu_1+\ldots+\tilde\nu_p=\nu_0 }} \!\!\!\!\!\!\!\!\!\!
\left( i\tilde\nu_{0,1} \right)^{p+q} f_{\tilde\nu_0} 
h^{(k_{1})}_{\tilde\nu_1} \ldots h^{(k_p)}_{\tilde\nu_p} , 
\label{eq:2.10a} \\
F^1_{\nu_0} 
& :=
\sum_{p=1}^{\io} \frac{1}{p!}
\sum_{\substack{ 1 \le k_{1},\ldots,k_p \le \gotm \\ k_{1}+\ldots+k_p \ge \gotm}} 
\!\!\!\!\!\ \e^{k_{0}+\ldots+k_p - \gotm} \!\!\!\!\!\!\!\!\!\!\!
\sum_{\substack{\tilde\nu_0 \in \ZZZ_*^2(N) \\
\tilde\nu_{1},\ldots, \tilde\nu_{p} \in \ZZZ_*^2(\gotm N) \\
\tilde\nu_0+\tilde\nu_{1}+\ldots+\tilde\nu_{p}=\nu_0 }} \!\!\!\!\!\!
\left( i\tilde\nu_{0,1} \right)^{p} f_{\tilde\nu_0}
h^{(k_{1})}_{\tilde\nu_1} \ldots h^{(k_p)}_{\tilde\nu_p} ,
\label{eq:2.10b} \\
G^{1}_{\nu_0} 
& := - \left( i \om_0\cdot\nu_0 \right) \ga \sum_{k=1}^{\gotm} \e^{k-1} h^{(k)}_{\nu_0} ,
\label{eq:2.10c}
\end{align}
\end{subequations}
%
%
%
then \eqref{eq:2.8} gives for $\nu\neq 0$
\vspace{-.1cm}
\begin{equation} \label{eq:2.11}
\left( \left( \om_{0}\cdot \nu \right)^2   - i \e^\gotm \ga \, \om_{0} \cdot \nu \right) \! H_{\nu}
= F^{1}_{\nu} 
+ G^{1}_{\nu} 
+ \sum_{q=1}^{\io}  \frac{1}{q!} 
\e^{(\gotm+1) q-\gotm} 
\!\!\!\!\!\!\!\!\!\! \sum_{\substack{ \nu_0,\nu_1,\ldots,\nu_q \in \ZZZ^2_* \\ \nu_0+\ldots+\nu_q=\nu }} \!\!\!\!\! \!\!\!\!\!
F^{0}_{\nu_0}(q) \, H_{\nu_1} \ldots H_{\nu_q} ,
\vspace{-.2cm}
\end{equation}
provided one has for $\nu=0$
\vspace{-.1cm}
\begin{equation} \label{eq:2.12}
\gamma \tilde \al =  F_{0}^{1} 
+ \sum_{q=1}^{\io} \frac{1}{q!} 
\e^{(\gotm +1)q-\gotm } 
\!\!\!\!\!  \sum_{\substack{ \nu_0,\nu_1,\ldots,\nu_q \in \ZZZ^2_* \\ \nu_0+\ldots+\nu_q=0}} \!\!\!\!\! 
F_{\nu_0}^{0}(q) \, H_{\nu_1} \ldots H_{\nu_q}  .
\vspace{-.2cm}
\end{equation}

In the following sections we prove that \eqref{eq:2.11} admits a solution in the form of a convergent series,
provided $\tilde\al$ is fixed according to \eqref{eq:2.12}.
For given $\al_0$, the counterterm $\tilde\al$ is uniquely defined by \eqref{eq:2.12}: this fixes $\al$
to the value $\al=\al_0+\e\tilde\al$.
However, the physical parameter is $\al$, so what we really need to do is the opposite,
i.e.~to find, for fixed $\al$, a value $\al_0$ such that $\al_0+\e\tilde\al=\al$.
This leads to solve an implicit function problem.

\zerarcounters 
\section{Construction of the perturbation series}
\label{sec:3} 

\subsection{Trees}
\label{sec:3.1} 

A \textit{rooted tree} $\vartheta$ is an oriented graph with no cycle,
such that all the lines are oriented toward a unique
vertex (the \textit{root}) which has only one incident line (the \emph{root line} $\ell_\vartheta$).
We call \textit{nodes} all the vertices in $\vartheta$ except the root.
The orientation of the lines in $\vartheta$ induces a partial ordering 
relation ($\preceq$) between the nodes. Given two nodes $v$ and $w$,
we shall write $w \prec v$ every time $v$ is along the path
(of lines) connecting $w$ to the root; we shall write $w\prec \ell$ if
$w\preceq v$, where  $v$ is the unique node that the line $\ell$ exits. If a line $\ell$
exits a node $v$ we may write $\ell=\ell_v$.
For any node $v$ denote by $q_{v}$ the number of lines entering $v$.

Given a rooted tree $\vartheta$ we denote by $N(\vartheta)$ the set of nodes,
by $E(\vartheta)$ the set of \textit{end nodes} (or \emph{levaes}), i.e.~nodes $v$ with $q_{v}=0$,
by $I(\vartheta)$ the set of \textit{internal nodes}, i.e.~nodes $v$ with $q_{v}\ge 1$,
and by $L(\vartheta)$ the set of lines; by definition $N(\vartheta)=E(\vartheta) \sqcup I(\vartheta)$.

We associate with each node $v\in N(\vartheta)$ a \textit{mode} label $\nu_v\in\ZZZ^2$ and with each line
$\ell\in L(\theta)\setminus\{\ell_\vartheta\}$ a \textit{momentum} label $\nu_\ell\in\ZZZ^2_*$, while
we associate with the root line $\ell_\vartheta$ a momentum $\nu_{\ell_\vartheta}\in\ZZZ^2$.
We also impose the \emph{conservation law}
\vspace{-.1cm}
\begin{equation} \nonumber 
\nu_\ell = \sum_{v\prec \ell} \nu_v \qquad \forall \ell \in L(\vartheta) . 
\vspace{-.2cm}
\end{equation}
%

\begin{lemma} \label{lem:3.1}
For any tree $\theta$ with $k$ nodes one has
\vspace{-.1cm}
\begin{equation} \nonumber
\sum_{v\in I(\vartheta)} \!\!\! q_v = k-1 , \qquad |I(\theta)| \le k-1 .
\vspace{-.1cm}
\end{equation}
\end{lemma}

\noindent\emph{Proof}. 
For each node $v$ there are $q_v$ entering lines and all lines enter a node, the only exception being the root line. 
Moreover there is at least one end node. \qed

\begin{lemma} \label{lem:3.2}
The number of unlabelled trees with $k$ nodes is bounded by $4^k$.
\end{lemma}

\noindent\emph{Proof}. 
The number of trees with $k$ nodes and with no labels is bounded by the number of random walks with $2k$ steps,
that is $2^{2k}$. \qed

\subsection{Tree representation for the function $\boldsymbol{h}$}
\label{sec:3.2} 

Here and in the following, we provide a graphical representation in terms of labelled trees,
as defined in Section \ref{sec:3.1}, for both the coefficients $h^{(k)}_{\nu}$ and $H^{(k)}_{\nu}$.
We start with the the coefficients of the approximate solution $h$.

Assume $\om_0$ to be non-resonant up to order $\gotm N$, so that $\be_0 \neq 0$.
In Fourier space \eqref{eq:2.5} gives, for $k=2,\ldots,\gotm$,
\begin{equation} \label{eq:3.1} 
f^{(k)}_{\nu} = 
\sum_{p=1}^{k-1} \frac{1}{p!}
\sum_{\substack{k_{0},\ldots,k_p \ge 1 \\ k_{0}+\ldots+k_p=k-1}}
\sum_{\substack{\nu_0 \in \ZZZ_*^2(N)}} 
(i\nu_{0,1})^p f_{\nu_0} \!\!\!\!\!\!\!\!
\sum_{\substack{\nu_1,\ldots,\nu_p  \in \ZZZ_*^2(\gotm N) \\ \nu_1+\ldots+\nu_p=\nu-\nu_0}} \!\!\!\!\!\!\!\!
h^{(k_{0})}_{\nu_1} \ldots h^{(k_p)}_{\nu_k} .
\end{equation}

A tree $\vartheta$ is said to be of order $k$ if one has $|N(\vartheta)|=k$.
Given a tree $\vartheta$, one associates
with each node $v\in N(\vartheta)$ a \textit{node factor} $\calI^{*}_v$
and with each line $\ell\in L(\vartheta)$ a \textit{propagator} $\calG^{*}_\ell$, where
$\calI^{*}_v$ and $\calG^{*}_\ell$ are defined as
%
\begin{equation} \nonumber
\calI^{*}_v := \frac{(i\nu_{v,1})^{q_v}}{q_v!} f_{\nu_v} , \qquad
\calG^{*}_\ell:= \begin{cases}
\displaystyle{ \frac{1}{(\om_{0}\cdot\nu_{\ell})^2} }, & \nu_{\ell} \in \ZZZ_*^2(\gotm N) , \\
\vspace{-.4cm} \\
1 , & \nu_{\ell} = 0 .
\end{cases}
\vspace{-.1cm}
\end{equation}
Define the \emph{value} of the tree $\vartheta$ as
\vspace{-.1cm}
\begin{equation} \nonumber 
\Val^{*}(\vartheta) := 
\Bigl( \prod_{v \in N(\vartheta)} \calI^{*}_v \Bigr)
\Bigl( \prod_{\ell \in L(\vartheta)} \calG^{*}_{\ell} \Bigr) ,
\end{equation}
and let $\gotT_{k,\nu}^{*}$ denote the set of the labelled trees of order $k$ and with
momentum $\nu$ associated with the root line.

Finally, set, for $\nu\in\ZZZ^2$ and $k\in\ZZZ_+$,
\begin{equation} \label{eq:3.2}
p(\nu) := \min \left\{ {p\in\ZZZ_+} : |\nu| \le (p+1) N \right\} , \qquad
p_{k}(\nu) = \max\{ p(\nu) - k , 0 \} ,
\qquad
\end{equation}
and define the \emph{small divisors}
\begin{equation} \label{eq:3.3}
D(\om_{0}\cdot\nu) := \left( \om_{0}\cdot \nu \right)^2  - i \e^\gotm \ga \, \om_{0} \cdot \nu - \e F^{0}_{0}(1) .
\end{equation}

\begin{lemma} \label{lem:3.3}
Let $p(\nu)$ and $p_k(\nu)$ be defined as in \eqref{eq:3.2}. Then one has
\begin{subequations} \label{eq:3.4}
\begin{align}
& p(\nu_1) + p(\nu_2) \ge p(\nu_1+\nu_2) - 1 ,
\label{eq:3.4a} \\
& p(\nu_1) + p_{k}(\nu_2) \ge p_{k}(\nu_1+\nu_2) - 1 .
\label{eq:3.4b}
\end{align} 
\end{subequations}
\end{lemma}

\noindent\emph{Proof}. 
One has
\[
 p(\nu_1+\nu_2) N < |\nu_1+\nu_2| \le |\nu_1| + |\nu_2| \le (p(\nu_1) +1) N + (p(\nu_2) + 1) N , 
\]
which yields $p(\nu_1+\nu_2) + 1 \le p(\nu_1)+p(\nu_2) + 2$,
since $p(\nu)$ is an integer for any $\nu\in\ZZZ^2$.  This proves \eqref{eq:3.4a}. 

If $p(\nu_2) \ge k$ and $p(\nu_1+\nu_2) \ge k$, then \eqref{eq:3.4b} follows from \eqref{eq:3.4a}.
If $p(\nu_1+\nu_2) < k$, 
\eqref{eq:3.4b} holds trivially, since
both $p(\nu_1)$ and $p_{k}(\nu_2)$ are non-negative. Finally, if $p(\nu_2)<k$ and $p(\nu_1+\nu_2) \ge k$,
one has $|\nu_1+\nu_2| \le |\nu_1| + |\nu_2| \le (p(\nu_1)+1) N + kN$ and hence
$p(\nu_1+\nu_2) \le p(\nu_1)+k+1$, which yields $p(\nu_1) \ge p_k(\nu_1+\nu_2)-1$. \qed

\begin{lemma} \label{lem:3.4}
Assume Hypothesis \ref{hyp:1}, and
assume $\om_0$ to be non-resonant up to order $\gotm N$.
For $k=1,\ldots,\gotm$, one has
\vspace{-.3cm}
\begin{equation} \nonumber 
f^{(k)}_{0} := \sum_{\vartheta \in \gotT_{k,0}^{*}} \Val^{*}(\vartheta) ,
\qquad\quad
h^{(k)}_{\nu} := \begin{cases}
\displaystyle{ \sum_{\vartheta \in \gotT_{k,\nu}^{*}} \Val^{*}(\vartheta)} , & \nn \in \ZZZ^2_*(kN) , \\
\vspace{-.3cm} \\
0 , & |\nu| > kN . \end{cases}
\end{equation}
Moreover, one has $f^{(k)}_{0} = 0$ for $k=1,\ldots,\gotm$ and hence the function \ref{eq:2.6} solves \eqref{eq:2.4}.
\end{lemma}

\noindent\emph{Proof}. 
The graphical representation in tems of trees is proved by induction,
using the recursive definition \eqref{eq:2.6} for the coefficients $h^{(k)}_{\nu}$ and the definition of the tree values. 

The vanishing of $f^{(k)}_{0}$ is a symmetry property, due to the fact that,
up to order $\gotm$ the system is conservative \cite{G,CF,GM}.
Given a tree $\vartheta$, if $v_0$ is the node the root line exits, call
$\vartheta(v)$ the tree obtained by detaching the root line from $v_0$ and reattaching to the node $v\in N(\vartheta)$,
and $\calmP(v_0,v)$ the unique path of lines connecting the node $v$ to the node $v_0$.
f $\FFF_0(\vartheta)$ denotes the family of all trees obtained from $\vartheta$ by varying the node $v$
which the root line is reattached to, one has
\begin{equation} \nonumber
\sum_{\vartheta \in \gotT_{k,0}} \Val^{*}(\vartheta) = \sum_{\vartheta \in \gotT_{k,0}}
\frac{1}{|\FFF_0(\vartheta)|} \sum_{\vartheta' \in \FFF_0(\vartheta)} \Val^{*}(\vartheta') .
\vspace{-.1cm}
\end{equation}
Since the mean of $f$ is zero, one can express $f(\theta,t)$ as the derivative with respect to $\theta$ of a function $\f(\theta,t)$
and write $f_{\nu}=i\nu_{1}\f_{\nu}$. Let $Q_w$ denote the number of lines incident on $w \in N(\vartheta)$,
without counting the root line for the node $v_0$;
for any $\vartheta(v)\in\FFF_0(\vartheta)$ one has $q_v=Q_v+1$ and $q_w=Q_w$ for any other node $w$.
Call $\vartheta_*$ the non-ordered graph obtained from $\vartheta$ by detaching the root line, and set
\begin{equation} \nonumber
\Val^{*}(\vartheta_*) =
\Bigl( \prod_{v \in N(\vartheta_*)} \frac{(i\nu_{v,1})^{Q_v} f_{\nu_v}}{Q_v!} \Bigr)
\Bigl( \prod_{\ell \in L(\vartheta_*)} \calG_{\ell} \Bigr) ,
\end{equation}
with $N(\vartheta_*):=N(\vartheta)$ and $L(\vartheta_*):=L(\vartheta)\setminus\{\ell_\theta\}$.
The value of tree $\theta(v)$ obtained from $\theta$ by reattaching the root line to $v$
differs from $\Val^{*}(\theta)$ because the momenta of the lines $\ell\in \calmP(v_0,v)$ are reversed
(i.e.~$\nu_{\ell}$ is replaced by $-\nu_{\ell}$) and the combinatorial factors $1/q_{v_0}!$ and $1/q_{v}!$
of the two nodes $v_0$ and $v$ are replaced by $1/(q_{v_0}-1)!$ and $1/(q_{v}+1)!$, respectively.
However, there are $q_{v_0}$ trees obtained from $\theta_*$ which have the same value as $\theta$
and $q_v+1$  trees obtained from $\theta_*$ trees which have the same value as $\theta(v)$,
so that 
\begin{equation} \nonumber
\sum_{\vartheta' \in \calI_0(\vartheta)} \Val^{*}(\vartheta') = 
\sum_{v\in N(\vartheta)} i\nu_{v,1} \, \Val^{*}(\vartheta_*) =
i \Val^{*}(\vartheta_*) \sum_{v\in N(\vartheta)} \nu_{v,1} = 0 ,
\end{equation}
which implies that $f^{(k)}_0=0$ for $k=1, \ldots,\gotm$.
\qed

\begin{lemma} \label{lem:3.5}
Assume Hypothesis \ref{hyp:1},
and assume $\om_0$ to be non-resonant up to order $\gotm N \!$.
For $k=1,\ldots,\gotm$, one has
\vspace{-.1cm}
\begin{equation} \nonumber
|h^{(k)}_{\nu}| \le (\be_0^{-2} \Phi_0)^{k} \quad \forall \nu \in \ZZZ^2_*(kN) ,
\qquad \qquad \Phi_0 := 16 N^2 \Phi ,
\vspace{-.1cm}
\end{equation}
with $\be_0$ as in \eqref{eq:2.3}.
Moreover one has $h^{(k)}_{\nu}=0$ for all $\nu\in\ZZZ_2^2$ such that $p(\nu) \ge k$.
\end{lemma}

\noindent\emph{Proof}. 
The value of any tree $\vartheta \in \gotT_{k,\nu}^{*}$, with $k \le \gotm$ and $\nu\neq 0$,
is bounded as
%
\begin{equation} \nonumber
|\Val^{*}(\vartheta)| \le N^{k-1} \Phi^k \be_0^{-2k} ,
\end{equation}
where Lemma \ref{lem:3.1} has been used.
The sum of the mode labels is bounded by $(4N^2)^{k}$ and the number of trees of order $k$
is bounded by $4^k$ by Lemma \ref{lem:3.2}. This yields the bound on $h^{(k)}_{\nu}$. 
Finally, by construction (cf.~Lemma \ref{lem:3.4}), $h^{(k)}_{\nu} \neq 0$ requires $|\nu| \le k N$, i.e.~$p(\nu) \le k - 1$. \qed
%
%

\subsection{Bounds on the the small divisors}
\label{sec:3.3} 

\begin{lemma} \label{lem:3.6}
Assume Hypothesis \ref{hyp:1},
and assume $\om_0$ to be non-resonant up to order $\gotm N \!$.
There exist positive constants $A_0$, $K_{0}$ and $\e_1$ such that, for all $|\e|\le \e_1$, one has
\begin{equation} \nonumber
\left| F_0^{0}(1) + A_0 \e \right| \le K_{0} |\e|^2 . 
\end{equation} 
\end{lemma}

\noindent\emph{Proof}. 
From \eqref{eq:2.10a}, with $q=1$ and $\nu_0=0$, one obtains
\begin{equation} \nonumber
F_{0}^{0} (1) = \e \!\!\!\!\! \sum_{{\tilde\nu_0 \in \ZZZ_*^2(N)  }} \!\!\!
\left( i\tilde\nu_{0,1} \right)^{2} f_{\tilde\nu_0} h^{(1)}_{-\tilde\nu_0}
+\tilde{F}_0^0(1)
\vspace{-.2cm}
\end{equation}
with
\vspace{-.2cm}
\begin{eqnarray}\nonumber
\tilde{F}_0^0(1)
& \!\!\!\! := \!\!\!\! & 
\sum_{{k=2 }}^\gotm \e^k \!\!\!
\sum_{{\tilde\nu_0 \in \ZZZ_*^2(N)  }}  
\left( i\tilde\nu_{0,1} \right)^{2} f_{\tilde\nu_0}  h^{(k)}_{-\tilde\nu_{0}} \\
& \!\!\!\! + \!\!\!\! & 
\sum_{p=2}^{\io} \frac{1}{p!}
\sum_{1 \le k_{1},\ldots,k_p \le \gotm} \!\!\!\!\!\ \e^{k_{1}+\ldots+k_p} \!\!\!\!\!\!\!\!\!\!\!
\sum_{\substack{\tilde\nu_0,\tilde\nu_{1},\ldots, \tilde\nu_{p} \in \ZZZ_*^2(\gotm N) \\ 
\tilde\nu_0+\tilde\nu_1+\ldots+\tilde\nu_p=\nu_0 }} \!\!\!\!\! 
\left( i\tilde\nu_{0,1} \right)^{p+1} f_{\tilde\nu_0} h^{(k_{1})}_{\tilde\nu_1} \ldots h^{(k_p)}_{\tilde\nu_p} ,
\nonumber
\end{eqnarray}
where we have (a) used that the first term in \eqref{eq:2.10a} vanishes and (b) separated
the rest of the sum, i.e.~$\tilde{F}_0^0(1)$, from the contribution of order $\e$.
The latter, using \eqref{eq:2.6}, becomes
\begin{equation} \nonumber
\e \!\!\!\!\! \sum_{\tilde\nu_0 \in \ZZZ_*^2(N)} \!\!\!\!\! 
( i\tilde\nu_{0,1})^{2} \frac{ f_{\tilde\nu_0} 
f_{-\tilde\nu_0}}{(\om_{0}\cdot\tilde\nu_0)^2} = -2  \e \, C_0 , 
\end{equation}
%
%
with
\begin{equation} \label{eq:3.5}
0 < C_0 := \frac{1}{2} \sum_{\tilde\nu_0 \in \ZZZ_*^2(N)} \!\!\!\!\! 
\left( \tilde\nu_{0,1} \right)^{2} \frac{|f_{\tilde\nu_0}|^2}{(\om_{0}\cdot\tilde\nu_0)^2} \le 2N^4\be_0^{-2}\Phi^{2} .
\end{equation}
%
%
%
%
%
Moreover one has, by using Lemma \ref{lem:3.5},
\begin{equation} \nonumber
\begin{aligned}
| \tilde{F}^0_0(1) | & \le
\sum_{k=2}^{\gotm} 4 N^{4} \Phi \left( |\e| \be_0^{-2}\Phi_0 \right)^{k} +
\sum_{p=2}^{\io} 
\sum_{k=p}^{p\gotm N} 
\sum_{\substack{ 1 \le k_{0},\ldots,k_p \le \gotm \\ k_{0}+\ldots+k_p = k}} 
\!\!\!\!\!\ N^{p+1} \Phi (4 \gotm^2 N^2)^{p} ( |\e |\be_0^{-2} \Phi_0)^k \\
& \le
4 \gotm N^{4} \Phi \left( |\e| \be_0^{-2}\Phi_0 \right)^{2} +
3 N \Phi \left( 8 \gotm^2 N^3 |\e| \be_0^{-2}\Phi_0 \right)^{2} \le K_{0} |\e|^2 , 
\end{aligned}
\end{equation}
 with
 \begin{equation} \label{eq:3.6}
 K_{0} := 4 N \Phi \Bigl( \frac{8 \gotm^2 N^3 \Phi_0}{\be_0^{2}} \Bigr)^2 ,
 \end{equation}
provided one takes $|\e|\le \e_1$, with
\begin{equation} \label{eq:3.7}
\e_1 := \frac{\be_0^2}{16 \gotm^2 N^3  \Phi_0} .
\end{equation}
Therefore, one obtains
$\left| F_0^{0}(1) + A_0 \e \right| \le K_{0} |\e|^2 $ for all $|\e| \le \e_1$. \qed

\begin{lemma} \label{lem:3.7}
Assume Hypothesis \ref{hyp:1}, 
and assume $\om_0$ to be non-resonant up to order $\gotm N \!$.
There exists positive constants $C_0$ and $\e_2$ such that, for all $|\e| \le \e_2$, one has
\begin{eqnarray} \nonumber
|D(\om_{0}\cdot\nu)| & \!\!\!\! \ge \!\!\!\! &
|\e|^{\gotm}\ga \, | \om_0\cdot\nu| , \qquad \forall \nu\in \ZZZ_*^2  , \\ [0.6ex]
|D(\om_{0}\cdot\nu)| & \!\!\!\! \ge \!\!\!\! &
(\om_{0}\cdot\nu)^2 +  C_0|\e|^2 \ge \max \left\{ C_0 |\e|^2 , \frac{\be_0^{2}}{16} \right\} \qquad \forall \nu\in \ZZZ_*^2 \nonumber \\ [0.6ex]
|D(\om_{0}\cdot\nu)| & \!\!\!\! \ge \!\!\!\! &
 \be_0^2 \qquad \forall \nu\in \ZZZ_*^2(4\gotm N) . \nonumber
\end{eqnarray}
%
%
%
\end{lemma}

\noindent\emph{Proof}. 
Define
\begin{equation} \label{eq:3.8}
\tilde \e_2 := \frac{C_0}{K_{0}} = \frac{C_0\,\e_1^2}{N\Phi} ,  \qquad \qquad
\e_2 := \min \{ \e_1 , \tilde \e_2 \} , 
\end{equation}
with $C_0$, $K_{0}$ and $\e_1$ defined in \eqref{eq:3.5}, in \eqref{eq:3.6} and in \eqref{eq:3.7}, respectively.

One easily check from \eqref{eq:2.10a} that $F^0_0(q)$ is real for any $q\in\ZZZ_+$, so that one can bound
$|D(\om_{0}\cdot\nu)| \ge |\e^{\gotm}\ga \, \om_0\cdot\nu|$. 

Moreover one has
\begin{equation} \nonumber
|D(\om_{0}\cdot\nu)| \ge |(\om_{0}\cdot\nu)^2 - \e F^0_0(1)| \ge 
(\om_{0}\cdot\nu)^2 +  C_0|\e|^2 \ge 
\max \left\{ (\om_{0}\cdot\nu)^2 , C_0 |\e|^2  \right\} ,
\end{equation}
for all $|\e| \le \e_2$.
Thus, one can also bound $D(\om_{0}\cdot\nu)$ by $C_0|\e|^2$,  whenever one has
$|\om\cdot\nu| < \be_0/4$.

Finally, one has $|D(\om_{0}\cdot\nu)| \ge (\om_{0}\cdot\nu)^2 \ge \be_0^2$
for all $\nu\in\ZZZ_*^2(4\gotm N)$ by \eqref{eq:2.3}. \qed

\begin{rmk} \label{rmk:3.8}
\emph{
By comparing \eqref{eq:3.8} with \eqref{eq:3.6} and \eqref{eq:3.7},
one finds
\begin{equation} \nonumber
\tilde\e_2 \le \frac{C_0 \e_1^2}{N\Phi} , \qquad 
\frac{C_0\e_1}{N\Phi} \le \frac{4N^4 \be_0^{-2}\Phi^2}{2N\Phi} \frac{\be_0^2}{16 \gotm^2 N^3 (16 N^2 \Phi)} \le
\frac{2}{(16 \gotm N)^2} ,
\end{equation}
so that $\tilde \e_2 \le 2 (16 \gotm N)^{-2} \e_1$.
Thus, one has $\e_2=\tilde\e_2$.
}
\end{rmk}

\subsection{Bounds on the the coefficients
$\boldsymbol{F_{\nu_0}^0(q)}$, $\boldsymbol{F_{\nu_0}^1 
}$ and $\boldsymbol{G^{1}_{\nu_0} 
}$.}
\label{sec:3.4} 

\begin{lemma} \label{lem:3.9}
Assume Hypothesis \ref{hyp:1},
and assume $\om_0$ to be non-resonant up to order $\gotm N \!$.
Let $\e_1$ be as in Lemma \ref{lem:3.6}.
For all $|\e| \le \e_1$, all $q\ge 1$ and all $\nu_0 \in\ZZZ^2$, one has
\begin{equation} \nonumber
|F_{\nu_0}^0(q)| \le 
3 N^q \Phi |\gotD_0\e|^{p(\nu_0)} , \qquad
\gotD_0 := 8 \gotm^2 N^3 \be_0^{-2} \Phi_0 ,
\end{equation}
where $\Phi_0$ is as in Lemma \ref{lem:3.5}.
\end{lemma}

\noindent\emph{Proof}. 
According to \eqref{eq:2.10a} and Lemma \ref{lem:3.5}, we have, if $p(\nu_0) \ge 1$ and hence $f_{\nu_0}=0$, 
\begin{eqnarray} \nonumber
|F_{\nu_0}^{0}(q)|
& \!\!\!\! \le \!\!\!\! & 
\sum_{p=1}^{\io} \frac{1}{p!} N^{p+q} \, \Phi 
\!\!\!\!\!
\sum_{\substack{\tilde\nu_1,\ldots,\tilde \nu_p \in \ZZZ_*^2(\gotm N)}} 
\Biggl( \, \prod_{i=1}^{p} \sum_{k_i \ge p(\tilde\nu_i) + 1} \left( |\e| \be_0^{-2} \Phi_0 \right)^{k_i} \Biggr) \\
& \!\!\!\! \le \!\!\!\! & 
N^{q} \Phi
\sum_{p=1}^{\io} N^p \left( 4 \gotm^2 N^2 \right)^{p} 2^p \left( |\e| \be_0^{-2} \Phi_0 \right)^{\max\{p(\nu_0) , p\}} ,
\nonumber
\end{eqnarray}
where we have used that $p(\tilde\nu_0)=0$, $p(\tilde\nu_1),\ldots,p(\tilde\nu_p) \ge0$ and
$p(\tilde\nu_0)+p(\tilde\nu_1)+\ldots+p(\tilde\nu_p) \ge p(\nu_0) -p$ by Lemma \ref{lem:3.3}, so that we obtain
\begin{eqnarray} \nonumber
|F_{\nu_0}^{0}(q)|
& \!\!\!\! \le \!\!\!\! & 
N^q \Phi  \Biggl(
\left( |\e| \be_0^{-2} \Phi_0 \right)^{p(\nu_0)}
\sum_{p=1}^{p(\nu_0)} \left( 8 \gotm^2 N^3 \right)^{p} + \!\!\!\!
\sum_{p=p(\nu_0)+1}^{\io} \left( 8 \gotm^2 N^3 |\e| \be^{-2} \Phi_0 \right)^{p} \Biggr) \\
& \!\!\!\! \le \!\!\!\! & 
N^q \Phi
\Bigl(
2 \left( 8 \gotm^2  N^3 |\e| \be_0^{-2} \Phi_0 \right)^{p(\nu_0)} +
\left( 8 \gotm^2 N^3 |\e| \be_0^{-2} \Phi_0 \right)^{p(\nu_0)} \Bigr) , \nonumber
\end{eqnarray}
provided we take $|\e| \le \e_1$, with $\e_1$ given by \eqref{eq:3.7}.

If $p(\nu_0)=0$, instead, we have
\begin{equation} \nonumber
|F_{\nu_0}^{0}(q)| \le 
N^q \Phi +
N^{q} \Phi
\sum_{p=1}^{\io} N^p \left( 4 \gotm^2 N^2 \right)^{p} 2^p \left( |\e| \be_0^{-2} \Phi_0 \right)^{p} \\
\le
2 N^q \Phi ,
\end{equation}
because we have to take into account also the first summand in \eqref{eq:2.10a}.
\qed

\begin{rmk} \label{rmk:3.10}
\emph{
For $\nu_0=0$, the bound in Lemma \ref{lem:3.9} might be improved into
$2 N^q \Phi |\gotD_0\e|$, since $p(\nu_0)=0$ and $f_{0}=0$ vanishes.
However, in the following, we do not need the improved bound, except when $q=1$,
for which Lemma \ref{lem:3.6} applies and, in particular, gives $| F^0_0(1) | \le 2A_0|\e|$,
if one assumes the stronger condition $|\e|\le \e_2$.
}
\end{rmk}

\begin{lemma} \label{lem:3.11}
Assume Hypothesis \ref{hyp:1},
and assume $\om_0$ to be non-resonant up to order $\gotm N \!$.
Define $\e_1$ and $\gotD_0$ as in Lemmas \ref{lem:3.6} and \ref{lem:3.9}, respectively,
and set $\Gamma:=2 \ga \be_0^{-1} \Phi_0\gotD_0^{-\gotm}$.
Then, for all $|\e| \!\le\! \e_1$ one has
\begin{eqnarray} \nonumber
|F_{\nu_0}^1
| 
& \!\!\! \le \!\!\! &
6 \Phi 
\gotD_0^{\gotm} |\gotD_0\e|^{p_{\gotm} (\nu_0)}  , 
\qquad \nu_0 \in
\ZZZ^2_* , \nonumber \\ [1ex]
|F_{0}^1 
| 
& \!\!\! \le \!\!\! &
6 \Phi \gotD_0^{\gotm} |\gotD_0\e| ,
\nonumber \\ [1ex]
|G^{1}_{\nu_0} 
| 
& \!\!\! \le \!\!\! & \Gamma \gotD_0^{\gotm} ,
\qquad \hspace{1.32cm}
\nu_0 \in \ZZZ^2_*(\gotm N) ,
\nonumber 
\end{eqnarray}
while one has $G^{1}_{\nu_0}
=0$ for $\nu_0 \in \ZZZ^2 \setminus \ZZZ^2_*(\gotm N)$.
\end{lemma}

\noindent\emph{Proof}. 
For $|\e| \le  \e_1$, with $\e_1$ as in Lemma \ref{lem:3.6},
according to \eqref{eq:2.10b} one has
\begin{eqnarray} \nonumber
|F_{\nu_0}^{1} 
|
& \!\!\!\! \le \!\!\!\! & 
\sum_{p=1}^{\io} \frac{1}{p!} N^{p} \Phi 
\!\!\!\!\! 
\sum_{k \ge \max\{\gotm , p \}} 
\sum_{\substack{k_{0},\ldots,k_p \ge 1 \\ k_{0} + \ldots + k_p = k}}  (4\gotm^2 N^2)^{p} 
|\e|^{k-\gotm} \left( \be_0^{-2} \Phi_0 \right)^{k} \\ [-1.6ex]
& \!\!\!\! \le \!\!\!\! & 
\Phi |\e|^{-\gotm}
\sum_{p=1}^{\gotm} \sum_{k=\gotm}^{\io}
(4\gotm^2 N^3)^{p} \left( 2 |\e| \be_0^{-2} \Phi_0 \right)^{k}
+ \Phi |\e|^{-\gotm}
\sum_{p=\gotm+1}^{\io} \sum_{k=p}^{\io}
(4\gotm^2 N^3)^{p} \left( 2 |\e| \be_0^{-2} \Phi_0 \right)^{k}
\nonumber \\ [-1.0ex]
& \!\!\!\! \le \!\!\!\! & 
2 \Phi |\e|^{-\gotm} \sum_{k=\gotm}^{\io}
(4\gotm^2 N^3)^{\gotm} \left( 2 |\e| \be_0^{-2} \Phi_0 \right)^{k}
+ 2 \Phi |\e|^{-\gotm}
\sum_{p=\gotm+1}^{\io} (4\gotm^2 N^3)^{p} \left( 2 |\e| \be_0^{-2} \Phi_0 \right)^{p} \nonumber \\
& \!\!\!\! \le \!\!\!\! & 
6 \Phi |\e|^{-\gotm} \left( 8 \gotm^2 N^3 |\e| \be_0^{-2} \Phi_0 \right)^{\gotm} , \nonumber
\end{eqnarray}
which yields the first bound for $\nu_0\in\ZZZ_*^2(\gotm N).$
%
%
Since the sum over the Fourier labels in \eqref{eq:2.10b} has the constraint $\nu_0=\tilde\nu_0+\tilde\nu_1+\ldots+\tilde\nu_p$,
by Lemma \ref{lem:3.4} one has
\begin{equation} \nonumber
|\nu_0| \le \sum_{p=0}^{k} |\tilde \nu_k| \le N + k_{0} N + \ldots + k_p N \le (k+1)N ,
\end{equation}
so that $k \ge p(\nu_0)$. If $|\nu_0| > \gotm N$, then $p(\nu_0) \ge \gotm$ and hence
\begin{eqnarray} \nonumber
|F_{\nu_0}^{1} 
| 
& \!\!\!\! \le \!\!\!\! & 
\sum_{p=1}^{\io} \frac{1}{p!} N^{p} \Phi \sum_{k \ge \max\{p(\nu_0) , p \}} 
\sum_{\substack{k_{0},\ldots,k_p \ge 1 \\ k_{0} + \ldots + k_p = k}}  (4\gotm^2 N^2)^{p}
|\e|^{k-\gotm} \left( \be_0^{-2} \Phi_0 \right)^{k} , \\ [-1.6ex]
& \!\!\!\! \le \!\!\!\! & 
2 \Phi |\e|^{-\gotm} \!\!\!\!
\sum_{k=p(\nu_0)}^{\io}
(4\gotm^2 N^3)^{p(\nu_0)} \left( 2 |\e| \be_0^{-2} \Phi_0 \right)^{k}
+ 2 \Phi |\e|^{-\gotm} \!\!\!\!
\sum_{p=p(\nu_0)+1}^{\io} (4\gotm^2 N^3)^{p} \left( 2 |\e| \be_0^{-2} \Phi_0 \right)^{p} \nonumber \\
& \!\!\!\! \le \!\!\!\! & 
6 \Phi |\e|^{-\gotm} \left( 8\gotm^2 N^2 |\e| \be_0^{-2} \Phi_0 \right)^{p(\nu_0)} , \nonumber
\nonumber
\end{eqnarray}
which shows that the bound still holds for $\nu_0 \in \ZZZ^2 \setminus \ZZZ^2(\gotm N)$.

For $\nu_0=0$, the contribution with $k_{0}+\dots+k_p=\gotm$ in \eqref{eq:2.10b} is
\begin{equation} \nonumber
f^{(\gotm+1)}_0 :=
\sum_{p=1}^{\io} \frac{1}{p!}
\sum_{\substack{ 1 \le k_{0},\ldots,k_p \le \gotm \\ k_{0}+\ldots+k_p =\gotm}} 
\sum_{\substack{\tilde\nu_0 \in \ZZZ_*^2(N) \\
\tilde\nu_{1},\ldots, \tilde\nu_{p} \in \ZZZ_*^2(\gotm N) \\
\tilde\nu_0+\tilde\nu_{1}+\ldots+\tilde\nu_{p}=0 }} \!\!\!\!\!\!
\left( i\tilde\nu_{0,1} \right)^{p} f_{\tilde\nu_0}
h^{(k_{0})}_{\tilde\nu_1} \ldots h^{(k_p)}_{\tilde\nu_p} ,
\end{equation}
which vanishes (the proof is the same as the proof that $f^{(k)}_0=0$ for $k=1,\ldots,\gotm$ in Lemma \ref{lem:3.5}),
while the contributions with $k_{0}+\dots+k_p \ge \gotm+1$ can be bounded by the same reasoning as was used
to obtain the first bound,
with the difference that now $k\ge\max\{ p, \gotm+1\}$. By collecting the two results one obtains the second bound.

The last bound follows from \eqref{eq:2.10c},
by noting that, for $|\e|\le \e_1$, one has, for $\nu_0\in\ZZZ_*^2(\gotm N)$,
\begin{equation} \nonumber
|G^{1}_{\nu_0} 
| \le \be_0^{-1} |\e|^{-1} \ga \sum_{k=1}^{\gotm} \be_0^{-2(k-1)} |\e|^k \Phi_0^{k} \le 
2 \ga \be_{0}^{-1} \Phi_0^2 ,
\end{equation}
since one has $|h^{(k)}_{\nu} | \le |\om_0\cdot\nu|^{-2} \be_0^{-2(k-1) }\Phi_0^k$
and $|\e|\be_0^{-2}\Phi_0 \le 1/2$. 
Finally, by Lemma \ref{lem:3.5}, one has $G^{1}_{\nu} 
=0$ if either $\nu=0$ or $|\nu|>\gotm N$.
\qed


\subsection{Perturbation series for the function $\boldsymbol{H}$}
\label{sec:3.5} 

In order to study \eqref{eq:2.11} and \eqref{eq:2.12}, it is convenient to introduce an auxiliary parameter $\mu$, by writing,
for $\nu\neq 0$,
\begin{equation} \nonumber
\begin{aligned}
& \left( \left( \om_{0}\cdot \nu \right)^2   - i \e^\gotm \ga \, \om_{0} \cdot \nu - \e 
F_{0}^{0}(1) \right) H_{\nu} \\
& \qquad \qquad
= \mu \Biggl( 
F_{\nu}^{1} 
+ G^{1}_{\nu} 
+ \e \sum_{\nu_0\in \ZZZ^2_*} F_{\nu_0}^{0}(1) \, H_{\nu-\nu_0} 
+ \sum_{q=2}^{\io} \e^{(\gotm +1)q-\gotm}
\!\!\!\!\!\!\!\! \sum_{\substack{ \nu_0 , \ldots, \nu_q \in \ZZZ^2 \\ \nu_0+\ldots+\nu_q=\nu}} \!\!\!\!\!\!\!\!
F_{\nu_0}^{0}(q) \, H_{\nu_1} \ldots H_{\nu_q} \Biggr) ,
\end{aligned}
\end{equation}
while requiring, for $\nu=0$, that
\begin{equation} \nonumber
\gamma \tilde \al = \mu \Biggl( 
F_{0}^{1} 
+  \e \sum_{\nu_0\in \ZZZ^2_*} F_{\nu_0}^{0}(1) \, H_{-\nu_0} + 
 \sum_{q=2}^{\io} \e^{(\gotm +1)q-\gotm} 
\!\!\!\!\! \sum_{\substack{ \nu_0,\nu_1,\ldots,\nu_q \in \ZZZ^2 \\ \nu_0+\ldots+\nu_q=0}} \!\!\!\!\! 
F_{\nu_0}^0(q) \, H_{\nu_1} \ldots H_{\nu_q} \Biggr)  .
\end{equation}
Note that, for $\mu=1$, the two equations reduce to \eqref{eq:2.11} and \eqref{eq:2.12}, respectively.

If one writes
\begin{equation} \label{eq:3.9}
H(\psi;\e) = \sum_{k=1}^{\io} \mu^{k} H^{(k)}(\psi) =
\sum_{k=1}^{\io} \sum_{\nu\in\ZZZ^2} \mu^{k} {\rm e}^{i\nu\cdot \psi} H^{(k)}_{\nu} , \qquad
\tilde \al = \sum_{k=1}^{\io} \mu^{k} \al^{(k)} ,
\end{equation}
then one obtains the recursive equations
\begin{subequations} \label{eq:3.10}
\begin{align}
D(\om_{0}\cdot\nu) \, H^{(1)}_{\nu} & = F_{\nu}^{1} 
+ G^{1}_{\nu} 
\label{eq:2.8a} \\
\gamma \al^{(1)} & = F_{0}^{1} 
,
\label{eq:2.8b}
\end{align}
\end{subequations}
for $k=1$,
\vspace{-.2cm}
\begin{subequations} \label{eq:3.11}
\begin{align}
D(\om_{0}\cdot\nu) \, H^{(2)}_{\nu} & = \e
\sum_{\nu_0\in \ZZZ^2_*} F_{\nu_0}^{0}(1) \, H^{(1)}_{\nu-\nu_0} , 
\label{eq:3.11a} \\
\gamma \al^{(2)} & = \e
\sum_{\nu_0\in \ZZZ^2_*} F_{\nu_0}^{0}(1) \, H^{(1)}_{-\nu_0} ,
\label{eq:3.119b}
\vspace{-.2cm}
\end{align}
\end{subequations}
for $k=2$, and, for $k \ge 3$.
\vspace{-.1cm} 
\begin{subequations} \label{eq:3.12}
\begin{align}
\hspace{-.3cm} 
D(\om_{0}\cdot\nu) \, H^{(k)}_{\nu} & \! = \! \e \!\!
\sum_{\nu_0\in \ZZZ^2_*} F_{\nu_0}^{0}(1) \, H^{(k-1)}_{\nu-\nu_0} + 
\sum_{q=2}^{\io} \e^{(\gotm+1)q-\gotm}
\!\!\!\!\!\!\!\! 
\sum_{\substack{ \nu_0 , \ldots, \nu_q \in \ZZZ^2 \\ \nu_0+\ldots+\nu_q=\nu}} \!\!\!\!\!\!\!\!
F_{\nu_0}^{0}(q)\!\!\!\!\!\!\!\!
\sum_{\substack{k_{1},\ldots,k_q \ge 1 \\ k_{1}+\ldots+k_q=k-1}} \!\!\!\!\!\!\!\!
H^{(k_{1})}_{\nu_1} \ldots H^{(k_q)}_{\nu_q} ,
\label{eq:3.12a} \\
\hspace{-.3cm} 
\gamma \al^{(k)} & \! = \! \e \!\!
\sum_{\nu_0\in \ZZZ^2_*} F_{\nu_0}^{0}(1) \, H^{(k-1)}_{-\nu_0} + 
\sum_{q=2}^{\io} \e^{(\gotm+1)q-\gotm} 
\!\!\!\!\!\!\!\! 
\sum_{\substack{ \nu_0,\nu_1,\ldots,\nu_q \in \ZZZ^2 \\ \nu_0+\ldots+\nu_q=0}} \!\!\!\!\!\!\!\! 
F^0_{\nu_0}(q) \!\!\!\!\!\!\!\!
\sum_{\substack{k_{1},\ldots,k_q \ge 1 \\ k_{1}+\ldots+k_q=k-1}} \!\!\!\!\!\!\!\!
H^{(k_{1})}_{\nu_1} \ldots H^{(k_q)}_{\nu_q}  .
\label{eq:3.12b}
\end{align}
\end{subequations}
%


The quantities $D(\om_{0}\cdot\nu)$ are bounded from below by Lemma \ref{lem:3.7}. 
By looking at the recursive equations \eqref{eq:3.10}-\eqref{eq:3.12}
and using that the coefficients $F_{\nu}^1(q)$ and $F_{\nu}^0(q)$ are well defined 
and decay exponentially in $\nu$
(by Lemmas \ref{lem:3.9} and \ref{lem:3.11}), one easily proves by induction that the coefficients 
$H^{(k)}_{\nu}$ are well defined for all $k\in\NNN$ and all $\nu\in\ZZZ^2_*$.
\qed

\subsection{Tree representation for the function $\boldsymbol{H}$}
\label{sec:3.6} 

Hereafter we consider trees which, in addition to the constraints considered in Section \ref{sec:1}, are such that
for any $v\in V(\theta)$, if $q_v=1$, then one has $\nu_{v} \neq 0$; this implies that
if there is only one line $\ell$ entering a node, then $\ell$ cannot have the same momentum as the line exiting that node.

We assign a further label $\lambda_v=F,G$ to each node $v\in E(\vartheta)$.
We call $\lambda_v$ the
\emph{end node label} and set $E_{\lambda}(\vartheta) :=\{v \in E(\vartheta) : \lambda_v = \lambda\}$,
so that $E(\vartheta)= E_{F}(\vartheta) \sqcup E_{G}(\vartheta)$.
Let $\gotT_{k,\nu}$ denote the set of all labelled trees with $|N(\vartheta)|=k$ such that the root line has momentum $\nu$.

We aim to show also that also the coefficients $H^{(k)}_{\nu}$ of the function $H$ admit a graphical representation in terms of trees,
with the further constraint mentioned above: of course, what changes, with respect to the coefficients of $h$,
are the node factors and the propagators to be associated with the nodes and the lines, respectively.

We associate with each node $v\in N(\vartheta)$ a \textit{node factor}
$\calI_v$, if $v\in I(\vartheta)$, and $\calE_v$, if $v\in E(\vartheta)$, and 
with each line $\ell\in L(\vartheta)$ \textit{propagator} $\calG_\ell$, with
\begin{equation} \nonumber 
\calI_v := 
\frac{1}{q_v!} F_{\nu_v}^0(q_v),  \qquad
\calE_v := 
\begin{cases}
F_{\nu_v}^{1} 
, & \lambda_v = F , \\ 
& \\
G^{1}_{\nu_v} 
, & \lambda_v = G , \end{cases}
%
%
\qquad
\calG_\ell:=\begin{cases}
\displaystyle{\frac{1}{D(\om_{0}\cdot\nu_\ell)}} , & \nu_{\ell} \neq 0 , \\ 
\vspace{-.4cm} \\
\displaystyle{ {\ga}^{-1}} & \nu_{\ell} = 0 .
\end{cases}
\end{equation}

Define the \emph{value} of the tree $\vartheta$ as
\begin{equation} \nonumber 
\Val(\vartheta) := 
\Bigl( \prod_{v \in I(\vartheta)} \e^{(\gotm+1)q_v - \gotm} 
\calI_v \Bigr) 
\Bigl( \prod_{v \in E(\vartheta)} 
\calE_v \Bigr)
\Bigl( \prod_{\ell \in L(\vartheta)} \calG_{\ell} \Bigr) .
\end{equation}
%
%
By construction only the root line may have zero momentum. 
Thus, an end node $v\in E(\vartheta)$ may have $\nu_v =0$ if and only if $\vartheta\in\gotT_{1,0}$,
and in such a case one has $\Val(\vartheta)=\ga^{-1}F^1_0 
$.

\begin{lemma} \label{lem:3.12}
Assume Hypothesis \ref{hyp:1},
and assume $\om_0$ to be non-resonant up to order $\gotm N \!$.
For $|\e| \le \e_2$, with $\e_2$ defined in Lemma \ref{lem:3.7},
and or any $k\in\NNN$ one has
\begin{equation} \nonumber 
H^{(k)}_{\nu} := \sum_{\vartheta \in \gotT_{k,\nu}}
 \Val(\vartheta), \qquad \nn \in \ZZZ^2_* , \qquad
\al^{(k)} := \sum_{\vartheta \in \gotT_{k,0}} \Val(\vartheta) .
\end{equation}
\end{lemma}

\noindent\emph{Proof}. 
By induction, as in the proof of Lemma \ref{lem:3.4}. \qed

\zerarcounters 
\section{Convergence of the perturbation series}
\label{sec:4} 

%

Define, for any tree $\vartheta\in\gotT_{k,\nu}$, with $(k,\nu) \in\NNN \times \ZZZ^2 \setminus \{(1,0)\}$,
\begin{equation} \nonumber 
\calW(\vartheta) := 
\Bigl( \prod_{v \in I(\vartheta)} |\gotD_0\e|^{(\gotm +1)q_v - \gotm + \frac{3}{4} p(\nu_v)} \Bigr) 
\Bigl( \prod_{v \in E(\vartheta)} |\gotD_0\e|^{\frac{3}{4} p_{\gotm}(\nu_v)} \Bigr)
\Bigl( \prod_{\ell \in L(\vartheta)} |\calG_{\ell}| \Bigr) ,
\end{equation}
so that, by taking into account Lemma \ref{lem:3.1}, one can bound, for any tree $\gotT_{k,\nu}$ with $(k,\nu)\neq (1,0)$,
\begin{equation} \label{eq:4.0}
\Val(\vartheta) \le \calW(\vartheta) \, \gotD_0^{\gotm-(k-1)} N^{k-1}
\Bigl( \!\!\! \prod_{v \in I(\vartheta)}  (\gotD_0 |\e|)^{\frac{1}{4} p(\nu_v)} 3\Phi  \Bigr) 
\Bigl( \!\!\! \prod_{v \in E_{F}(\vartheta)}
 (\gotD_0 |\e|)^{\frac{1}{4} p_{\gotm}(\nu_v)} 6\Phi \Bigr)
\Bigl( \!\!\! \prod_{v \in E_{G}(\vartheta)} \Gamma \Bigr) .
\end{equation}
%

\begin{rmk} \label{rmk:4.1}
\emph{
If $|\e|\le \e_2$, with $\e_2$ as in Lemma \ref{lem:3.7}, then one has $|\gotD_0\e|^{\frac18} \le 1/2$.
Indeed, by definition of $\gotD_0$ and $\e_1$, one has $\gotD_0\e_1=1/2$, so that,
for $|\e| \le \e_2$, by using Remark \ref{rmk:3.8},
one obtains $|\gotD_0\e| \le 2 \,16^{-2} |\gotD_0\e_1| \le 16^{-2}$.
}
\end{rmk}

\begin{lemma} \label{lem:4.2}
Assume Hypothesis \ref{hyp:1},
and assume $\om_0$ to be non-resonant up to order $\gotm N \!$.
Let $\gotD_0$ be defined as in Lemma \ref{lem:3.9}.
There exists a positive constant $\gotB$
such that, for all $k\in\NNN$, all $\nu\in\ZZZ^2_*$ and all $|\e| \le \e_2$, with $\e_2$ as in Lemma \ref{lem:3.7}, 
one has
\begin{equation} \nonumber
\calW(\vartheta) \le 
\gotB^{k} |\gotD_0\e|^{\frac{k-1}{4}}
\end{equation}
for any tree $\vartheta\in\gotT_{k,\nu}$.
\end{lemma}

\noindent\emph{Proof}. 
The proof is by induction on $k$.
If $\vartheta\in\gotT_{1,\nu}$,
then $I(\vartheta)$ is empty, $E(\vartheta)$ contains only one node $v_0$ and
$L(\vartheta)$ contains only the root line $\ell_0$. One has 
\begin{equation} \nonumber
\calW(\vartheta) = |\gotD_0\e|^{\frac{3}{4} p_{\gotm}(\nu)} \, | \calG_{\ell_0} | ,
\end{equation}
where $\nu=\nu_{v_0}$. If $p(\nu) < 4\gotm$, one has $|\calG_{\ell_0}| \le 16\be_0^{-2}$, while if 
$p(\nu) \ge 4\gotm$ (and hence
$p_\gotm(\nu) \ge 3\gotm$), 
\begin{equation} \nonumber
| \gotD_0\e|^{\frac{3}{4} p_\gotm(\nu)} |\calG_{\ell_0}| \le
|\gotD_0\e|^{{\frac{9}{4}\gotm}} C_0^{-1} \gotD_0^2 |\gotD_0 \e|^{-2} < C_0^{-1} \gotD_0^{2} ,
\end{equation}
so that the bound holds if one takes
%
$\gotB \ge \max\{ 16\be_0^{-2}, C_0^{-1} \gotD_0^{2}  \} = C_0^{-1} \gotD_0^{2}$.

If $\vartheta\in\gotT_{k,\nu}$, with $k>1$, let $v_0$ be the node which the root line exits and let
$\vartheta_1,\ldots,\vartheta_q$, with $q\ge 1$, denote the subtrees which have $v_0$ as root. By construction one has
\begin{equation} \nonumber
\calW(\vartheta) = |\gotD_0\e|^{(\gotm+1)q -\gotm + \frac{3}{4} p(\nu_{v_0})} |\calG_{\ell_0}| \prod_{i=1}^{q} \calW(\vartheta_i) .
\end{equation}
\begin{itemize}[leftmargin=.5cm]
\itemsep0em
\item[1.] If $q\ge 2$, one uses the inductive hypothesis to bound
\begin{equation} \nonumber
\prod_{i=1}^{q} \calW(\vartheta_i) \le 
\prod_{i=1}^{q} \gotB^{k_i} |\gotD_0\e|^{\frac{k_i-1}{4}} = \gotB^{k-1} |\gotD_0\e|^{\frac{k-1-q}{4}} ,
\end{equation}
since one has $k_{1}+\ldots+k_q=k-1$. Thus, for $q\ge 2$, one obtains
\begin{equation} \nonumber
\calW(\vartheta) \le  |\gotD_0\e|^{(\gotm+1)q -\gotm } \max\{ 16\be_0^{-2} , C_0^{-1} \gotD_0^2 |\gotD_0\e|^{-2} \}
\, \gotB^{k-1} |\gotD_0\e|^{\frac{k-1-q}{4}} .
\end{equation}
By using that
\begin{equation} \nonumber
(\gotm+1)q - \gotm - 2 + \frac{k-1-q}{4} = \Bigl( \gotm+\frac{3}{4} \Bigr) q - (\gotm+2) + \frac{k-1}{4} > \frac{k-1}{4} ,
\end{equation}
the bound follows by requiring
%
$\gotB \ge C_0^{-1} \gotD_0^2$.
%
\item[2.]
If $q=1$ one has
\begin{equation} \nonumber
\calW(\vartheta) =|\gotD_0\e|^{1 + \frac{3}{4} p(\nu_{v_0})} |\calG_{\ell_0}| \calW(\vartheta_1) ,
\end{equation}
so that, if $|\om_{0}\cdot\nu| \ge \be_0/4$, the bound follows once more by the inductive hypothesis, since 
\begin{equation} \nonumber
\calW(\vartheta) \le  |\gotD_0\e| 16\be_0^{-2} \gotB^{k-1} |\gotD_0\e|^{\frac{k-2}{4}} 
< 16\be_0^{-2} \gotB^{k-1} |\gotD_0\e|^{\frac{k-1}{4}} ,
\end{equation}
which implies the bound provided one has
%
$\gotB \ge 16\be_0^{-2}$.
%
If $|\om_{0}\cdot\nu|<\be_0/4$ (but still $|\calG_{\ell_0}| \le C_0^{-1}|\e|^{-2}$ by Lemma \ref{lem:3.7}),
then one has $p(\nu)\ge4\gotm$ and hence $p_{\gotm}(\nu) \ge 3\gotm$.
Let $\ell_1$ denote the root line of $\vartheta_1$, $v_1$ the node which $\ell_1$ exits, and
$\vartheta_1',\ldots,\vartheta_{q'}'$ the subtrees which have $v_1$ as root.
\begin{itemize}[leftmargin=.31cm]
\itemsep0em
\item[2.1.] If $q'=0$ (i.e.~$k=2$) we distinguish between the following cases.
\begin{itemize}[leftmargin=.3cm]
\itemsep0em
\item[2.1.1.] If $|\om_0\cdot \nu_{\ell_1}| \ge \be_0/4$, then, by Lemma \ref{lem:3.3}, one bounds
\begin{equation} \nonumber
\calW(\vartheta) \le|\gotD_0\e|^{-1+\frac{3}{4}(p(\nu_{v_0})+p_{\gotm}(\nu_{\ell_1}))} C_0^{-1}\gotD_0^{2}16\be_0^{-2}
\le \gotB^2 |\gotD_0\e|^{\frac{9\gotm - 7}{4}} < \gotB^2 |\gotD_0\e|^{\frac{1}{4}} ,
\end{equation}
if $\gotB \ge C_0^{-1}\gotD_0^2$, since $p(\nu_{v_0}) + p_{\gotm}(\nu_{\ell_1}) \ge
p(\nu_{v_0} + \nu_{\ell_1}) - 1 \ge 3\gotm - 1$.
\item[2.1.2.] If $|\om_0\cdot \nu_{\ell_1}|<\be_0/4$, then $p(\nu_{\ell_1}) \ge 4\gotm$ and hence
$p_{\gotm}(\nu_{\ell_1}) \ge 3\gotm$. Moreover
\begin{equation} \nonumber
|\om_0\cdot\nu_{v_0}|=|\om_0\cdot(\nu-\nu_{\ell_1})|\le|\om\cdot\nu|+|\om\cdot\nu_{\ell_1}|<\frac{\be_0}{2}
\end{equation}
implies $p(\nu_{v_0})\ge 4\gotm$, so that one bounds
\begin{equation} \nonumber
\calW(\vartheta) \le|\gotD_0\e|^{1+\frac{3}{4}(p(\nu_{v_0})+p_{\gotm}(\nu_{\ell_1}))} C_0^{-2} |\e|^{-4}
\le |\gotD_0\e|^{\frac{21\gotm -12}{4}} C_0^{-2} \gotD_0^4 < \gotB^2 |\gotD_0\e|^{\frac{1}{4}} ,
\end{equation}
if $\gotB \ge C_0^{-1}\gotD_0^2$.
\end{itemize}
\item[2.2.] If $q'\ge1$ we distinguish among the following three cases.
\begin{itemize}[leftmargin=.3cm]
\itemsep0em
\item[2.2.1.] If $|\om_{0}\cdot\nu_{\ell_1}| < \be_0/4$, one bounds $|\calG_{\ell_1}| \le C_0^{-1} |\e|^{-2}$ and
one has $p(\nu_{v_0}) \ge 4\gotm$, otherwise one would find
\begin{equation} \nonumber
\be_0 \le |\om_{0}\cdot\nu_{v_0}| \le |\om_{0}\cdot\nu| + |\om_{0}\cdot\nu_{\ell_1}| < \frac{\be_0}{2} ,
\end{equation}
which would lead to a contradiction. Thus one finds 
\begin{eqnarray} \nonumber
\calW(\vartheta)
& \!\!\! \le \!\!\! &  
C_0^{-1} |\e|^{-2}|\gotD_0\e|^{1+\frac34 p(\nu_{v_0})} \, C_0^{-1}|\e|^{-2} \, |\gotD_0\e|^{(\gotm+1)q'-\gotm +\frac34 p(\nu_{v_1})}
\prod_{i=1}^{q'} \calW(\vartheta_i') \nonumber \\
& \!\!\! \le \!\!\! &
C_0^{-2} \gotD_0^4  \gotB^{k-2} |\gotD_0\e|^{(\gotm +1)q'+2\gotm-3+\frac{k-2-q'}{4}} , \nonumber
\end{eqnarray}
where one has $\gotm,q'\ge 1$ and hence
\begin{equation} \nonumber
(\gotm +1)q'+2\gotm-3+\frac{k-2-q'}{4}  \ge \frac{k-1}{4} ,
\end{equation}
which yields the bound if one requires
%
$\gotB \ge C_0^{-2} \gotD_0^4 \gotB^{-1}$.
%
\item[2.2.2.] If $|\om_{0}\cdot\nu_{\ell_1}| \ge \be_0/4$ and $q'\ge 2$, one has
\begin{eqnarray} \nonumber
\calW(\vartheta) 
& \!\!\! \le \!\!\! &
C_0^{-1} \gotD_0^2 |\gotD_0\e|^{-2} |\gotD_0\e|^{1+\frac34(p(\nu_{v_0})+p(\nu_{v_1})) + (\gotm+1)q'-\gotm}
 \, 16\be_0^{-2} \prod_{i=1}^{q'} \left| \calW(\vartheta_i') \right| \nonumber \\
& \!\!\! \le \!\!\! &
C_0^{-1} \gotD_0^2 16\be_0^{-2} \gotB^{k-2} |\gotD_0\e|^{(\gotm+1)q'-\gotm-1+\frac{k-2-q'}{4}} , \nonumber
\end{eqnarray}
where one has
\begin{equation} \nonumber
(\gotm+1)q'-(\gotm+1)+\frac{k-2-q'}{4} = \Bigl( \gotm +\frac34 \Bigr) q' - (\gotm+1) + \frac{k-2}{4} \ge \frac{k-1}{4} ,
\end{equation}
so that the bound follows provided one has
%
$\gotB \ge C_0^{-1} \gotD_0^2 16 \be_0^{-2} \gotB^{-1}$.
%
\item[2.2.3.] If $|\om_{0}\cdot\nu_{\ell_1}| \ge \be_0/4$ and $q'=1$, let 
$\ell_1'$ denote the root line of $\vartheta_1'$, $v_1'$ the node which $\ell_1'$ exits, and
$\vartheta_1'',\ldots,\vartheta_{q''}''$
the subtrees which have $v_1'$ as root. 
\begin{itemize}[leftmargin=.6cm]
\itemsep0em
\item[2.2.3.1.\hspace{.3cm}] \hspace{-.4cm}
If $q''=0$ (i.e. $k=3$), we distinguish between two cases.
%
\item[2.2.3.1.1.] If $|\om_0\cdot\nu_{\ell_1'}|\ge \be_0/4$ then, since $p(\nu) \ge 4\gotm$ and hence $p_{\gotm}(\nu) \ge 3\gotm$,
one obtains, if $\gotB \ge C_0^{-1} \gotD_0^2$,
\vspace{-.2cm}
\begin{equation} \nonumber
\calW(\vartheta) \le C_0^{-1} \gotD_0^2 \, (16 \be_0^{-2})^{2} |\gotD_0\e|^{
\frac{3}{4}(p(\nu_{v_0}) + p(\nu_{v_1}) + p_{\gotm}(\nu_{v'_1}))}\le \gotB^{3}  |\gotD_0\e|^{\frac{3}{4}} <
\gotB^{3}  |\gotD_0\e|^{\frac{1}{2}} ,
\vspace{-.1cm}
\end{equation}
because, by Lemma \ref{lem:3.3}, 
one has $p(\nu_{v_0}) + p(\nu_{v_1}) + p_{\gotm}(\nu_{v'_1})) \ge p_{\gotm}(\nu) - 2 \ge 3\gotm-2$.
\item[2.2.3.1.2.] If $|\om_0\cdot\nu_{\ell_1'}|\le \be_0/4$, then one has $p_{\gotm}(\nu_{v_1'})\ge3\gotm$. Moreover one has
\vspace{-.2cm}
\begin{equation} \nonumber
|\om_0\cdot (\nu_{v_0}+\nu_{v_1})| = |\om_{0}\cdot(\nu - \nu_{\ell_1'})| 
 \le |\om_0\cdot \nu| +   |\om_0\cdot \nu_{\ell_1'}| < \frac{\be_0}{2} ,
\vspace{-.1cm}
\end{equation}
so that, since, by Lemma \ref{lem:3.3}, one has
$p(\nu_{v_0}) + p(\nu_{v_1}) \ge p(\nu_{v_0}+\nu_{v_1}) - 1 \ge 4\gotm - 1$, one obtains
\begin{equation} \nonumber
\calW(\vartheta) \!\le \! C_0^{-2} \gotD_0^4 16 \be_0^{-1} |\gotD_0\e|^{-2+\frac{3}{4}(p(\nu_{v_0}) + 
p(\nu_{v_1})+p_{\gotm}(\nu_{v'_1}))} \!\le \! \gotB^3 |\gotD_0\e|^{\frac{21\gotm-11}{4}}
\!\!\! < \! \gotB^3 |\gotD_0\e|^{\frac{1}{2}} ,
\end{equation}
provided that $\gotB \ge C_0^{-1} \gotD_0^2$.
\item[2.2.3.2.\hspace{.31cm}]
If one has $q''\ge1$, we still have two cases to discuss.
\item[2.2.3.2.1.] If one has $|\om_{0}\cdot\nu_{\ell_1'}| \ge \be_0/4$, one obtains
\vspace{-.2cm}
\begin{equation} \nonumber
\calW(\vartheta) \le C_0^{-1} \gotD_0^2 |\gotD_0\e|^{-2} |\gotD_0\e|^{2} \, (16\be_0^{-2})^2 
\gotB^{k-3} |\gotD_0\e|^{(\gotm+1)q''-\gotm+\frac{k-3-q''}{4}} ,
\vspace{-.1cm}
\end{equation}
where one has
\vspace{-.2cm}
\begin{equation} \nonumber
- 2 + 2 + (\gotm+1)q'' - \gotm + \frac{k-3-q''}{4} = \Bigl(\gotm+ \frac{3}{4} \Bigr) q''  -\gotm + \frac{k-3}{4} > \frac{k-1}{4} ,
\vspace{-.1cm}
\end{equation}
so that the bound holds as soon as one requires
%
$\gotB \ge C_0^{-1} \gotD_0^2 ( 16 \be_0^{-2})^2 \gotB^{-2}$.
%
\item[2.2.3.2.2.]
If, instead, one has $|\om_{0}\cdot\nu_{\ell_1'}| < \be_0/4$, one obtains
\vspace{-.2cm}
\begin{equation} \nonumber
\calW(\vartheta) \! \le \! ( C_0 |\e|^{2})^{-2} 16\be_0^{-2}
|\gotD_0 \e|^{2+\frac{3}{4} (p(\nu_{v_0}) + p(\nu_{v_1}))} 
\gotB^{k-3} 
|\gotD_0\e|^{(\gotm+1)q'' \! -\gotm+\frac{k-3-q''}{4}} , 
\end{equation}
and, since $\nu=\nu_{v_0}+\nu_{v_1}+\nu_{\ell_1'}$, one has
\vspace{-.2cm}
\begin{equation} \nonumber
|\om_{0}\cdot(\nu_{v_0}+\nu_{v_1})| = |\om_{0}\cdot(\nu - \nu_{\ell_1'})| 
\le |\om_{0}\cdot\nu| + |\om_{0}\cdot\nu_{\ell_1'}| < \frac{\be_0}{2} ,
\vspace{-.1cm}
\end{equation}
which in turns implies $p(\nu_{v_0}+\nu_{v_1})\ge 4\gotm$ and hence $p(\nu_{v_0})+p(\nu_{v_1}) \ge 4\gotm-1$ by 
Lemma \ref{lem:3.3}. Eventually one finds
\vspace{-.2cm}
\begin{equation} \nonumber
\calW(\vartheta) \le C_0^{-2} \gotD_0^4 16\be_0^{-2} \gotB^{k-3} 
|\gotD_0\e|^{\frac34(4\gotm-1)-2+(\gotm+1)q''-\gotm +\frac{k-3-q''}{4}} ,
\vspace{-.1cm}
\end{equation}
where, using that $\gotm,q''\ge 1$, one has
\vspace{-.2cm}
\begin{equation} \nonumber
\frac34(4\gotm-1)-2+(\gotm+1)q''-\gotm +\frac{k-3-q''}{4} > \frac{k-1}{4} ,
\vspace{-.1cm}
\end{equation}
so that the bound follows once more provided that one requires
%
$\gotB \ge C_0^{-2} \gotD_0^4 16 \be_0^{-2} \gotB^{-2}$.
%
%
\end{itemize}
\end{itemize}
\end{itemize}
\end{itemize}
By collecting together the conditions required on the constant $\gotB$ and taking into
account Remark \ref{rmk:4.1}, one may fix $\gotB$ as
\begin{equation} \label{eq:4.1}
\gotB = C_0^{-1} \gotD_0^{2}.
\end{equation}
%
%
This concludes the proof. \qed

\begin{lemma} \label{lem:4.3}
Assume Hypothesis \ref{hyp:1},
and assume $\om_0$ to be non-resonant up to order $\gotm N \!$.
Let $\gotD_0$ be defined as in Lemma \ref{lem:3.9}.
For all $k\in\NNN$, with $k\ge 2$, and all $|\e| \le \e_2$, with $\e_2$ as in Lemma \ref{lem:3.7}, one has
\begin{equation} \nonumber
\calW(\vartheta) \le 
\ga^{-1} \gotB^{k-1} |\gotD_0\e|^{\frac{k+2}{4}} 
\end{equation}
for any tree $\vartheta\in\gotT_{k,0}$, with $\gotB$ the same constant as in Lemma \ref{lem:4.2}.
\end{lemma}

\noindent\emph{Proof}. 
Given a tree $\vartheta\in\gotT_{k,0}$, with $k\ge 2$, let $v_0$ be the node which the root line of $\vartheta$ exits
and let $\vartheta_1,\ldots,\vartheta_q$, with $q\ge 1$, be the subtrees which have $v_0$ as root. One has
%
\begin{equation} \nonumber
\calW(\vartheta) =  \gamma^{-1} |\gotD_0\e|^{(\gotm+1)q-\gotm + \frac{3}{4} p(\nu_{v_0})} \prod_{i=1}^{q} \calW(\vartheta_i) ,
\end{equation}
where each $\calW(\vartheta_i)$ can be bounded by using Lemma \ref{lem:4.2}, since the momentum of the
corresponding root line is different from zero. Thus, one finds
%
\begin{equation} \nonumber
\calW(\vartheta)\le \gamma^{-1} \gotB^{k-1} 
|\gotD_0 \e|^{(\gotm+1)q-\gotm + \frac{k-1-q}{4}} \le
\ga^{-1} \gotB^{k-1} |\gotD_0 \e|^{\frac{k+2}{4}} , 
\end{equation}
where we have 
%
%
used that 
%
\begin{equation} \nonumber
(\gotm+1)q - \gotm + \frac{k-1-q}{4}  = \frac{k -1+ 3q }{4}+\gotm(q-1) \ge \frac{k+2}{4} ,
\end{equation}
which implies the bound. 
\qed

\begin{rmk} \label{rmk:4.4}
\emph{
The bound in Lemma \ref{lem:4.3} may be improved for $k\ge 3$. First of all, note that, for $q\ge 2$,
the last inequality in the proof implies that the exponent of  $|\gotD_0\e|$ is bounded by $(k+9)/4$ for any $k$.
If $q=1$ and $k\ge 3$, let $\ell_1$ be the line entering $v_0$, $v_1$ be the node which $\ell_1$ exits,
and $\theta_1',\ldots,\theta_{q'}'$, with $q'\ge 1$, the subtrees which have $v_1$ as root.
We distinguish between two cases.
\vspace{-.2cm}
\begin{itemize}[leftmargin=.5cm]
\itemsep0em
%
%
%
%
\item[1.] If $|\om_0\cdot\nu_{\ell_1}| < \be_0/4$, one has
%
\begin{equation} \nonumber
\begin{aligned}
\calW(\vartheta) & \le \ga^{-1} |\gotD_0\e|^{1 + \frac34 p(\nu_{v_0}) - 2 + (\gotm+1)q'-\gotm} C_0^{-1} \gotD_0^2
\gotB^{k-2} |\gotD_0\e|^{\frac{k-2-q'}{4}} \\
& \le \ga^{-1} |\gotD_0\e|^{3\gotm} C_0^{-1} \gotD_0^2 \gotB^{k-2} |\gotD_0\e|^{\frac{k-3}{4}} \le
\ga^{-1} \gotB^{k-1} |\gotD_0\e|^{\frac{12\gotm + k-3}{4}} \le \ga^{-1} \gotB^{k-1} |\gotD_0\e|^{\frac{k+9}{4}} ,
\end{aligned}
\end{equation}
because one has $p(\nu_{\ell_1}) \ge 4\gotm$ and $\nu_{v_0}+\nu_{\ell_1}=0$.
\item[2.] If $|\om_0\cdot\nu_{\ell_1}| \ge \be_0/4$, one has
%
\begin{equation} \nonumber
\calW(\vartheta) \le \ga^{-1} |\gotD_0\e|^{1 + (\gotm+1)q'-\gotm} 16 \be_0^{-2} \gotB^{k-2} |\gotD_0\e|^{\frac{k-2-q'}{4}} 
\le \ga^{-1} \gotB^{k-1} |\gotD_0\e|^{\frac{k+5}{4}} .
\end{equation}
\end{itemize}
Therefore, for $k\ge 3$ and for any tree $\vartheta\in\TT_{k,0}$, one can bound
$\calW(\vartheta) \le \ga^{-1} \gotB^{k-1} |\gotD_0\e|^{\frac{k+5}{4}}$.
}
\end{rmk}

\begin{lemma} \label{lem:4.5}
There exists a positive constant $\gotM_0$
such that, for all $k\in\NNN$, 
one has
\begin{equation} \nonumber
\Bigl| \Bigl( \prod_{v \in I(\vartheta)} 3\Phi \Bigr) 
\Bigl( \prod_{v \in E_{F}(\vartheta)} 6\Phi \Bigr)
\Bigl( \prod_{v \in E_{G}(\vartheta)} \Gamma  \Bigr) \Bigr| 
\le \gotM_0^k .
%
\end{equation}
%
%
%
\end{lemma}

\noindent\emph{Proof}. 
%
%
Define
\begin{equation} \label{eq:4.2}
\gotM_0 := \max\{ 3\Phi, 6\Phi ,\Gamma \} = 
\max\{ 6\Phi , 
2 \ga \be_0^{-1} \Phi_0 \gotD_0^{-\gotm} \} ,
\end{equation}
and use that
%
$|I(\vartheta)|+|E(\vartheta)|=|N(\vartheta)|=k$.
%
Then the assertion follows with $\gotM_0$ given by \eqref{eq:4.2}. \qed

\begin{lemma} \label{lem:4.6}
Assume Hypothesis \ref{hyp:1},
and assume $\om_0$ to be non-resonant up to order $\gotm N \!$.
Let $p(\nu)$ and $p_{\gotm}(\nu)$ be defined as in \eqref{eq:3.2}.
There exists a positive constant 
$\gotR_0$ 
such that for all $|\e| \le \e_2$, with $\e_2$ as in Lemma \ref{lem:3.7}, one has
\begin{equation} \nonumber
\sum_{\nu\in\ZZZ^2_*} | \gotD_0 \e|^{\frac18 p(\nu)} \le 
\sum_{\nu\in\ZZZ^2_*} | \gotD_0 \e|^{\frac18 p_{\gotm}(\nu)} \le \gotR_{0} .
\end{equation}
\end{lemma}

\noindent\emph{Proof}. 
For $|\e|\le \e_2$ one has $|\gotD_0\e| \le 1/2^8$ (cf.~Remark \ref{rmk:4.1}), so that the first sum is less than the second one.
For any $x\in [0,1)$ one has
%
%
%
%
%
\begin{eqnarray} \nonumber
\sum_{\nu\in\ZZZ^2_*} x^{p_{\gotm}(\nu)}
& \!\!\! \le \!\!\! &
 \sum_{\substack{\nu\in\ZZZ^2_* \\ p(\nu)\le \gotm}}1
+ \sum_{p=\gotm+1}^{\io} \sum_{\substack{\nu\in\ZZZ^2_* \\ p(\nu)=p}} x^{p-\gotm} 
\le 4 (\gotm+1)^2 N^2 + \sum_{p=0}^{\io} x^{p+1} 4 (p + 1 + (\gotm+1)) N^2
\nonumber \\
& \!\!\! \le \!\!\! &
4 N^2 \left( (\gotm+1)^2 + \frac{x}{(1-x)^2} + \frac{(\gotm + 1)x}{(1-x)}\right) ,
\nonumber
\end{eqnarray}
so that, if one sets
\begin{equation} \label{eq:4.4}
\gotR_{0} := 4 N^2 
\left( (\gotm+1)^2 + \gotm + 3 \right) ,
\end{equation}
the second sum is bounded by $\gotR_0$. 
%
%
Note that $4N^2(\gotm+1)^2 \le \gotR_0 \le 8N^2 (\gotm+1)^2$.
%
\qed

\begin{lemma} \label{lem:4.7}
For all $k\in\NNN$, all $\nu\in\ZZZ^2_* $ and all $|\e| \le \e_2$, with $\e_2$ as in Lemma \ref{lem:3.7}, one has
\begin{equation} \nonumber
\sum_{v \in V(\theta)} p(\nu_v) + \sum_{v \in E(\theta)} p_{\gotm}(\nu_v) \ge \max \left\{ p(\nu)-\gotm k - (k-1) , 0 \right\} 
\end{equation}
for all $\vartheta\in \gotT_{k,\nu}$.

\end{lemma}

\noindent\emph{Proof}. 
Since $p(\nu) \ge0$ and $p_{\gotm}(\nu)\ge 0$ for all $\nu\in\ZZZ_*^2$, the bound is proved if we show that one has
\begin{equation} \nonumber
\sum_{v \in I(\vartheta)} p(\nu_v) + \sum_{v \in E(\vartheta)} p_{\gotm}(\nu_v) \ge p(\nu) - \gotm k - (k-1)
\end{equation}
for all $|\nu| \ge (\gotm k + (k - 1))N$. By iteratively applying Lemma \ref{lem:3.3}, one finds
\begin{equation} \nonumber
\sum_{v \in I(\vartheta)} p(\nu_v) + \sum_{v \in E(\vartheta)} p(\nu_v) \ge p(\nu)- (k-1) ,
\end{equation}
since $|I(\vartheta)|+|E(\vartheta)|=k$, so that 
\begin{equation} \nonumber
\sum_{v \in I(\vartheta)} \!\! p(\nu_v) + \sum_{v \in E(\vartheta)} \!\! p_{\gotm}(\nu_v) \ge 
\sum_{v \in I(\vartheta)} \!\! p(\nu_v) + \sum_{v \in E(\vartheta)} \!\! p(\nu_v) - \gotm |E(\vartheta)| \ge
p(\nu)- \left( k - 1 + \gotm |E(\vartheta)| \right) ,
\end{equation}
with $|E(\vartheta)| \le k$. (In fact, one has $E(\vartheta) \le k-1$ if $k\ge 2$). Hence the assertion follows.
\qed

\begin{lemma} \label{lem:4.8}
Assume Hypothesis \ref{hyp:1},
and assume $\om_0$ to be non-resonant up to order $\gotm N \!$.
Set 
$\gotC_0 \!:=\! 2 \gotD_0^{-1} N \gotM_0 \gotR_0 \gotB$, with
the constants $\gotD_0$, $\gotB$, $\gotM_0$ and $\gotR_0$  
as in Lemmas \ref{lem:3.9}, \ref{lem:4.2}, \ref{lem:4.5} and \ref{lem:4.6}, respectively.
For all $k\in\NNN$ and all $|\e| \le \e_2$, with $\e_2$ as in Lemma \ref{lem:3.7}, one has
\begin{equation} \nonumber
| H^{(k)}_{\nu} | \le 
\gotD_0^{\gotm+1} N^{-1} 
\gotC_0^{k} 
|\gotD_0\e|^{\frac{k-1}{4} +
\max\{ \frac{1}{8}(p(\nu)-\gotm k -(k-1)) , 0 \} } 
\end{equation}
for all $\nu\in\ZZZ^2_*$ and
\vspace{-.2cm}
\begin{equation} \nonumber
|\al^{(1)}| 
\le \ga^{-1} 6 \Phi \gotD_0^{\gotm} |\gotD_0\e| , \qquad \quad 
|\al^{(k)}| 
\le \ga^{-1} \gotD_0^{\gotm+1} N^{-1} \gotB^{-1}
\gotC_0^{k} 
|\gotD_0\e|^{\frac{k+2}{4}} , \qquad k \ge 2 .
\end{equation}
\end{lemma}

\noindent\emph{Proof}. 
According to Lemmas \ref{lem:3.12} and \ref{lem:4.5}, we find
\begin{equation} \nonumber
H^{(k)}_{\nu} \! \le \!\!\! \sum_{\vartheta \in \gotT_{k,\nu}} \!\! |\VV(\vartheta)| , \quad
|\VV(\vartheta)| \! \le \! \calW(\vartheta)
 \, \gotD_0^{\gotm-(k-1)} N^{k-1} \gotM_0^k \Bigl( \prod_{v \in I(\vartheta)} \!\!  |\gotD_0\e|^{\frac{1}{4} p(\nu_v)} \Bigr) 
\Bigl( \!\! \prod_{v \in E_{F}(\vartheta)} |\gotD_0\e|^{\frac{1}{4} p_{\gotm}(\nu_v)} \Bigr),
\end{equation}
where we use have \eqref{eq:4.0}, and bound $\calW(\vartheta)$ as in Lemma \ref{lem:4.2}.
The sum over the end node labels produces at most a factor $2^{k-1}$ and, writing
%
\begin{equation} \nonumber
\frac{1}{4} p(\nu_v) = \frac{1}{8} p(\nu_v) + \frac{1}{8} p(\nu_v) , 
\qquad \frac{1}{4} p_{\gotm}(\nu_v) = \frac{1}{8} p_{\gotm}(\nu_v) + \frac{1}{8} p_{\gotm}(\nu_v),
\end{equation}
the sum over the mode labels of the product
%
\begin{equation} \nonumber
\Bigl( \prod_{v \in I(\vartheta)} |\gotD_0\e|^{\frac{1}{8} p(\nu_v)} \Bigr) 
\Bigl( \prod_{v \in E_{F}(\vartheta)} |\gotD_0\e|^{\frac{1}{8} p_{\gotm}(\nu_v)} \Bigr) ,
\end{equation}
neglecting the constraint that the sum of the modes equals $\nu$,
is dealt with by using Lemma \ref{lem:4.6}.
The remaining product is bounded as
\begin{equation} \nonumber
\Bigl( \prod_{v \in I(\vartheta)} |\gotD_0\e|^{\frac{1}{8} p(\nu_v)} \Bigr) 
\Bigl( \prod_{v \in E_{F}(\vartheta)} |\gotD_0\e|^{\frac{1}{8} p_{\gotm} (\nu_v)} \Bigr) \le |\gotD_0\e|^{\frac{1}{8}(p(\nu)-\gotm k - (k-1))} ,
\end{equation}
by Lemma \ref{lem:4.7}, if $p(\nu) > \gotm k - (k-1)$, and as 1 otherwise. The bound for $H^{(k)}_{\nu}$ then follows.

The bounds on the coefficients $\al^{(k)}$, for $k \ge 2$, are discussed in a similar way, relying on Lemma \ref{lem:4.3},
while the bound on $\al^{(1)}$ follows immediately from Lemma \ref{lem:3.11},
recalling that $\Val(\vartheta)=\ga^{-1}F^1_0 
$ if $\vartheta\in\gotT_{1,0}$.
\qed

\begin{lemma} \label{lem:4.9}
Assume Hypothesis \ref{hyp:1},
and assume $\om_0$ to be non-resonant up to order $\gotm N \!$.
Let $\gotD_0$ and $\gotC_0$ be defined as in Lemma \ref{lem:3.9} and \ref{lem:4.8}, respectively.
For all $k\in\NNN$ and all $|\e| \le \e_2$, with $\e_2$ as in Lemma \ref{lem:3.7}, 
the functions $H^{(k)}(\psi)$ in \eqref{eq:3.9},
%
%
for $\psi\in\TTT^2$, are bounded as
\begin{equation} \nonumber
| H^{(k)}(\psi) | \le  k^2 \gotC \, \gotD_0^{\gotm+1} N^{-1} \gotC_0^k 
|\gotD_0\e|^{\frac{k-1}{4}} ,
\end{equation}
for a suitable constant $\gotC$.
\end{lemma}

\noindent\emph{Proof}. 
Write $H^{(k)}(\psi)$ as
\begin{equation} \nonumber
H^{(k)}(\psi) = 
\!\!\!\!\!\!\!\!\!\!\! \sum_{\substack{ \nu\in\ZZZ^2_* \\ p(\nu) < \gotm k + (k-1)}} \!\!\!\!\!\!\!\!\!\!\!
{\rm e}^{i\nu\cdot \psi} H^{(k)}_{\nu} + \!\!\!\!\!\!\!\!\!\!\!
\sum_{\substack{ \nu\in\ZZZ^2_* \\ p(\nu) \ge \gotm k + (k-1)}} \!\!\!\!\!\!\!\!\!\!
{\rm e}^{i\nu\cdot \psi} H^{(k)}_{\nu} .
\end{equation}
For real $\psi$, the first sum is bounded by
\begin{equation} \nonumber
4 \left( \gotm k + k \right)^2N^2
\gotD_0^{\gotm+1} N^{-1}
\gotC_0^k 
|\gotD_0\e|^{\frac{k-1}{4}} ,
\end{equation}
by Lemma \ref{lem:4.8},
while the second sum is bounded by
\begin{equation} \nonumber
\gotD_0^{\gotm+1} N^{-1} 
\gotC_0^{k} 
|\gotD_0\e|^{\frac{k-1}{4}}
\!\!\!\!\!\!\!\!\!\!
\sum_{\substack{ \nu\in\ZZZ^2_* \\ p(\nu) \ge \gotm k + (k-1)}} 
\!\!\!\!\!\!
|\gotD_0 \e|^{\frac{1}{8}(p(\nu)- \gotm k - (k-1))}  ,
\end{equation}
which can be dealt with as in the proof of Lemma \ref{lem:4.6} by noting that
\begin{eqnarray} \nonumber
\sum_{\substack{ \nu\in\ZZZ^2_* \\ p(\nu) \ge p_0 }} x^{p(\nu)-p_0} 
& \!\!\! \le \!\!\! &
\sum_{p=0}^{\io} \sum_{\substack{\nu\in\ZZZ^2_* \\ p(\nu)=p+p_0}} x^{p} \le 
\sum_{p=0}^{\io} x^{p} 4 (p+p_0+1)N^2
\le 4N^2 \left( \frac{1}{(1-x)^2} + \frac{p_0}{1-x} \right)  , 
\nonumber
\end{eqnarray}
with $x=|\gotD_0\e|^{\frac18}$ and $p_0=\gotm k + (k-1)$. Eventually, by bounding 
\begin{equation} \nonumber 
4N^2 \left(  \left( \gotm k + k \right)^2
+ \frac{1}{(1-2^{-1})^2}  + \frac{\gotm k + (k-1)}{1-2^{-1}}\right)  \le \gotC \, k^2 ,
\end{equation}
with 
\begin{equation} \label{eq:4.5}
\gotC := 
4N^2 \left( (\gotm+1)^2 + 
2 + 2(\gotm+1) 
\right) ,
%
\end{equation}
the assertion follows.
%
%
\qed

\begin{lemma} \label{lem:4.10}
Assume Hypothesis \ref{hyp:1},
and assume $\om_0$ to be non-resonant up to order $\gotm N \!$.
There exists a positive constant $\e_3$ such that, for all $|\e|\le \e_3$,
the series expansions \eqref{eq:3.9} converge when setting $\mu=1$.
\end{lemma}

\noindent\emph{Proof}. 
The power series \eqref{eq:3.9} converge for $|\mu| \le \mu_0$, with $\mu_0$ such that
\begin{equation} \nonumber
\gotC_0 |\gotD_0\e|^{\frac{1}{4}} \mu_0 = 
2 \gotD_0^{-1} N \gotB \gotM_0 \gotR_0 |\gotD_0\e|^{\frac{1}{4}} \mu_0 < 1 .
\end{equation}
Therefore, if one chooses $|\e|\le\e_3$, with
\begin{equation} \label{eq:4.6}
\tilde\e_3 := \frac{1}{(2\gotC_0)^4\gotD_0} =
\frac{\gotD_0^3}{(4N \gotB \gotM_0 \gotR_0)^4} , \qquad \e_3 := \min\{ \e_2, \tilde\e_3 \} =
\min\{\tilde\e_2,\tilde\e_3\} ,
\end{equation}
%
one can take $\mu=\mu_0=1$.
\qed

\begin{rmk} \label{rmk:4.11}
\emph{
By using that one has $\gotB \ge C_0^{-1} \gotD_0^2 = 2A_0^{-1}\gotD_0^2$ (cf.~\eqref{eq:4.1}),
$\gotM_0 \ge 6\Phi $ (cf.~\eqref{eq:4.2}) and $\gotR_0 \ge 4N^2(\gotm+1)^2$
(cf.~the proof of Lemma \ref{lem:4.6}), one finds
\begin{equation} \nonumber
4 N \gotB \gotM_0 \gotR_0 \ge \frac{8 N \Phi \gotD_0^2}{A_0} 24 N^2 (\gotm+1)^2 \ge
\frac{24 N^2 (\gotm+1)^2}{\tilde\e_2}
\end{equation}
and hence, by Remark \ref{rmk:4.1},
\begin{equation} \nonumber
\tilde\e_3 \le \frac{\gotD_0^3 \, \tilde\e_2^4}{(24 N^2 (\gotm+1)^2)^4}
\le \frac{1}{8(24 N^2 (\gotm+1)^2)^4} \left( \frac{\tilde\e_2}{\e_1} \right)^3 \tilde \e_2 \le
\frac{1}{8(24 N^2 (\gotm+1)^2)^4} \left( \frac{2}{(16\gotm N)^2} \right)^3 \tilde \e_2 .
\end{equation}
Therefore, in \eqref{eq:4.6}, one has $\e_3 = \tilde\e_3 < \tilde \e_2 < \e_1$.
}
\end{rmk}

\begin{rmk} \label{rmk:4.12}
\emph{
So far, constructing a solution to \eqref{eq:2.1} of the form \eqref{eq:2.7},
 we have fixed the parameter $\al_0$. In the following we need to consider $\al_0$ varying
in a small interval around a fixed value $\al$. Write
$\be_0=\be_0(\al_0)$ and $C_0=C_0(\al_0)$
in \eqref{eq:2.3}, \eqref{eq:3.5} and \eqref{eq:3.8}, respectively,
to stress the dependence on $\al_0$; then, for given $\be>0$, if one defines
\begin{equation} \nonumber
\calA(\be) := \{ \al_0 \in \RRR : 
\be_0(\al_0) \ge \be/2 \} ,
\qquad
C_0^*(\be) := \min \{ C_0(\al_0) : \al_0 \in \calA(\be) \} ,
\end{equation}
%
all the bounds of this and the previous sections still hold for all $\al_0\in\calA(\be)$,
with $C_0$ replaced with $C_0^*(\be)$ and $\be_0$ replaced with $\be/2$.
Indeed, for all $\al_0\in\calA(\be)$, when $\nu\in\ZZZ_*^2$ is such such that $|\om_0\cdot\nu| < \be_0/4$,
one can first bound $|D(\om_0\cdot\nu)|\ge C_0(\al_0)|\e|^2$ -- a property used at  length
in the proof of Lemma \ref{lem:4.2} -- and then $C_0(\al_0)|\e|^2 \ge C_0^*(\be)|\e|^2$,
while all  the remaining factors appearing in $\Val(\vartheta)$ are proportional to inverse powers of $\be_0$
(cf.~Section \ref{sec:3.6}), so that, once $C_0(\al_0)$ has been replaced with $C_0^*(\be)$,
the bounds which hold for $\al_0$ such that $\beta_0(\al_0)=\be/2$ hold as well for all $\al_0$
such that $\be_0(\al_0)\ge \be/2$. This allows us to have uniform bounds for all $\al_0\in \calA(\be)$.
%
%
%
}
\end{rmk}

\zerarcounters 
\section{The implicit function problem}
\label{sec:5} 

By construction the counterterm $\tilde\al$, as well as $H(\psi;\e)$, depends on $\e$ and $\al_0$
(through the frequency vector $\om_0$); thus, we write $\tilde\al=\tilde\al(\e,\al_0)$.

Given $\al$ in \eqref{eq:2.1}, set $\om:=(\al,1)$ and define
\begin{equation} \label{eq:5.1}
\be: = \min\{ |\om\cdot\nu| : \nu \in \ZZZ^2_{*}(4\gotm N)  \} . 
\end{equation}
%

\begin{lemma} \label{lem:5.1}
Assume Hypotheses \ref{hyp:1} and \ref{hyp:2}, and let $\be$ be defined as in \eqref{eq:5.1}.
There exist positive constants $\e_4$, $\eta_0$ and $a_0$,
such that for all $|\e| \le \e_4$ and all $|\al_0- \al| \le \eta_0$ the counterterm
$\tilde \al (\e,\al_0)$ depends continuously on both $\e$ and $\al_0$, and one has $|\tilde \al (\e,\al_0) | \le a_0 |\e|$.
\end{lemma}

\noindent\emph{Proof}. 
For fixed $\al_0$, let $\be_0=\be_0(\al_0)$ be defined as in \eqref{eq:2.3}.
The value $\tilde \e_3$ defined in \eqref{eq:4.6} depends on $\be_0$ through the quantities
$\gotD_0$, $\gotB$ and $\gotM_0$. To make explicit such a dependence we write
$\gotD_0=\gotD_0(\be_0)$, $\gotB=\gotB(\be_0)$ and $\gotM_0=\gotM_0(\be_0)$, and set
\begin{equation} \nonumber
\tilde \e_3(\be_0) :=
\frac{1}{(2\gotC_0)^4\gotD_0(\be_0)} =
\frac{\gotD_0^3(\be_0)}{(4N \gotB \gotM_0(\be_0) \gotR_0)^4} .
\end{equation}
Define also
\begin{equation} \nonumber
\gotB_*(\be_0) := \max\{ 16\be_0^{-2}, (C_0^*(\be))^{-1} \gotD_0^{2}(\be_0)  \} .
\end{equation}
with $C_0^*(\be)$ defined in Remark \ref{rmk:4.12}.
For fixed $\be$, as long as $\be_0 \ge \be/2$, the series in \eqref{eq:3.9} which defines the counterterm
converges provided one has $|\e| \le \e_4$, with
\begin{equation} \label{eq:5.2}
\e_4 = \e_4(\be) := \frac{\gotD_0^3(\be/2)}{(4N \gotB_*(\be/2) \gotM_0(\be/2) \gotR_0)^4} , 
\end{equation}
as discussed in Remark \ref{rmk:4.12}. Note that $\e_4$ differs from $\tilde\e_3(\be_0)$
as, first, $\gotB=\gotB(\be_0)$ is replaced with $\gotB_*(\be/2)$ and, then,
the remaining factors $\be_0$ are replaced with $\be/2$. By construction one has
$\tilde\e_4(\be) \le \min\{ \tilde \e_3(\be_0) : \al_0 \in \calA(\be) \}$.

Set $\gotD_{1}:=\gotD_0(\be/2)$, $\gotM_{1}:=\gotM_0(\be/2)$, $\gotB_{1}:=\gotB_*(\be/2)$.
For all $\al_0 \in \calA(\be)$, with $\calA(\be)$ defined in Remark \ref{rmk:4.12},
by Lemma \ref{lem:4.8} (and Remark \ref{rmk:4.12}) the coefficients $\al^{(k)}$ are bounded as
\begin{equation} \nonumber
| \al^{(1)} | \le \ga^{-1} 6 \Phi \gotD_{1}^{\gotm} |\gotD_{1}\e| , 
\qquad
| \al^{(k)} | \le \ga^{-1} \gotB_1^{-1} \gotD_{1}^{\gotm+1} N^{-1} \, \gotC_{1}^k 
|\gotD_{1}\e|^{\frac{k+2}{4}} ,  \quad k \ge 2 ,
\end{equation}
with $\gotC_{1} := 2 \gotD_{1}^{-1}N \gotM_{1} \gotR_0 \gotB_{1}$.
By Lemma \ref{lem:4.10} and Remark \ref{rmk:4.12},
the series expansion \eqref{eq:3.9} for $\tilde\al$ is well defined and converges for $\mu=1$, provided that $|\e| \le \e_4$.
Indeed, one has
\begin{equation} \nonumber
|\tilde\al| \le \ga^{-1} 6\Phi \gotD_{1}^{\gotm} |\gotD_{1}\e| + 
(\ga \gotB_1 N)^{-1} \gotD_{1}^{\gotm+1} \sum_{k=2}^{\io} \gotC_{1}^{k} |\gotD_{1}\e|^{\frac{k+2}{4}} 
\le
\ga^{-1} \gotD_{1}^{\gotm} ( 6 \Phi +  8 N \gotB_1 \gotD_{1}^{-1} \gotM_{1}^2 \gotR_0^2 ) \, |\gotD_{1}\e| ,
\end{equation}
since $\gotC_{1}|\gotD_{1}\e|^{1/4}<1/2$ if $|\e| \le \e_4$.
Then, defining
\begin{equation} \label{eq:5.3}
a_0 := \ga^{-1} 6 \Phi \gotD_{1}^{\gotm+1} \left(  1 + \frac{4 N \gotB_1 \gotM_{1}^2 \gotR_0^2}{3\Phi \gotD_{1}} \right) ,
\end{equation}
the bound $|\tilde \al | \le a_0|\e|$ follows as long as $\al_0\in\calA(\be)$. 

On the other hand, setting $\om=(\al,1)$ and $\om_0=(\al_0,1)$, with $|\al_0-\al| \le \eta_0$, one has
\begin{equation} \nonumber
|\om_0\cdot\nu| \ge |\om\cdot\nu| - |\al-\al_0| \, |\nu_1| \ge \be - 4 \gotm N \eta_0
\end{equation}
for all $\nu\in\ZZZ^2_*(4\gotm N)$. Thus, defining
\begin{equation} \label{eq:5.4}
\eta_0 := \frac{\be}{8\gotm N} ,
\end{equation}
for all $\al_0 \in I(\eta_0) := [\al-\eta_0,\al+\eta_0]$ one has $\be_0(\al_0) \ge \be/2$ and hence $\al_0\in\calA(\be)$.
Therefore, if one fixes $\al$ so that $\be>0$, for all $\al_0\in I(\eta_0)$ and $\e\in[-\e_4,\e_4]$,
the function $\tilde\al(\e,\al_0)$ depends continuously on $\al_0$ and satisfies the bound
$|\al(\e,\al_0)| \le a_0|\e|$, with $a_0$ as in \eqref{eq:5.3}. Continuity on $\e$ trivially holds for $\e\neq0$
and follows from the bounds of Lemma \ref{lem:4.8} for $\e=0$.
\qed

\begin{rmk} \label{rmk:5.2}
\emph{
The value of $a_0$ given in Lemma \ref{lem:5.1} might be improved by computing explicitly the contribution $\al^{(2)}$,
and using the bounds in Lemma \ref{lem:4.7} for the contributions with $k \ge 3$,
which are bounded proportionally to $|\gotD_0\e|^{5/4}$.
}
\end{rmk}

\begin{lemma} \label{lem:5.3}
Assume Hypotheses \ref{hyp:1} and \ref{hyp:2}, and let $\be$ be defined as in \eqref{eq:5.1}.
Let $\e_4$ and $\eta_0$ be as in Lemma \ref{lem:5.1}.
There exist positive constants $a_1$ and $a_2$ such that, for all $|\e| \le  \e_4$
and all $|\al_0 - \al| \le \eta_0$, one has
%
\begin{equation} \nonumber
\left| \frac{\partial \tilde\al}{\partial \al_0}(\e,\al_0) \right| \le a_1 |\e| , \qquad
\left| \frac{\partial \tilde\al}{\partial \e}(\e,\al_0) \right| \le a_2 .
\end{equation}
\end{lemma}

\noindent\emph{Proof}. 
The counterterm $\tilde\al(\e,\al_0)$ is given by \eqref{eq:3.9}, with the coefficients $\al^{(k)}$ as in Lemma \ref{lem:3.12}.
\vspace{-.3cm}
\begin{itemize}[leftmargin=.5cm]
\itemsep0em
\item[1.] In order to bound the derivative of $\al^{(k)}$ with respect to $\al_0$, one has to bound the derivatives
\vspace{-.1cm}
\begin{equation} \nonumber
\frac{\partial F_{\nu_0}^1 
}{\partial \al_0} , \qquad
\frac{\partial F_{\nu_0}^0(q)}{\partial \al_0} , 	\qquad
\frac{\partial G^{1}_{\nu_0} 
}{\partial \al_0} , \qquad
\frac{\partial \calG_{\ell}}{\partial \al_0} .
\vspace{-.1cm}
\end{equation}
\begin{itemize}[leftmargin=.31cm]
\itemsep0em
\item[1.1.] 
All quantities above depend on the coefficients $h^{(k)}_{\nu}$, which in turn depend on $\al_0$ through the
propagators $\calG_{\ell}^*$ (cf.~Lemma \ref{lem:3.4}).
For any tree $\vartheta\in\gotT_{k,\nu}^{*}$ and any $\ell\in L(\vartheta)$,
\begin{equation} \nonumber
\left| \frac{\partial \calG_{\ell}^*}{\partial \al_0} \right| \le |\calG_{\ell}^*| \frac{2|\nu_1|}{|\om_0\cdot\nu|} \le
2 \be_0^{-1} k N \, |\calG_{\ell}^*| \le 2 \be_0^{-1} \gotm N \, |\calG_{\ell}^*| ,
\end{equation}
so that, by reasoning as in the proof of Lemma \ref{lem:3.5}, one finds $h^{(k)}_{\nu}=0$ if $k\le p(\nu)$ and
\begin{equation} \nonumber
\left| \frac{\partial
h^{(k)}_{\nu}
}{\partial \al_0} \right|
\le 2 \be_0^{-1} k^2 N (\be_0^{-2} \Phi_0)^{k} 
\le 2 \be_0^{-1} \gotm^2 N (\be_0^{-2} \Phi_0)^{k} , 
\qquad
k \ge p(\nu)+1 .
\end{equation}
\item[1.2.] 
With respect to the bound for $F_0^0(1)$ discussed in the proof of Lemma \ref{lem:3.6},
when bounding its derivative with respect to $\al_0$, one has to take into account that, in \eqref{eq:2.10a} one has
\vspace{-.2cm}
\begin{equation} \nonumber
\frac{\partial}{\partial \al_0}
 \left( h^{(k_{0})}_{\tilde\nu_1} \ldots h^{(k_p)}_{\tilde\nu_p} \right)
= \sum_{i=1}^{p} \left( \frac{\partial h^{(k_i)}_{\tilde\nu_i} }{\partial \al_0} \right) 
\prod_{\substack{j=1 \\ j\neq i}}^{p} h^{(k_j)}_{\tilde\nu_j}
\vspace{-.2cm}
\end{equation}
so that one obtains, for $|\e|\le \e_3$,
%
\begin{equation} \nonumber
\begin{aligned}
\left| \frac{\partial F_{0}^0(1)}{\partial \al_0} \right| 
& \le \left| 2 \e \frac{\partial C_0}{\partial \al_0} \right| + \left| \frac{\partial \tilde F^0_0(1)}{\partial \al_0} \right|
\le  4 \be_0^{-1}N C_0 |\e| 
+ 2 \be_0^{-1} \gotm^2 N \sum_{k=2}^{\gotm} | \e|^{k} 4 N^{4} \Phi \left( \be_0^{-2}\Phi_0 \right)^{k}
\\
& + 2 \be_0^{-1} \gotm^2 N \sum_{p=2}^{\io} p N^{p+1} \Phi (4 \gotm^2 N^2)^{p} (2 |\e |\be_0^{-2} \Phi_0)^p \\
& \le 4 \be_0^{-1} N C_0 |\e| + 8 \be_0^{-1} \gotm^2 N \, K_{0} |\e|^2 \le  8 \be_0^{-1} N  C_0 |\e| ,
\end{aligned}
\end{equation}
with $K_{0}$ as in \eqref{eq:3.6}, since one has (cf.~Remark \ref{rmk:4.11})
%
\begin{equation} \label{eq:5.5}
\e_3 < \frac{C_0}{4 \gotm^2 K_{0}} = \frac{1}{4\gotm^2} \tilde\e_2 ,
\end{equation}
with $\tilde\e_2$ as in \eqref{eq:3.8}.
This, together with the bounds in Lemma \ref{lem:3.7}, implies the bound
\begin{equation} \nonumber
\begin{aligned}
\left| \frac{\partial \calG_{\ell}}{\partial \al_0} \right|
& \le
|\calG_{\ell}|^2 \! \left( 
|\nu_{\ell,1}| 
| 2 \om_0\cdot\nu_{\ell} - i \e^{\gotm} \ga |  + \left|  \e \frac{\partial F_0^0(1)}{\partial\al_0} \right| \right)
\\ & 
\le
|\calG_{\ell}| 
\frac{2\,|\om_0\cdot\nu_{\ell}| + |\e|^{\gotm} \ga + 8 \be_0^{-1} N C_0 |\e|^2}{(\om_0\cdot\nu_\ell)^2 + |\e|^2 C_0} 
|\nu_{\ell}| \\
& \le
|\calG_{\ell}| \left( \frac{1}{|\e|\sqrt{C_0}} + \frac{|\e|^{\gotm}\ga}{|\e|^2C_0} + 8\be_0^{-1} N \right) |\nu_{\ell}| 
\le \frac{Q_0}{|\gotD_0\e|} |\calG_{\ell}| \, |\nu_{\ell}| ,
\end{aligned}
\end{equation}
with
\begin{equation} \nonumber
Q_0 := \gotD_0 
\left( \frac{1}{\sqrt{C_0}} + \frac{\ga}{C_0} |\e|^{\gotm-1} + 8 \be_0^{-1} N |\e|  \right).
\end{equation}
The derivatives of $F^0_{\nu_0}(q)$ can be bounded in the same way as was done for $F^0_0(1)$. Again one uses the fact that,
when acting on some $h^{(k_i)}_{\tilde\nu_i}$, the derivative with respect to $\al_0$ produces
a further factor $2\be_0^{-1}\gotm^2 N$, so that one proceeds as in the proof of Lemma \ref{lem:3.9},
by taking into account that, in the sum over $p$, there is an extra factor $2 \be_0^{-1}\gotm^2 N p$.
Eventually, instead of the bound for $F^0_{\nu_0} (q)$ of Lemma \ref{lem:3.9}, one finds,
for all $|\e| \le \e_3$,
\begin{equation} \nonumber
\left| \frac{\partial F_{\nu_0}^0(q)}{\partial \al_0} \right| \le 
8 \be_0^{-1} N \gotm^2 \, 3 N^{q} \Phi |\gotD_0\e|^{p(\nu_0)} .
\end{equation}
The same argument applies to $F^1_{\nu_0}(0)$ as well. With respect to the bounds in Lemma \ref{lem:3.11},
the derivatives are bounded as
\begin{eqnarray} \nonumber
\left| \frac{\partial
F_{\nu_0}^1 
}{\partial\al_0} \right|
& \!\!\! \le \!\!\! &
8 \be_0^{-1} N \gotm^2 \, 6 \Phi \gotD_0^{\gotm}
|\gotD_0\e|^{p_{\gotm}(\nu_0)}  , 
\qquad \nu_0 \in 
\ZZZ_*^2 ,
\nonumber \\
\left| \frac{\partial
F_{0}^1 
}{\partial\al_0} \right|
& \!\!\! \le \!\!\! &
8 \be_0^{-1} N \gotm^2 \, 6 \Phi \gotD_0^{\gotm} \gotD_0|\e| .
\nonumber
\end{eqnarray}
Finally, the derivative of $G^0_{\nu}(0)$ is bounded as
\begin{equation} \nonumber
\left| \frac{\partial
G^{1}_{\nu_0} 
}{\partial\al_0} \right| 
\le 2\be_0^{-1} \gotm^2 N \, 2\ga \be_0^{-1} \Phi_0 .
\end{equation}
\item[1.3.] When evaluating the derivative with respect to $\al_0$ of $\Val(\vartheta)$, for any tree $\vartheta\in\gotT_{k,0}$,
one obtains $2k-1$ contributions:
(a) $k-1$ of them are produced when the derivative acts on the propagator $\calG_{\ell}$ of a line different from
the root line, and for each of them there is an extra factor 
$Q_0|\nu_{\ell}|/|\gotD_0\e|$
with respect to the bound for  $\Val(\vartheta)$;
(b) the other $k$ are produced when the derivative acts on a node factor, and each of them admits
the same bound as before times a factor which is at most $48 \be_0^{-1}\gotm^2 N$.
Moreover, for any line $\ell\in L(\vartheta)$, one has
\begin{equation} \nonumber
|\nu_{\ell}| \le \NN(\vartheta) , \qquad \NN(\vartheta) := \sum_{v \in N(\vartheta)}  |\nu_v| .
\end{equation}
%
%
If one writes in the bound of $\Val(\vartheta)$
\begin{equation} \nonumber
\frac{1}{4} p(\nu_v) = \frac{1}{8} p(\nu_v) + \frac{1}{8} p(\nu_v) ,
\qquad 
\frac{1}{4} p_{\gotm}(\nu_v) = \frac{1}{8} p_{\gotm}(\nu_v) + \frac{1}{8} p_{\gotm}(\nu_v) ,
\end{equation}
and uses the fact that, by Lemmas \ref{lem:3.1} and \ref{lem:3.3},
\begin{equation} \nonumber
\NN(\vartheta) \! \le \! (\calP(\vartheta) + \gotm |E(\theta)| + k) N \le (\calP(\vartheta) + (\gotm+1) k) N , \quad
\calP(\vartheta) \! := \!\!\!\! \sum_{v \in I(\vartheta)} \!\! p(\nu_v) + \!\!\!\! \sum_{v \in E(\vartheta)} \!\! p_{\gotm}(\nu_v) ,
\end{equation}
one can bound
\begin{equation} \nonumber
\NN(\vartheta) \Bigl( \prod_{v \in I(\vartheta)} \!\!\! |\gotD_0\e|^{\frac{1}{8} p(\nu_v)} \Bigr)
\Bigl( \prod_{v \in E(\vartheta)} \!\!\! |\gotD_0\e|^{\frac{1}{8} p_{\gotm}(\nu_v)} \Bigr) 
\le (\gotm+2)Nk .
\end{equation}
\item[1.4.] 
The derivative of $\al^{(1)}=\ga^{-1}F^1_0(0)$ can be computed explicitly, and gives
\begin{equation} \nonumber
\left| \frac{\partial
\al^{(1)} 
}{\partial\al_0} \right|
\le
8 \be_0^{-1} N \gotm^2 \, 6 \ga^{-1} \Phi \gotD_0^{\gotm} \, |\gotD_0\e| .
\end{equation}
\item[1.5.] 
When considering the coefficients $\al^{(k)}$, for $k\ge 2$, first of all one notes that,
for $|\e| \le \e_3$, by using also Remark \ref{rmk:4.11} and \eqref{eq:3.5}, one has
\begin{equation} \nonumber
48\be_0^{-1} \gotm^2 N |\gotD_0 \e| \le
48\be_0^{-1} \gotm^2 N \gotD_0 \e_1 =
\frac{3\be_0\gotD_0}{16 N^4\Phi} \le
\frac{3\sqrt{2}}{16N^2} \frac{\gotD_0}{\sqrt{C_0}} \le  Q_0
\end{equation}
and hence $48\be_0^{-1}\gotm^2 N \le Q_0/|\gotD_0\e|$.
\begin{itemize}[leftmargin=.3cm]
\item[1.5.1.] 
For $k\ge 3$, it is convenient to use the improved bound in Remark \ref{rmk:4.4}, which gives
\begin{equation} \nonumber
\left| \frac{\partial
\al^{(k)} 
}{\partial\al_0} \right|
\le
(\gotm+2) N k (2k-1) \ga^{-1} Q_0 \gotB^{-1}\gotD_0^{\gotm+1}N^{-1} \gotC_0^k |\gotD_0\e|^{\frac{k+1}{4}} , \qquad k \ge 3 .
\end{equation}
\item[1.5.2.] 
The contribution with $k=2$ needs a more careful analysis with respect to the discussion in Lemma \ref{lem:4.3}
(cf.~also Remark \ref{rmk:5.2}).
The two nodes of any tree $\vartheta$ whose value contributes to $\al^{(2)}$ have opposite mode labels
$-\nu_{\ell}$ and $\nu_{\ell}$ (since the sum vanishes), where $\nu_\ell$ is the momentum of the line exiting the end node, so that
$|\om_0\cdot\nu_\ell| <\be_0/4$ is possible only if $p(\nu_\ell)\ge 4\gotm$ and in that case $\calW(\vartheta)$
contains an extra factor $|\gotD_0\e|^{\frac34p(\nu_\ell)} \le |\gotD_0\e|^{3\gotm} \le |\gotD_0\e|^{3}$.
If $|\om_0\cdot\nu_{\ell}|\ge \be_0/4$, then the derivative of the propagator
admits a better bound, i.e.
\begin{equation} \nonumber
\begin{aligned}
\left| \frac{\partial \calG_{\ell}}{\partial \al_0} \right|
& \le
|\calG_{\ell}| 
\frac{2 |\om_0\cdot\nu_{\ell}| + |\e|^{\gotm} \ga + 8 \be_0^{-1} N C_0 |\e|^2}{(\om_0\cdot\nu_\ell)^2 + |\e|^2 C_0} 
|\nu_{\ell}| \\
& \le
8 |\calG_{\ell}| \left( \be_0^{-1} + 2\be_0^{-2}|\e|^{\gotm}\ga + \be_0^{-1} N \right) |\nu_{\ell}| ,
\end{aligned}
\end{equation}
where one easily checks that $8 (\be_0^{-1} +2\be_0^{-2}|\e|^{\gotm}\ga + \be_0^{-1} N ) \le Q_0$.
This means that, with respect to $\calW(\vartheta)$, the derivative with respect to $\al_0$ produces an extra factor
which is less than $8 \be_0^{-1} N \gotm^2$, when acting on a node factor, and less than
$\max\{ Q_0|\gotD_0\e|^2 , Q_0 \} = Q_0$, when acting on the propagator of the line exiting the end node.
In conclusion, one obtains
\begin{equation} \nonumber
\left| \frac{\partial
\al^{(2)} 
}{\partial\al_0} \right|
\le
(\gotm+2) N 6 \ga^{-1} Q_0 \gotB^{-1} \gotD_0^{\gotm+1}N^{-1} \gotC^2 |\gotD_0\e| ,
\end{equation}
which, together with the previous bounds, imply
\begin{equation} \nonumber 
\hspace{-.8cm}
\left| \frac{\partial
\tilde\al 
}{\partial\al_0} \right|
\le 
\ga^{-1} \gotD_0^{\gotm} \left( \frac{48 \gotm^2 N \Phi}{\be_0} +
\frac{48 (\gotm+1) Q_0 N^2 \gotB \gotM_0^2 \gotR_0^2}{\gotD_0}
\left( 1 + \frac{16 N \gotB \gotM_0 \gotR_0}{3 \gotD_0}  \right) \right) |\gotD_0\e| .
\end{equation}
\end{itemize}
\item[1.6.] 
Finally, in order to obtain uniform bounds for $|\al_0-\al| < \eta_0$,
one has to reason as in the proof of Lemma \ref{lem:5.1}: one requires $|\e| \le \e_4$ and one has to replace
$\be_0^{-1}$, $\gotB$, $\gotD_0$ and $\gotM_0$ with, respectively,
$2\be^{-1}$, $\gotB_1$, $\gotD_{1}$ and $\gotM_{1}$, as defined in the proof of Lemma \ref{lem:5.1}, and $Q_0$ with
\begin{equation} \nonumber
Q_1 := 2 \gotD_{1} \left( \frac{1}{2 \sqrt{C_0^*(\be)}} + \frac{\ga}{2 C_0^*(\be)} |\e|^{\gotm-1} + 16 \be^{-2} N |\e|  \right) .
\end{equation}
Then, provided one has $|\al_0-\al| < \eta_0$ and $|\e| \le \e_4$, the first bound follows with 
\begin{equation} \label{eq:5.6}
\hspace{-.8cm}
a_1 := 
\ga^{-1} \gotD_{1}^{\gotm+1} \left( \frac{96 \gotm^2 N \Phi}{\be} +
\frac{48 (\gotm+1) Q_1 N^2 \gotB_1 \gotM_{1}^2 \gotR_0^2}{\gotD_{1}^2} \left( 
1 + \frac{16 N \gotB_1 \gotM_{1} \gotR_0}{3\gotD_{1}} \right) \right) .
\end{equation}
\end{itemize}
\item[2.] 
The derivative with respect to $\e$ of $\al^{(k)}$ can be bounded in a similar way.
\vspace{-.2cm}
\begin{itemize}[leftmargin=.31cm]
\item[2.1.]
First of all we recall that $h^{(k)}_{\nu}=0$ if $k\le p(\nu)$ and observe that
\begin{equation} \nonumber 
\left| \frac{\partial
\, \e^k h^{(k)}_{\nu}
}{\partial \e} \right|
\le \frac{k}{|\e|} (|\e|\be_0^{-2} \Phi_0)^{k} \le \frac{\gotm}{|\e|} (|\e|\be_0^{-2} \Phi_0)^{k} , 
\qquad
k \ge p(\nu) + 1 ,
\end{equation}
as Lemma \ref{lem:3.5} and the definition of the function $h$ immediately imply.
Therefore, in order to bound the derivatives with respect to $\e$ of the coefficients
$F_{\nu}^0(q)$, $F_{0}^1(0)$, $F_{\nu}^1(0)$, one observes that,
with respect to the bounds in the proofs of Lemmas \ref{lem:3.9} and \ref{lem:3.11},
a further factor $p\gotm/|\e|$ appears for each summand in the sum over $p$,
so that eventually, one obtains
\begin{eqnarray} \nonumber
\left| \frac{\partial F_{\nu_0}^0(q)}{\partial \e} \right| 
& \!\!\! \le \!\!\! &
\frac{4\gotm}{|\e|}
3 N^{q} \Phi
|\gotD_0\e|^{p(\nu_0)} , \\
%
\left| \frac{\partial
F_{\nu_0}^1 
}{\partial\e} \right|
& \!\!\! \le \!\!\! &
\frac{4 \gotm}{|\e|}
6 \Phi \gotD_0^{\gotm} |\gotD_0\e|^{p_{\gotm}(\nu_0)}  , 
\quad \nu_0 \in \ZZZ^2_* 
, \qquad \qquad
\left| \frac{\partial
F_{0}^1 
}{\partial\e} \right|
\le
\frac{4\gotm}{|\e|} 
6 \Phi |\gotD_0|^{\gotm} |\gotD_0\e| . 
\nonumber
\end{eqnarray}
In a similar way one obtains, for $|\e| \le \e_3$, using also \eqref{eq:5.5},
\begin{equation} \nonumber
\left| \frac{\partial
G_{\nu_0}^1 
}{\partial\e} \right| 
\le \frac{2\gotm}{|\e|} 
\ga \be_0^{-1} \Phi_0 ,  
%
%
\qquad
\left| \frac{\partial F_{0}^0(1)}{\partial \e} \right| 
\le
2 C_0 +  4 \gotm  K_{0} |\e| \le 4 C_0 ,
\end{equation}
with $K_{0}$ as in \eqref{eq:3.6}. The latter bound also implies
\begin{eqnarray} \nonumber
\left| \frac{\partial \calG_{\ell}}{\partial \e} \right|
& \!\!\!\! \le \!\!\!\! &
|\calG_{\ell}|^2
\left| -i \gotm \e^{\gotm-1} \ga \, \om_0\cdot\nu_{\ell}  - \e \frac{\partial F_0^0(1)}{\partial\e} - F_0^0(1) \right|
\\
& \!\!\!\! \le \!\!\!\! &
 |\calG_{\ell}| 
\frac{ \gotm |\e|^{\gotm-1} \ga |\om_0\cdot\nu_\ell |+ 4 C_0 |\e| + 4 C_0|\e| }{\max
\{ (\om_0\cdot\nu_\ell)^2 + |\e|^2 C_0, |\e|^{\gotm} \ga |\om_0\cdot\nu_\ell| \}}
\le \frac{2\gotm + 8
}{|\e|} |\calG_{\ell}| ,
\nonumber
\end{eqnarray}
%
%
%
\item[2.2.]
Thus, the derivative with respect to $\e$ of $\Val(\vartheta)$, for any tree $\vartheta\in\gotT_{k,0}$,
one obtains $2k-1$ contributions:
$k-1$ of them are produced when the derivative acts on a propagator, and for each of them there
is an extra factor $\tilde Q_2/|\e|$ with respect to the bound for  $\Val(\vartheta)$;
the other $k$ are produced when the derivative acts on a node factor, and each of them admits
the same bound as before times a factor which is bounded by
\begin{equation} \nonumber
\gotD_0 \left( \frac{4\gotm}{|\gotD_0\e|} + \frac{(\gotm +1)q-\gotm}{|\gotD_0\e|} \right) =
\frac{\gotD_0 \left( (\gotm +1 ) q  + 3 \gotm \right)}{|\gotD_0\e|} ,
\end{equation}
where the factor $q$ cancels out thanks to the factorial $q!$ appearing in the node factor $\calI_v$.
If $|\e| < \e_3$, then, 
for any $\vartheta\in\gotT_{k,0}$, the derivative 
with respect to $\e$ of $\Val(\vartheta)$
admits the same bound as found for $\Val(\vartheta)$ in Lemma \ref{lem:4.8}
times an extra factor $(2k-1)Q_2/|\e|$, with
\vspace{-.1cm}
\begin{equation} \nonumber
Q_2 
= \max\{2\gotm + 8 , 4\gotm+1\} .
\vspace{-.2cm}
\end{equation}
Therefore one obtains
\begin{equation} \nonumber 
\left| \frac{\partial
\tilde\al 
}{\partial\e} \right|
\le 4 \gotm \, 6 \ga^{-1} \Phi \gotD_0^{\gotm+1} + 24 \,  Q_2  N \ga^{-1} \gotD_0^{\gotm}
\gotB \gotM_0^2 \gotR_0^2 .
\vspace{-.1cm}
\end{equation}
\item[2.3.]
Uniform bounds for $|\al_0-\al|\le\eta_0$ are given by replacing $\gotB$, $\gotD_0$ and $\gotM_0$ with $\gotB_1$,
$\gotD_{1}$ and $\gotM_{1}$, respectively. 
%
%
In conclusion, the second bound follows as well, with
\vspace{-.1cm}
\begin{equation} \label{eq:5.7}
a_2 := 24 \ga^{-1} \gotm \Phi \gotD_{1}^{\gotm+1} 
\left( 1 + \frac{Q_2 N \gotB_1 \gotM_{1}^2 \gotR_0^2}{\gotm \Phi \gotD_{1}} \right) ,
\vspace{-.2cm}
\end{equation}
provided one has $|\al_0-\al| \le \eta_0$ and $|\e| \le \e_4$.
\end{itemize}
\vspace{-.3cm}
\end{itemize}
Therefore the bounds follow for all $|\al_0 -\al | \le \eta_0$ and all $|\e| \le \e_4$. \qed

\begin{rmk} \label{rmk:5.4}
\emph{
A more careful analysis of the coefficient $\al^{(3)}$ would permit an improvement to the bound of the derivative
with respect to $\al_0$. Indeed the derivative of $\al^{(3)}$  turns out to be more than linear, 
and the derivatives of the coefficients $\al^{(k)}$, with $k\ge 4$, are bounded proportionally to $|\gotD_0\e|^{5/4}$,
as it follows from Remark \ref{rmk:4.4} and from the fact that the derivative with respect to $\al_0$
produces at most the loss of a power $\e$.
}
\end{rmk}

\begin{lemma} \label{lem:5.5}
Assume Hypotheses \ref{hyp:1} and \ref{hyp:2}, and let $\be$ be defined as in \eqref{eq:5.1}.
Then there exist constants $\al_0$ and $\e_0$ such that, setting
\vspace{-.1cm}
\begin{equation} \nonumber
\om_0 :=(\al_0,1) , \qquad  \be_0 : = \min\{ |\om_0\cdot\nu| : \nu \in \ZZZ^2_{*}(4\gotm N)  \} ,
\vspace{-.1cm}
\end{equation}
for all $|\e|\le \e_0$ there exists a counterterm $\tilde\al(\e,\al_0)$, differentiable in both $\e$ and $\al_0$,
and a function $H(\psi;\e)$, continuous in $\e$ and analytical in $\psi$, such that
$\al_0 + \e \tilde \al(\e,\al_0)=\al$, $\be_0 \ge \be/2$ and $H(\om_0 t;\e)$ solves \eqref{eq:2.8}.
%
%
\end{lemma}

\noindent\emph{Proof}. 
If $\al$ is such that $\be$ is positive, the counterterm $\tilde\al(\e,\al_0)$ is well defined and differentiable
as long as $|\e| \le \e_4$ and $|\al_0-\al| \le \eta_0$, by Lemma \ref{lem:5.3}. Define
\vspace{-.1cm}
\begin{equation} \label{eq:5.8}
\e_0 = \min\{\e_4,\tilde\e_5\} , \qquad \tilde\e_5 := \min\left\{  \sqrt{\frac{1}{2a_1}} , \sqrt{\frac{\eta_0}{4(a_0+a_2)}} \right\} .
\vspace{-.1cm}
\end{equation}
For fixed $\al$, one needs $\al_0$ to be such $\al_0+\e \tilde\al(\e,\al_0)=\al$.
Consider the implicit function problem
\vspace{-.1cm}
\begin{equation} \nonumber
F(\e,\al_0) = 0 , \qquad
F(\e,\al_0)  := \al_0+\e \tilde\al(\e,\al_0) - \al .
\vspace{-.1cm}
\end{equation}
One has $F(0,\al)=0$ and, if one takes $|\e| \le \e_0$, with $\e_0$ as in \eqref{eq:5.8},
\vspace{-.1cm}
\begin{equation} \nonumber
\begin{aligned}
\frac{\partial F}{\partial \al_0}(\e,\al_0) & = 1 + \e \frac{\partial\tilde\al}{\partial\al_0}(\e,\al_0)
\ge 1 - a_1 |\e|^2 \ge \frac{1}{2} ,  \\
\left| \frac{\partial F}{\partial \e}(\e,\al_0) \right| & = \left| 
\tilde\al(\e,\al_0) + 
\e \frac{\partial\tilde\al}{\partial\e}(\e,\al_0) \right| \le (a_0 + a_2) |\e| ,
\end{aligned}
\end{equation}
with $a_0$, $a_1$ and $a_2$ given by \eqref{eq:5.3}, \eqref{eq:5.6} and \eqref{eq:5.7}, respectively.

One has, for $|\e| \le \e_0$ and
for suitable constants $\al_*\in(\al,\al+\eta_0)$ and $\e_*\in(0,\e)$,
\begin{equation} \nonumber
\begin{aligned}
F(\e,\al + \eta_0) & = F(\e,\al +\eta_0) - F(0,\al) = F(\e,\al + \eta_0) - F(\e,\al) + F(\e,\al) - F(0,\al) \\
& = \eta_0 \, \frac{\partial F}{\partial \al_0} (\e,\al_*) +  \e \, \frac{\partial F}{\partial \e} (\e_*,\al) \ge
\frac{\eta_0}{2} - (a_0+a_2) \e_0^2 \ge \frac{\eta_0}{4} ,
\end{aligned}
\end{equation}
by \eqref{eq:5.8}.
%
%
In a similar way one proves that $F(\al-\eta_0,\e) \le -\eta_0/4$ for $|\e| \le \e_0$.
Therefore, one concludes that for any
$|\e|\le \e_0$, there exists a unique value $\al_0=\al_0(\e)$ such that $F(\e,\al_0(\e))=0$.

The function $H(\psi;\e)$, given by \eqref{eq:3.9} with $\mu=1$,
for all $|\e|\le \e_0$ and $|\al_0-\al|\le \eta_0$, is well defined
and continuous in $\al_0$ as well as in $\e$, by construction, and,
since $\al_0+\e\tilde\al=\al$, it solves \eqref{eq:2.8}.
Lemma \ref{lem:4.9}, together with Lemma \ref{lem:4.10}, only proves that $H(\psi;\e)$,
as a function of $\psi$, is bounded. To prove analyticity, one relies once more on the bound
\begin{equation} \nonumber
| H^{(k)}_{\nu} | \le 
\gotD_1^{\gotm} N^{-1} \gotC_1^{k} 
|\gotD_1 \e|^{\frac{k-1}{4} +
\max\{ \frac{1}{8}(p(\nu)-\gotm k -(k-1)) , 0 \} }  , 
\end{equation}
which follows, for all $k\in \NNN$ and all $\nu\in\ZZZ_*^2$, from Lemma \ref{lem:4.8} and Remark \ref{rmk:4.12}.
Write the Fourier coefficients of $H(\psi;\e)$ as
\begin{equation} \nonumber
H_{\nu} =  \sum_{k=1}^{\io} H^{(k)}_{\nu} ,
\end{equation}
and, for any $\nu\in\ZZZ_*^2$, define
\begin{equation} \nonumber
k(\nu) 
:= \left\lfloor \frac{p(\nu)+1}{\gotm+1} \right\rfloor .
\end{equation}
%
Then one has, if $k(\nu) \ge 1$,
\begin{eqnarray} \nonumber
| H_{\nu} | 
&\!\!\!\! \le \!\!\!\! &
\sum_{k=1}^{k(\nu)} H^{(k)}_{\nu} + \sum_{k=k(\nu)+1}^{\io} H^{(k)}_{\nu} \\ [-1ex]
&\!\!\!\! \le \!\!\!\! &
\gotD_1^{\gotm}N^{-1}
\sum_{k=1}^{k(\nu)} \gotC_1^k |\gotD_1\e|^{\frac{k-1}{4} + \frac{\gotm+1}{8}(k(\nu)-k)}
+ \gotD_1^{\gotm}N^{-1} \!\!\!\!\! \sum_{k=k(\nu)+1}^{\io} \gotC_1^k |\gotD_1\e|^{\frac{k-1}{4} } \nonumber \\ [-1.6ex]
&\!\!\!\! \le \!\!\!\! &
2 \gotD_1^{\gotm}N^{-1} |\gotD_1\e|^{\frac{\gotm+1}{8}k(\nu) - \frac{1}{4}}
\left( \gotC_1 |\gotD_1\e|^{ - \frac{1}{8}(\gotm-1)} \right)^{k(\nu)}
+ 2\gotD_1^{\gotm}N^{-1} \gotC_1 \left( \gotC_1 |\gotD_1\e|^{\frac{1}{4} } \right)^{k(\nu)} \nonumber \\
&\!\!\!\! \le \!\!\!\! &
2 \gotD_1^{\gotm}N^{-1} \gotC_1^{k(\nu)} |\gotD_1\e|^{\frac{k(\nu) - 1}{4}} 
+ 2 \gotD_1^{\gotm}N^{-1} \gotC_1^{k(\nu)+1} |\gotD_1\e|^{\frac{k(\nu)}{4}} \nonumber \\
&\!\!\!\! \le \!\!\!\! &
2 \gotD_1^{\gotm}N^{-1} \gotC_1^{k(\nu)} |\gotD_1\e|^{\frac{k(\nu) - 1}{4}}
\left( 1 + \gotC_1 |\gotD_1\e|^{\frac{1}{4}} \right) \nonumber \\
&\!\!\!\! \le \!\!\!\! &
4 \gotD_1^{\gotm}N^{-1} \gotC_1^{k(\nu)} |\gotD_1\e|^{\frac{k(\nu) - 1}{4}} , \nonumber
\end{eqnarray}
which implies, for $|\nu| \ge 2N(\gotm+1)$, that
\begin{equation} \nonumber
| H_{\nu} | \le 4 \gotD_1^{\gotm}N^{-1} \gotC_1 ( \gotC_1^4|\gotD_1\e|)^{\frac{|\nu|}{2N(\gotm+1)}}.
\end{equation}
For fixed $\e$, the Fourier coefficients decay, for $|\nu|$ large enough, as
\begin{equation} \nonumber
| H_{\nu} | \le 4 \gotD_1^{\gotm}N^{-1} {\rm e}^{-\xi_{\e} |\nu|} , \qquad
\xi_{\e} := \frac{1}{2N(\gotm+1)} \log \frac{1}{\gotC_1^4|\gotD_1\e|} ,
\end{equation}
so that, one finds that, uniformly for all $|\e| \le \e_0$, the function $H(\psi;\e)$ is analytic in $\psi$
in a strip of width
\begin{equation} \nonumber
\xi := \frac{2 \log 2}{N(\gotm+1)} ,
\end{equation}
since $\gotC_1^4|\gotD_1\e| \le 1/16$. Then the assertion follows. 
\qed

\begin{rmk} \label{rmk:5.7}
\emph{
%
%
According to \eqref{eq:5.8}, $\e_0=\min\{\e_4,\tilde\e_5\}$.
For fixed $\e$ and $\gotm$, if $\ga$ tends to $0$,
$\e_4$ tends to a finite value, while $\tilde\e_5$ tends to 0.
}
\end{rmk}

\begin{lemma} \label{lem:5.6}
Assume Hypotheses \ref{hyp:1} and \ref{hyp:2}, and let $\be$ be defined as in \eqref{eq:5.1}.
Let $\e_0$ be defined as in Remark \ref{rmk:5.7}.
Then for all $|\e| \le \e_0$ there is a multi-periodic solution to \eqref{eq:2.1},
with $f$ given by \eqref{eq:2.2}, 
of the form $\theta(t)=\theta_0+\al_0 t + \Theta(\om_0 t;\e)$, with $\theta_0\in\TTT^2$ and $\om_0=(\al_0,1)$, such that 
\begin{equation} \nonumber
\lim_{\e\to0} \al_0 = \al , \qquad
\lim_{\e\to0} \Theta(\psi ; \e) = 0 .
\end{equation}
\end{lemma}

\noindent\emph{Proof}. 
The result follows from Lemma \ref{lem:5.5}, fixing $\tilde\al$ accordingly 
and using that $h(\om_0 t;\e)$ and $H(\om_0 t;\e)$ solve \eqref{eq:2.4} and \eqref{eq:2.8}, respectively.
\qed

\begin{rmk} \label{rmk:5.8}
\emph{
Assuming Hypothesis \ref{hyp:2}, one has $\be >0$ and hence Lemma \ref{lem:5.6} applies. However,
the two frequencies $\al$ and $1$ are either commensurate or incommensurate according to the value of $\e$.
Therefore, for a given value of $\e$ we cannot exclude the possibility that the solution is periodic. Of course,
for a full measure set of values of $\e$ in $[-\e_0,\e_0]$ the solution is expected to be quasi-periodic.
}
\end{rmk}

\zerarcounters 
\section{Proof of Theorem \ref{thm:2}}
\label{sec:6} 

The value $\e_0$ in Lemma \ref{lem:5.6} depends on $\ga$ through the constants $\gotM_1=\gotM_0(\be/2)$,
appearing in $\tilde\e_4$ in \eqref{eq:5.2}, and the constants $a_0$, $a_1$ and $a_2$ -- which in turn depend on
$\gotM_1$ too --  appearing in $\tilde\e_5$ in \eqref{eq:5.8}.
In particular, for $\ga$ large enough, $\gotM_1$ is proportional to $\ga$ and hence $\tilde\e_4$ is proportional to $\ga^{-4}$.
However the bounds in Lemma \ref{lem:4.8} may be improved with a more careful analysis, showing that,
with the exception of the tree of order 1,
each time an end node $v$ carries a label $\lambda_{v}=G$ a further power of $\e$ arises.
This will allow us to write $\ga=|\e|^{-a} \bar \ga$, with $a\in [0,1)$ and $\bar \ga$ independent of $\e$,
so that Theorem \ref{thm:2} follows as soon as one fixes $a$ so that $\gotm: = s + a$ is an integer.

First of all the bound in Lemmas \ref{lem:4.2} and \ref{lem:4.3} can be refined as follows.

\begin{lemma} \label{lem:6.1}
Assume Hypotheses \ref{hyp:1}, 
and assume $\om_0$ to be non-resonant up to order $\gotm N \!$.
Let $\gotD_0$ and $\e_2$ be defined as in Lemma \ref{lem:3.9} and in Lemma \ref{lem:3.7}, respectively,
and assume that $|\e| \le \e_2$.
For all $\vartheta\in\gotT_{k,\nu}$, with $k \ge 1$ and $\nu\in\ZZZ^2_*$ and $|E_G(\vartheta)|\ge 1$, one has
\begin{equation} \nonumber
\calW(\vartheta) \le 
\begin{cases}
\gotB^{k} |\gotD_0\e|^{\frac{k-1}{4}} |\gotD_0\e|^{|E_G(\theta)|} , & \quad k \ge 3 , \; |E_G(\vartheta) | \ge 1 , \\
\gotB^{2} |\gotD_0\e| ,  & \quad k =2 , \; |E_G(\vartheta)| = 1 ,\\
\gotB & \quad k = 1 , \; |E_G(\vartheta)| = 1 . \\
\end{cases}
\end{equation}
\end{lemma}

\noindent\emph{Proof}. 
The bounds are checked by direct computation for $k=1,2,3$, and are proved by induction for $k\ge 4$.
If $v_0$ is the node which the root line $\ell_0$ of the tree $\vartheta$ exits,
let $\ell_1,\ldots,\ell_q$ denote the $q$ lines entering $v_0$, with $q\ge 0$, and
let $\vartheta_1,\ldots,\vartheta_q$ be the subtrees whose root lines are $\ell_1,\ldots,\ell_q$, respectively.
Up to a reordering, one can assume that 
the first $q_1$ of such subtrees contain either at least two end nodes with end node labels $G$ or none;
the following $q_2$ subtrees contain only one node with end node label $G$ and
are of order 2;
the last $q_3$ subtrees contain only one node with end node label $G$ and are of order 1.
By construction, one has $q = q_1 + q_2 + q_3$ and
$|E_G(\vartheta)| = q_2 + q_3$ if $q_1=0$, otherwise, if $q_1\ge 1$,
\begin{equation} \label{eq:6.1}
|E_G(\vartheta)| = \sum_{i=1}^{q_1} |E_G(\vartheta_i)| + q_2 + q_3 .
\end{equation}
\vspace{-.6cm}
\begin{itemize}[leftmargin=.5cm]
\itemsep0em
\item[1.]
For $k=1$ 
one has $\calW(\vartheta) \le 16\be_0^{-2} \le \gotB$, 
since $|\nu_{\ell_0}| \le \gotm N$.
\item[2.]
For $k=2$, one has $q=1$ and the set $N(\vartheta)$ contains the nodes $v_0$ and $v_1$,
where $v_1$ is the end node which $\ell_1$ exits, so that $\lambda_{v_1}=G$.
Then $|\nu_{\ell_1}| \le \gotm N$ and hence $|\om_0\cdot\nu_{\ell_1}| \ge \be_0$. Therefore,
if $|\om_0\cdot\nu_{\ell_0}| \ge \be_0/4$, one finds
\begin{equation} \nonumber
\calW(\vartheta) \le 16\be_0^{-2} |\gotD_0\e|^{1+\frac{3}{4}p(\nu_{v_0})} 16\be_0^{-2} \le \gotB^2 |\gotD_0\e| ,
\end{equation}
while, if $|\om_0\cdot\nu_{\ell_0}| < \be_0/4$, one obtains
\begin{equation} \nonumber
\calW(\vartheta) \le C_0^{-1} \gotD_0^2 |\gotD_0\e|^{-2 + 1+\frac{3}{4}p(\nu_{v_0})} 16\be_0^{-2} 
\le \gotB^2 |\gotD_0\e|^{\frac{5}{4}} ,
\end{equation}
since $p(\nu_{v_0}) \ge 3\gotm$ in such a case. Therefore one can bound 
$\calW(\vartheta) \le \gotB^2 |\gotD_0\e|$.
\item[3.]
For $k = 3$, one has either $q=1$ or $q=2$. 
\vspace{-.2cm}
\begin{itemize}[leftmargin=.31cm]
\itemsep0em
\item[3.1.] 
If $q=2$, the set $N(\vartheta)$ contains $v_0$ and the two end nodes $v_1$ and $v_2$
which the two lines $\ell_1$ and $\ell_2$ exit, respectively, such that $\lambda_{v_2}=G$.
\begin{itemize}[leftmargin=.3cm]
\itemsep0em
\item[3.1.1.] 
If $\lambda_{v_1}=F$ and $\lambda_{v_2}=G$, one has
\begin{equation} \nonumber
\calW(\vartheta) \le |\calG_{\ell_0}| |\gotD_0\e|^{2(\gotm+1)-\gotm+\frac{3}{4} (p(\nu_{v_0}) + p_{\gotm}(\nu_{v_1}))}
|\calG_{\ell_1}| \gotB .
\end{equation}
We distinguish among four possible subcases.
\begin{itemize}[leftmargin=.31cm]
\itemsep0em
\item[3.1.1.1.] 
If $|\om_0\cdot\nu_{\ell_0}|, |\om_0\cdot\nu_{\ell_1}| \ge \be_0/4$ one has
\begin{equation} \nonumber
\calW(\vartheta) \le (16\be_0^{-2})^2 \gotB |\gotD_0\e|^{2(\gotm+1)-\gotm}
\le \gotB^3 |\gotD_0\e|^{\gotm+2} \le \gotB^3 |\gotD_0\e|^3 .
\end{equation}
\item[3.1.1.2.] If $|\om_0\cdot\nu_{\ell_0}| \ge \be_0/4$ and $|\om_0\cdot\nu_{\ell_1}| < \be_0/4$, one has
\begin{equation} \nonumber
\calW(\vartheta) \le C_0^{-1} \gotD_0^2 \,16\be_0^{-2} \gotB |\gotD_0\e|^{-2 + 2(\gotm+1)-\gotm + \frac{9}{4}\gotm} \le
\gotB^3 |\gotD_0\e|^{\frac{13}{4}} ,
\end{equation}
since $|\nu_{\ell_1}| > 4\gotm N$ implies $p_{\gotm}(\nu_{v_1}) \ge 3\gotm$.
\item[3.1.1.3.] 
If $|\om_0\cdot\nu_{\ell_0}| < \be_0/4$ and $|\om_0\cdot\nu_{\ell_1}| \ge \be_0/4$, one has
\begin{equation} \nonumber
\calW(\vartheta) \le C_0^{-1} \gotD_0^2 \, 16\be_0^{-2} \gotB |\gotD_0\e|^{-2 + 2(\gotm+1)-\gotm + \frac{3}{4}(2\gotm-1)} \le
\gotB^3 |\gotD_0\e|^{\frac{7}{4}} ,
\end{equation}
since $|\nu_{\ell_0}| > 4\gotm N$ and $|\nu_{\ell_2}| \le \gotm N$ imply 
$|\nu_{v_0}+\nu_{v_1}| > 3\gotm N$ and hence $p(\nu_{v_0}+\nu_{v_1}) \ge 3\gotm$,
so that, by Lemma \ref{lem:3.3}, one deduces $p(\nu_{v_0})+p_{\gotm}(\nu_{v_1}) \ge 2\gotm-1$.
\item[3.1.1.4.] 
If both $|\om_0\cdot\nu_{\ell_0}|$ and $|\om_0\cdot\nu_{\ell_1}|$ are less than $\be_0/4$ then one has
(a) $|\om_0\cdot(\nu_{v_0} + \nu_{v_2})| \le \be_0/2$ and hence $|\nu_{v_0} + \nu_{v_2}| > 4\gotm N$,
which, since $|\nu_{v_2}| \le \gotm N$,  implies $|\nu_{v_0}| > 3\gotm N$, i.e.~$p(\nu_{v_0}) \ge 3\gotm$,
and (b) $|\om_0\cdot\nu_{v_1}| > 4 \gotm N$, i.e. $p(\nu_{v_1})\ge 4\gotm$ and hence $p_{\gotm}(\nu_{v_1})\ge 3\gotm$.
Thus one has
\begin{equation} \nonumber
\calW(\vartheta) \le (C_0^{-1} \gotD_0^2)^2 \gotB |\gotD_0\e|^{-4 + 2 (\gotm+1) - \gotm + \frac{9}{4}\gotm + \frac{9}{4}\gotm}
\le \gotB^3 |\gotD_0\e|^{\frac{7}{2}} ,
\end{equation}
\end{itemize}
\item[3.1.2.] 
If $\lambda_{v_1}=\lambda_{v_2}=G$, one has
\begin{equation} \nonumber
\calW(\vartheta) \le |\calG_{\ell_0}| |\gotD_0\e|^{2(\gotm+1)-\gotm+\frac{3}{4}p(\nu_{v_0})} \gotB^2 .
\end{equation}
%
We distinguish between two possible subcases. 
\begin{itemize}[leftmargin=.31cm]
\itemsep0em
\item[3.1.2.1.] 
If $|\om_0\cdot\nu_{\ell_0}| \ge \be_0/4$ one has
\begin{equation} \nonumber
\calW(\vartheta) \le 16\be_0^{-2} \gotB^2 |\gotD_0\e|^{2(\gotm+1)-\gotm} \le
\gotB^3 |\gotD_0\e|^{3} .
\end{equation}
\item[3.1.2.2.] 
If $|\om_0\cdot\nu_{\ell_0}| < \be_0/4$ one has
\begin{equation} \nonumber
\calW(\vartheta) \le C_0^{-1} \gotD_0^2 \gotB^2 |\gotD_0\e|^{-2 + 2(\gotm+1)-\gotm + \frac{3}{2}\gotm } \le
\gotB^3 |\gotD_0\e|^{\frac{5}{2}} .
\end{equation}
since both $|\nu_{\ell_1}|$ and $|\nu_{\ell_2}|$ are $\le \gotm N$, while
$|\nu_{\ell_0}| > 4\gotm N$, so that one has $p(\nu_{v_0}) \ge 2\gotm$.
\end{itemize}
\end{itemize}
\item[3.2.] 
If $q=1$, the set $N(\vartheta)$ contains $v_0$, a second internal node $v_1$ and an end node $v_1'$, such that
$\lambda_{v_1'}=G$. One has
\begin{equation} \nonumber
\calW(\vartheta) \le |\calG_{\ell_0}|
|\gotD_0\e|^{\gotm+1-\gotm+\frac{3}{4} p(\nu_{v_0})}
|\calG_{\ell_1}|
|\gotD_0\e|^{\gotm+1-\gotm+\frac{3}{4} p(\nu_{v_1})}
\gotB .
\end{equation}
We distinguish between four subcases.
\begin{itemize}[leftmargin=.3cm]
\itemsep0em
\item[3.2.1.]
If $|\om_0\cdot\nu_{\ell_0}|, |\om_0\cdot\nu_{\ell_1}| \ge \be_0/4$, one obtains
\begin{equation} \nonumber
\calW(\vartheta) \le \gotB^3 |\gotD_0\e|^{2} .
\end{equation}
\item[3.2.2.]
If $|\om_0\cdot\nu_{\ell_0}| < \be_0/4$ and $|\om_0\cdot\nu_{\ell_1}| \ge \be_0/4$, then one has
$p(\nu_{v_0}) + p(\nu_{v_1}) \ge 3\gotm - 1$, so that 
\begin{equation} \nonumber
\calW(\vartheta) \le \gotB^3 |\gotD_0\e|^{- 2 + 2 + \frac{3}{4}(3\gotm -1)}\le \gotB^3 |\gotD_0\e|^{\frac{3}{2}}. 
\end{equation}
\item[3.2.3.]
If $|\om_0\cdot\nu_{\ell_0}| > \be_0/4$ and $|\om_0\cdot\nu_{\ell_1}| < \be_0/4$, then one has
$p(\nu_{v_1}) \ge 3\gotm$, which gives
\begin{equation} \nonumber
\calW(\vartheta) \le \gotB^3 |\gotD_0\e|^{- 2 + 2 + \frac{9}{4}\gotm } \le \gotB^3 |\gotD_0\e|^{\frac{9}{4}}. 
\end{equation}
\item[3.2.4.]
Finally, if $|\om_0\cdot\nu_{\ell_0}|, |\om_0\cdot\nu_{\ell_1}| < \be_0/4$, then one has
$p(\nu_{v_0})\ge 4\gotm$ and $p(\nu_{v_1}) \ge 3\gotm$, so that 
\begin{equation} \nonumber
\calW(\vartheta) \le \gotB^3 |\gotD_0\e|^{- 4 + 2 + 3\gotm + \frac{9}{4} \gotm} 
\le \gotB^3 |\gotD_0\e|^{\frac{13}{4}}. 
\end{equation}
\end{itemize}
\end{itemize}
Summarising, for $k=3$, one has $\calW(\vartheta) \le \gotB^3 |\gotD_0\e|^{\frac{1}{2} + |E_G(\vartheta)|}$ for all possible
values of $|E_G(\vartheta)|$, i.e.~$|E_G(\vartheta)|=1$ and $|E_G(\vartheta)|=2$.

\item[4.]
If $k \ge 4$, one has, by the inductive hypothesis and \eqref{eq:6.1},
\begin{equation} \nonumber
\begin{aligned}
\calW(\vartheta) & \le |\calG_{\ell_0}| |\gotD_0\e|^{(\gotm+1)q-\gotm + \frac{3}{4}p(\nu_{v_0}) } \gotB^{k-1}
|\gotD_0\e|^{\frac{k-1-q}{4}} |\gotD_0\e|^{|E_G(\vartheta)|-q_2-q_3} |\gotD_0\e|^{\frac{3}{4}q_2} , \\
& \le \gotB^{k} |\gotD_0\e|^{\frac{k-1}{4}} |\gotD_0\e|^{|E_G(\vartheta)|}
\left( \gotB^{-1} |\calG_{\ell_0}|
 |\gotD_0\e|^{\frac{3}{4}p(\nu_{v_0})  + \rho(q_1,q_2,q_3)} \right) ,
\end{aligned}
\end{equation}
with
\begin{equation}\nonumber 
\begin{aligned}
\rho(q_1,q_2,q_3)
& := (\gotm+1)q-\gotm - \frac{q}{4} - q_2 - q_3 + \frac{3}{4}q_2 \\
& =
\Bigl( \gotm + \frac{3}{4} \Bigr) q_1 + \Bigl( \gotm + \frac{1}{2} \Bigr) q_2 + \Bigl( \gotm - \frac{1}{4} \Bigr) q_3 - \gotm .
\end{aligned}
\end{equation}
%
\begin{itemize}[leftmargin=.31cm]
\itemsep0em
\item[4.1.] 
If $|\om_0\cdot\nu_{\ell_0}| \ge \be_0/4$, then one has $\gotB^{-1} |\calG_{\ell_0}|\le 1$. Moreover,
if $q_1+q_2=0$, one needs $q_3\ge 3$ in order to have
$k\ge 4$, so that $\rho(0,0,q_3) \ge 2\gotm-3/4>0$, while, if $q_1+q_2\ge 1$, one has
$\rho(q_1,q_2,q_3)\ge 1/2$. 
\item[4.2.] If $|\om_0\cdot\nu_{\ell_0}| < \be_0/4$, one has $\gotB^{-1} |\calG_{\ell_0}|\le |\gotD_0\e|^{-2}$,
we distinguish among the following cases.
%
%
\begin{itemize}[leftmargin=.3cm]
\itemsep0em
\item[4.2.1.] 
If $q_1+q_2=0$, once more $k\ge 3$ yields $q_3\ge 3$. If $q_3 \ge 4$, then
$\rho(0,0,q_3) - 2\ge 3(\gotm-1) \ge 0$; if $q_3=3$, then the fact that one has
$|\nu_{\ell_1}|,|\nu_{\ell_2}|,|\nu_{\ell_3}| \le \gotm N$ and $|\nu_{\ell_0}| > 4\gotm N$ implies
that $|\nu_{v_0}|> \gotm N$ and hence $p(\nu_{v_0})\ge \gotm$, so that
\begin{equation} \nonumber
\frac34 p(\nu_{v_0})+\rho(0,0,3) - 2 \ge \frac{11(\gotm -1)}{4} \ge 0 .
\end{equation}
\item[4.2.2.] 
If $q_1+q_2\ge 2$, one has $\rho(q_1,q_2,q_3) -2 \ge \gotm -1 \ge 0$.
\item[4.2.3.] 
If $q_1+q_2=1$, we distinguish between two cases.
\begin{itemize}[leftmargin=.6cm]
\itemsep0em
\item[4.2.3.1.\hspace{.3cm}] \hspace{-.4cm}
If $q_3 \ge 2$, still one has $\rho(0,0,q_3) - 2 \ge 2(\gotm -1) \ge 0$.
\item[4.2.3.2.\hspace{.3cm}] \hspace{-.4cm}
If $q_3=0,1$, then call $\ell_1',\ldots,\ell_{q'}'$ the $q'$ lines entering $v_1$, with $q'\ge 1$ (since $k\ge 4$), and
$\vartheta_1',\ldots,\vartheta_{q'}'$ the subtrees whose root lines are $\ell_1',\ldots,\ell_{q'}'$, respectively.
Once more, assume, without loss of generality, that 
the first $q_1'$ of such subtrees contain either at least two end nodes with end node labels $G$ or none;
the following $q_2'$ subtrees contain only one node with end node label $G$ and
are of order 2;
the last $q_3'$ subtrees contain only one node with end node label $G$ and are of order 1.
By construction, one has $q'= q_1' + q_2' + q_3'$, while
$|E_G(\vartheta)| = q_3 + |E_G(\vartheta_1')| + \ldots + |E_G(\vartheta_{q_1'}')| + q_2' +q_3'$,
if $q_1'\ge 1$, and $|E_G(\vartheta)| = q_3 + q_2' +q_3'$, if $q_1'=0$;
moreover the sum of the orders of the $q'$ subtrees equals $k-2-q_3$.
Then, using $(\gotm +1)q-\gotm=(\gotm+1)q_3+1$ for $q_1+q_2=1$ and $q_3=0,1$, one obtains
\begin{equation} \nonumber
\calW(\vartheta) 
\le \gotB^{k} |\gotD_0\e|^{\frac{k-1}{4}} |\gotD_0\e|^{|E_G(\vartheta)|}
\left( \gotB^{-1} |\calG_{\ell_1}|
|\gotD_0\e|^{ \frac{3}{4}(p(\nu_{v_0}) + p(\nu_{v_1})) + \rho_2(q_1',q_2',q_3',q_3) } \right) ,
\end{equation}
with
\begin{equation} \nonumber
\begin{aligned}
\rho_1(q_1',q_2',q_3',q_3) := & - 2 + (\gotm +1) q_3 + 1 + \rho(q_1',q_2',q_3') - \frac{1}{4} - \frac{5}{4}q_3  \\
= &
\Bigl( \gotm + \frac{3}{4} \Bigr) q_1' + \Bigl( \gotm + \frac{1}{2} \Bigr) q_2' + 
\Bigl(  \gotm - \frac{1}{4} \Bigr) q_3' - \gotm - \frac{5}{4} + \Bigl( \gotm -  \frac{1}{4} \Bigr) q_3 .
\end{aligned}
\end{equation}
\item[4.2.3.2.1.] 
If $|\om_0\cdot\nu_{\ell_1}|\ge \be_0/4$, then one has $\gotB^{-1} |\calG_{\ell_1}|\le 1$.
If $q_1'+q_2'=0$ and either $q_3=1$ (so that $q_3'\ge 1$) or $q_3=0$ (so that $q_3'\ge 2$),
one finds, for $q_3+q_3'=2$,
\begin{equation} \nonumber
\rho_1(q_1',q_2',q_3',q_3) + \frac{3}{4} (p(\nu_{v_0})+p(\nu_{v_1})) \ge 
2 \Bigl( \gotm - \frac{1}{4} \Bigr) - \gotm - \frac{5}{4} + \frac{3}{4} (2\gotm - 1) \ge \frac{5}{2}(\gotm-1)  ,
\end{equation}
since in such a case one has $|\nu_{v_0} + \nu_{v_1}| > 2\gotm$ and hence $p(\nu_{v_0})+p(\nu_{v_1}) \ge 2\gotm - 1$,
and, for $q_3+q_3' \ge 3$,
\begin{equation} \nonumber
\rho_1(q_1',q_2',q_3',q_3) + \frac{3}{4} (p(\nu_{v_0})+p(\nu_{v_1}) \ge 
3 \Bigl( \gotm - \frac{1}{4} \Bigr) - \gotm - \frac{5}{4}  \ge 2(\gotm-1) .
\end{equation}
If $q_1'+q_2' \ge 2$, one has
\begin{equation} \nonumber
\rho_1(q_1',q_2',q_3',q_3) \ge 
2 \Bigl( \gotm + \frac{1}{2} \Bigr) - \gotm - \frac{5}{4} \ge \gotm - \frac{1}{4} \ge \frac{3}{4} .
\end{equation}
If $q_1'+q_2' =1$, then either one has $q_3' \ge 1$, so that
\begin{equation} \nonumber
\rho_1(q_1',q_2',q_3',q_3) + \frac{3}{4} (p(\nu_{v_0})+p(\nu_{v_1}) \ge 
\Bigl( \gotm + \frac{1}{2} \Bigr) + \Bigl( \gotm - \frac{1}{4} \Bigr) 
- \gotm - \frac{5}{4}  \ge \gotm - 1 , 
\end{equation}
or there is only one line $\ell_1'$ entering $v_1$. In the latter case,
one iterates the construction once more:
call $\ell_1'',\ldots,\ell_{q''}''$ the $q''$ lines entering the node $v_1'$ which $\ell_1'$ exits, with $q''\ge 1$,
and $\vartheta_1'',\ldots,\vartheta_{q''}''$ the subtrees whose root lines are $\ell_1'',\ldots,\ell_{q''}''$, respectively.
Again, one may assume that the first $q_1''$ subtrees contain either at least two end nodes with end node labels $G$ or none;
the following $q_2'$ subtrees contain only one node with end node label $G$ and are of order 2;
the last $q_3'$ subtrees containi only one node with end node label $G$ and are of order 1.
By construction, one has $q''= q_1'' + q_2'' + q_3''$, while
$|E_G(\vartheta)| = q_3 + |E_G(\vartheta_1'')| + \ldots + |E_G(\vartheta_{q_1''}'')| + q_2'' +q_3''$,
if $q_1''\ge 1$, and $|E_G(\vartheta)| = q_3 + q_2'' +q_3''$, if $q_1''=0$;
finally the sum of the orders of the $q''$ subtrees equals $k-3-q_3$. Then one finds
\begin{equation} \nonumber
\calW(\vartheta) 
\le \gotB^{k} |\gotD_0\e|^{\frac{k-1}{4}} |\gotD_0\e|^{|E_G(\vartheta)|}
\left( \gotB^{-1} |\calG_{\ell_1'}|
|\gotD_0\e|^{ \frac{3}{4}(p(\nu_{v_0}) + p(\nu_{v_1})) + \rho_1(q_1'',q_2'',q_3'',q_3) + 1 - \frac{1}{4}} \right) .
\end{equation}
%
If $|\om_0\cdot\nu_{\ell_1'}|\ge \be_0/4$,
one has $\gotB^{-1} |\calG_{\ell_1'}|\le 1$ and, if $q_1''+q_2'' \ge 1$, one finds
\begin{equation} \nonumber
\rho_2(q_1'',q_2'',q_3'',q_3) + 1 - \frac{1}{4}  \ge 
\Bigl( \gotm + \frac{1}{2} \Bigr) - \gotm - \frac{5}{4} + 1 - \frac{1}{4} \ge 0,
\end{equation}
while, if $q_1''+q_2'' =0$, either one has $q_3''+q_3 \ge 2$, so that one finds
\begin{equation} \nonumber
\rho_2(q_1'',q_2'',q_3'',q_3) + 1 - \frac{1}{4}  \ge 
2 \Bigl( \gotm - \frac{1}{4} \Bigr) - \gotm - \frac{5}{4} + 1 - \frac{1}{4} \ge \gotm - 1 ,
\end{equation}
or $q_3''+q_3=1$, so that one has $|\nu_{v_0}+\nu_{v_1}+\nu_{v_1'}| > 3\gotm N$ and hence
$p(\nu_{v_0}) + p(\nu_{v_1}) + p(\nu_{v_1'}) \ge 3\gotm -2$, which gives
\begin{equation} \nonumber
\rho_2(q_1'',q_2'',q_3'',q_3) + 1 - \frac{1}{4}  \ge 
\Bigl( \gotm - \frac{1}{4} \Bigr) - \gotm - \frac{5}{4} + 1 - \frac{1}{4} + \frac{3}{4} (3\gotm-2) \ge \frac{9}{4} (\gotm-1) .
\end{equation}
If instead $|\om_0\cdot\nu_{\ell_1'}| < \be_0/4$, then one has $\gotB^{-1} |\calG_{\ell_1'}|\le |\gotD_0\e|^{-2}$,
but, in such a case, in order to have $|\nu_{\ell_1'}| < 4\gotm N$,
either $q_1''+q_2'' \ge 1$ or $p(\nu_{v_1}) \ge \max\{(4-q_3'')\gotm,0\}$.
Moreover, if $q_3=1$, then
$|\nu_{v_0} + \nu_{v_1} 
+ \nu_{\ell_2}| > 4\gotm N$, which, since $|\nu_{\ell_2}| \le \gotm N$, implies
$|\nu_{v_0} + \nu_{v_1}| > 3\gotm N$ and hence $p(\nu_{v_0}) + p(\nu_{v_0}) \ge 3\gotm - 1$, by Lemma \ref{lem:3.3},
while, if $q_3=0$, then
$|\nu_{v_0} + \nu_{v_1}| > 4\gotm N$ and hence $p(\nu_{v_0}) + p(\nu_{v_0}) \ge 4\gotm -1 $, so that,
for both values of $q_3$, one finds
\begin{equation} \nonumber
\begin{aligned}
& - 2 + \rho_2(q_1'',q_2'',q_3'',q_3) + 1 - \frac{1}{4}  + \frac{3}{4} \left( p(\nu_{v_0}) + p(\nu_{v_1}) \right) \\
& \qquad
\ge - 2 + \min \left\{ \gotm + \frac{1}{2} ,  \frac{3}{4} (4-q_3'')\gotm
+ \Bigl( \gotm - \frac{1}{4}  \Bigr) q_3'' \right\}
- \gotm - \frac{5}{4} + \Bigl( \gotm - \frac{1}{4}  \Bigr) q_3  \\
& \qquad \qquad + 1 - \frac{1}{4} + \frac{3}{4} \left( (4 - q_3) \gotm - 1 \right)
\ge \frac{12 \gotm - 11}{4} + \frac{\gotm-1}{4} q_3 \ge \frac{1}{4} .
\end{aligned}
\end{equation}
\item[4.2.3.2.2.] 
If $|\om_0\cdot\nu_{\ell_1}| < \be_0/4$,
then one has $\gotB^{-1} |\calG_{\ell_1}|\le |\gotD_0\e|^{-2}$. On the other hand
one has $\nu_{\ell_0} = \nu_{v_0} + \nu_{\ell_1} + \nu_{\ell_2}$
if $q_3=1$, and $\nu_{\ell_0} = \nu_{v_0} + \nu_{\ell_1}$ if $q_3=0$; in the first case, since one has 
$|\nu_{\ell_2}|\le \gotm N$ and $|\nu_{v_0}+\nu_{\ell_2}| > 4\gotm N$, one finds
$|\nu_{v_0}| > 3\gotm N$ and hence $p(\nu_{v_0}) \ge 3\gotm$,
while in the second case one $|\nu_{v_0}| > 4\gotm N$ and hence $p(\nu_{v_0}) \ge 4\gotm$.
Moreover, one has either $q_1'+q_2' \ge 1$ or $q_3'+q_3 \ge 2$
(because $k \ge 4$) and $p(\nu_{v_1}) \ge \max\{(4-q_3')\gotm,0\}$.
In summary, one finds
\begin{equation} \nonumber
\begin{aligned}
& - 2 + \rho_2(q_1',q_2',q_3',q_3) + \frac{3}{4}(p(\nu_{v_0}) + p(\nu_{v_1})) \\
& \qquad \ge
\min \left\{ \Bigl( \gotm + \frac{3}{4} \Bigr) q_1' + \Bigl( \gotm + \frac{1}{2} \Bigr) q_2' ,
\frac{3}{4} (4-q_3') \gotm + \Bigl(  \gotm - \frac{1}{4} \Bigr) q_3' \right\}  \\
& \qquad \qquad
- \gotm - \frac{13}{4}
+ \Bigl( \gotm -  \frac{1}{4} \Bigr) q_3 + \frac{3}{4}  (4-q_3) \gotm ,
\end{aligned}
\end{equation}
so that one obtains, for both values of $q_3$,
\begin{equation} \nonumber
- 2 + \rho_2(q_1',q_2',q_3',q_3) + \frac{3}{4} \left( p(\nu_{v_0}) + p(\nu_{v_1}) \right) \ge
\Bigl(  \gotm + \frac{1}{2} \Bigr) + 2\gotm -  \frac{13}{4} \ge \frac{12\gotm-11}{4} .
\end{equation}
\end{itemize}
\end{itemize}
\end{itemize}
\end{itemize}
This concludes the proof. \qed

\begin{lemma} \label{lem:6.2}
Assume Hypotheses \ref{hyp:1}. 
Let $\gotD_0$ and $\e_2$ be defined as in Lemmas \ref{lem:3.9} and \ref{lem:3.11}, respectively,
and assume that $|\e| \le \e_2$.
For all $\vartheta\in\gotT_{k,0}$, with $k\ge 2$, $\nu\in\ZZZ^2_*$ and $|E_G(\vartheta)|\ge 1$, one has
\begin{equation} \nonumber
\calW(\vartheta) \le 
\begin{cases}
\ga^{-1} \gotB^{k-1} |\gotD_0\e|^{\frac{k+1}{4}} |\gotD_0\e|^{|E_G(\theta)|} , & \quad k \ge 3 , \; |E_G(\vartheta) | \ge 
1 ,\\
\ga^{-1} \gotB |\gotD_0\e| , 
& \quad k = 2 , \; |E_G(\vartheta)| = 1 .
\end{cases}
\end{equation}
\end{lemma}

\noindent\emph{Proof}. 
Let $\ell_1,\ldots,\ell_q$ denote the $q$ lines entering $v_0$ which the root line $\ell_0$ exits, and
$\vartheta_1,\ldots,\vartheta_q$ the subtrees whose root lines are $\ell_1,\ldots,\ell_q$, respectively.
As in the proof of Lemma \ref{lem:6.1}, assume for definiteness that 
the first $q_1$  subtrees contain either at least two end nodes with end node labels $G$ or none;
the following $q_2$ subtrees contain only one node with end node label $G$ and
are at least of order 2;
the last $q_3$ subtrees contain only one node with end node label $G$ and are of order 1.
By construction, one has $q = q_1 + q_2 + q_3 \ge 1$ (since $E_G(\vartheta) \ge 1$) and
$|E_G(\vartheta)|$ equals $q_2 + q_3$ if $q_1=0$, while it is given by \eqref{eq:6.1} if $q_1\ge 1$.
\begin{itemize}[leftmargin=.5cm]
\itemsep0em
\item[1.]
For $k=2$ one has $q=1$ and the $\ell_1$ exits an end node $v_1$ with $\lambda_{v_1}=G$.
Then one has $\calW(\vartheta) \le \ga^{-1} |\gotD_0\e|\gotB$.
\item[2.]
%
%
%
%
%
%
%
%
%
For $k\ge 3$ and $|E_G(\vartheta)|\ge 1$, one has
\begin{equation} \nonumber
\calW(\vartheta) =  \gamma^{-1} |\gotD_0\e|^{(\gotm+1)q - \gotm}  \prod_{i=1}^{q} \calW(\vartheta_i) 
\le \gamma^{-1} \gotB^{k-1} |\gotD_0\e|^{\frac{k}{4} + \rho_0(q_1,q_2,q_3)} |\gotD_0\e|^ {E_G(\vartheta)| }
\end{equation}
with
\begin{equation} \nonumber
\begin{aligned}
\rho_0(q_1,q_2,q_3) & := (\gotm+1)q - \gotm - \frac{1+q}{4} - q_3 - \frac{1}{4}q_2 \\
& = \Bigl( \gotm+\frac{3}{4} \Bigr) q_1 + \Bigl( \gotm+\frac{1}{2} \Bigr) q_2 +
\Bigl( \gotm - \frac{1}{4} \Bigr) q_3 - \gotm - \frac{1}{4} .
\end{aligned}
\end{equation}
If $q_1+q_2 \ge 1$ one has
\begin{equation} \nonumber
\rho_0(q_1,q_2,q_3) \ge \Bigl( \gotm + \frac{1}{2} \Bigr) (q_1+q_2) - \gotm - \frac{1}{4} 
\ge \gotm + \frac{1}{2} - \gotm - \frac{1}{4} \ge \frac{1}{4} ,
\end{equation}
while, if $q_1+q_2=0$ then, in order to have $k\ge 3$, one needs $q_3 \ge 2$ and hence
\begin{equation} \nonumber
\rho_0(q_1,q_2,q_3) \ge \Bigl( \gotm - \frac{1}{4} \Bigr) q_3 - \gotm - \frac{1}{4}
\ge \Bigl( \gotm - \frac{1}{4} \Bigr) 2 - \gotm - \frac{1}{4} \ge \gotm - \frac{3}{4} \ge \frac{1}{4} ,
\end{equation}
so that the last bound follows in both cases. \qed
\end{itemize}


Define $\TT_{k_{1},k_{2},\nu}$ as the set of labelled trees $\theta \in \gotT_{k_{1}+k_{2},\nu}$ such that
$|E_F(\vartheta)|+|I(\vartheta)|=k_{1}$ and $|E_G(\vartheta)|=k_{2}$. Write
\begin{equation} \nonumber 
H^{(k_{1},k_{2})}_{\nu} := \!\!\!\!\!\sum_{\vartheta \in \gotT_{k_{1},k_{2},\nu}} \!\!\!\!\!
 \Val(\vartheta), \qquad \nn \in \ZZZ^2_* , \qquad
\al^{(k_{1},k_{2})} :=  \!\!\!\!\! \sum_{\vartheta \in \gotT_{k_{1},k_{2},0}}  \!\!\!\!\! \Val(\vartheta) ,
\end{equation}
so that \eqref{eq:3.9} is replaced with
\begin{equation} \label{eq:6.2}
H(\psi) = \sum_{\substack{k_{1},k_{2} \ge 1 \\ k_{1} + k_{2} \ge 1}}
\sum_{\nu\in\ZZZ^2} \mu_1^{k_{1}} \mu_2^{k_{2}}{\rm e}^{i\nu\cdot \psi} H^{(k_{1},k_{2})}_{\nu} , \qquad
\tilde \al = \sum_{\substack{k_{1},k_{2} \ge 1 \\ k_{1} + k_{2} \ge 1}}
\mu_1^{k_{1}} \mu_2^{k_{2}}  \al^{(k_{1},k_{2})} ,
\end{equation}
where $\mu_1$ and $\mu_2$ are two parameters that eventually have to be set equal to 1.
As previously for $\mu$, this is proved to be allowed by showing that the radius of convergence
in both variables $\mu_1$ and $\mu_2$ is greater than 1, provided that $\e$ is small enough.

Finally, set 
$\gotC_{1}^{*} \!:=\! 2\gotD_{1}^{-1} N \gotR_0 \gotB_1$,
with the constants $\gotD_{1}$, $\gotB_1$ and $\gotR_0$
as in the proof of Lemma \ref{lem:5.1},  and define
\begin{equation} \label{eq:6.3}
\e_6 := \frac{1}{(12\Phi\gotC_{1}^{*})^4\gotD_{1}} = \frac{\gotD_{1}^3}{(24N \Phi \gotB_1 \gotR_0)^4} .
\end{equation}
Note that $\e_6$ does not depend on $\ga$. Also set
$\Gamma_1 :=4 \ga\be^{-1} \Phi_0\gotD_{1}^{-\gotm}$, and write $\Gamma_1=\bar\Gamma_1 |\gotD_{1}\e|^{-a}$,
with $\bar\Gamma_1 :=4\bar\ga\be^{-1} \Phi_0\gotD_{1}^{-\gotm}$ independent of $\e$ and $a=\gotm - s$, 
where $\bar \ga$ and $s$ are as in Theorem \ref{thm:2}.

Lemmas \ref{lem:4.2}, \ref{lem:4.3}, \ref{lem:6.1} and \ref{lem:6.2} immediately imply  the following result.

\begin{lemma} \label{lem:6.3}
Assume Hypotheses \ref{hyp:1} and \ref{hyp:2}, and let $\be$ be defined as in \eqref{eq:5.1}.
For all $k_{1}\in\NNN$, all $|\e| \le \e_6$ and all $|\al_0-\al | \le \eta_0$,
with $\e_6$ and $\eta_0$ defined in \eqref{eq:6.3} and in \eqref{eq:5.4}, respectively, one has,
\vspace{-.2cm}
\begin{equation} \nonumber
\begin{aligned}
| H^{(k_{1},0)}_{\nu} | & \le 
\gotD_{1}^{\gotm+1} N^{-1} (\gotC_{1}^{*})^{k_{1}} (6\Phi)^{k_{1}} 
|\gotD_{1}\e|^{\frac{k_{1}-1}{4}} , \\
| H^{(0,1)}_{\nu} | & \le 
\gotD_{1}^{\gotm+1} N^{-1} \gotC_{1}^{*} \Gamma_1 
, \\
| H^{(1,1)}_{\nu} | & \le 
\gotD_{1}^{\gotm+1} N^{-1} (\gotC_{1}^{*})^{2} 6\Phi \Gamma_1
|\gotD_{1}\e| , \\ 
%
| H^{(k_{1},k_{2})}_{\nu} | & \le 
\gotD_{1}^{\gotm+1} N^{-1} (\gotC_{1}^{*})^{k_{1}+k_{2}} (6\Phi)^{k_{1}} \Gamma_1^{k_{2}}
|\gotD_{1}\e|^{\frac{k_{1}+k_{2}-1}{4}} |\gotD_{1}\e|^{k_{2}} ,
\qquad k_{2} \ge 1 , \\
\end{aligned}
\end{equation}
for all $\nu\in\ZZZ^2_*$ with $p(\nu) < \gotm (k_{1}+k_{2}) + k_1 + k_2 -1$. The same bounds hold, but multiplied by
$|\gotD_{1}\e|^{\frac{1}{8}(p(\nu)-(\gotm+1)(k_{1}+k_{2})-(k_1+k_2-1)}$, if $p(\nu) \ge (\gotm+1)(k_1+k_2) + (k_{1}+k_{2}-1)$, and
%
\begin{equation} \nonumber
\begin{aligned}
|\al^{(1,0)}| & 
\le \ga^{-1} 6 \Phi \gotD_{1}^{\gotm} |\gotD_{1}\e| , \qquad \quad \\
|\al^{(k_{1},0)}| & 
\le \ga^{-1} B_1^{-1} \gotD_{1}^{\gotm+1}N^{-1}  (\gotC_{1}^{*})^{k_{1}} (6\Phi)^{k_{1}} 
|\gotD_{1}\e|^{\frac{k_{1}+2}{4}}  ,
\hspace{2.80cm} \qquad k_{1}\ge 2 , \\ 
|\al^{(1,1)}| & 
\le \ga^{-1} B_1^{-1} \gotD_{1}^{\gotm+1} N^{-1} (\gotC_{1}^{*})^2 6\Phi \Gamma_1 |\gotD_{1}\e| , \qquad \quad \\
%
%
|\al^{(k_{1},k_{2})}| & 
\le \ga^{-1} B_1^{-1} \gotD_{1}^{\gotm+1}N^{-1}  (\gotC_{1}^{*})^{k_{1}+k_{2}} (6\Phi)^{k_{1}} \Gamma_1^{k_{2}}
|\gotD_{1}\e|^{\frac{k_{1}+k_{2}+1}{4}} |\gotD_{1}\e|^{k_{2}} , 
\qquad 
k_{2} \ge 2 .
\end{aligned}
\end{equation}
Moreover there exists a positive constant $\e_7 \le \e_6$ such that, for all $|\e| \le \e_7$,
the series \eqref{eq:6.3} converge with $\mu_1=\mu_2=1$.
\end{lemma}

\noindent\emph{Proof}. 
For fixed $\al_0$, one reasons as in the proof of Lemma \ref{lem:4.8}, using the bounds in Lemmas \ref{lem:6.1} and \ref{lem:6.2}.
The result is then extended to any $\al_0\in\calA(\be)$ as outlined in Remark \ref{rmk:4.12}.
The proof that the radius of convergence of the series \eqref{eq:3.9} in both
variables $\mu_1$ and $\mu_2$ is greater than $1$, proceeds as in the proof of Lemma \ref{lem:4.10} -- by taking into
account once more Remark \ref{rmk:4.12} to extend the argument to any $\al_0\in\calA(\be)$ -- provided that
$|\e| \le \e_7$, with
\vspace{-.2cm}
\begin{equation} \label{eq:6.4}
\e_7 := \min \{ \tilde \e_7, \e_6 \} , \qquad \tilde \e_7 := \frac{1}{(2 \gotC_{1}^{*}\bar\Gamma_1)^4\gotD_{1}} ,
\vspace{-.2cm}
\end{equation}
with the condition $|\e| \le \tilde\e_7$ guaranteeing summability over $k_2$.
\qed

\begin{lemma} \label{lem:6.4}
Assume Hypotheses \ref{hyp:1} and \ref{hyp:2}, and let $\be$ be defined as in \eqref{eq:5.1}.
Let $\eta_0$ and $\e_7$ be as in Lemmas \ref{lem:5.1} and \ref{eq:6.3}, respectively.
There exist positive constants $a_0^*$, $a_1^*$ and $a_2^*$ such that, for all $|\e|< \e_7$
and all $|\al_0 - \al| < \eta_0$, one has
\vspace{-.2cm}
\begin{equation} \nonumber
|\tilde\al(\e,\al_0)| \le a_0^* |\e| , \qquad
\left| \frac{\partial \tilde\al}{\partial \al_0}(\e,\al_0) \right| \le a_1^* |\e| , \qquad
\left| \frac{\partial \tilde\al}{\partial \e}(\e,\al_0) \right| \le a_2  .
\vspace{-.1cm}
\end{equation}
\end{lemma}

\noindent\emph{Proof}. 
The first bound is proved as in Lemma \ref{lem:5.1}, the only difference being that now there is a double sum,
over $k_{1}$ and $k_{2}$, with $k_{1} \ge 1$ and $k_{2}\ge 0$. The two sums converge, for all $a\in[0,1)$,
provided that  $|\e|\le \e_7$, as discussed in the proof of Lemma \ref{lem:6.3}. Thus, the first bound follows with
\vspace{-.2cm}
\begin{equation} \label{eq:6.5}
a_0^* := \ga^{-1} 6 \Phi \gotD_{1}^{\gotm+1} \left(  1 + \frac{\gotK_0 N \gotB_1 \gotR_0^2}{3\Phi \gotD_{1}} 
\right) ,
\vspace{-.2cm}
\end{equation}
where 
%
\begin{equation} \nonumber
\gotK_0 := 2 (6\Phi)^2 + 6 \Phi \Gamma_1 
+ 4 (6\Phi) \Gamma_1 |\gotD_1\e| .
\end{equation}
%
Note that $\ga^{-1}$ and $\gotK_0$ are proportional to $|\e|^{a}$ and
to $|\e|^{-a}$, respectively, so that $a_0$ is bounded by an $\e$-independent constant.

The second and third bounds are established as in the proof of Lemma \ref{lem:5.3}. The only \emph{caveat}
is that, when differentiating with respect to $\al_0$ the propagator $\calG_\ell$ of a line $\ell\in L(\vartheta)$,
even though one still bounds in general (i.e.~for any $\nu_{\ell}$)
\begin{equation} \nonumber
\left| \frac{\partial \calG_{\ell}}{\partial \al_0} \right|
\le
2 |\calG_{\ell}| \left( \frac{1}{|\e|\sqrt{2A_0}} + \frac{|\e|^{\gotm}\ga}{|\e|^{2}A_0} + 4\be_0^{-1} N \right)
\le \frac{Q_0}{|\gotD_0\e|} |\calG_{\ell}| \, |\nu_{\ell}| ,
\end{equation}
with $Q_0$ as in the proof of Lemma \ref{lem:5.3}, hence proportional to $|\e|^{\gotm-1-a}$,
%
%
however, when $|\om_0\cdot\nu_{\ell}| \ge \be_0/4$, the bound can be replaced by
\begin{equation} \nonumber
\left| \frac{\partial \calG_{\ell}}{\partial \al_0} \right| \le
Q_0^* |\calG_{\ell}| |\nu_{\ell}| , \qquad
Q_0^* := \frac{8}{\be_0} + \frac{16|\e|^{\gotm}\ga}{\be_0^2} + \frac{4N 
}{\be_0} ,
\end{equation}
with $Q_0^*$ uniformly bounded (in $\e$) for $a\in[0,1)$.
On the other hand, if a tree $\vartheta$ of order $k$ is such that $E_G(\vartheta) \ge 1$, then
either $k=2$ and the line exiting the end node has momentum $\le \gotm N$ or
$k\ge 3$ and in that case $\VV(\vartheta)$ contains at least an extra factor $|\gotD_0\e|$. 

Eventually this leads to the bound $a_1^*|\e|$, with
\begin{equation} \label{eq:6.6}
a_1^* := 
\ga^{-1} \gotD_{1}^{\gotm+1} \left( \frac{96 \gotm^2 N \Phi}{\be} +
\frac{24 (\gotm+1) N^2 \gotB_1 \gotR_0^2}{\gotD_{1} } \left( 
\gotK_{0} Q_0^* + \frac{8 \gotK_2 Q_1 N \gotB_1 (6\Phi) \gotR_0}{3\gotD_{1}} 
\right) \right) ,
\end{equation}
where $Q_1$ is defined as in the proof of Lemma \ref{lem:5.3}, with $\ga=\bar\ga|\e|^{-a}$, and
%
%
%
\begin{equation} \nonumber
\gotK_{0} := (6\Phi)^2 + 6\Phi \, \Gamma_1 , \qquad
\gotK_2 := 2 (6\Phi)^2 
+ 4 (6\Phi) \Gamma_1 |\gotD_1\e| .
\end{equation}
In particular, the constants $\ga$, $\gotK_{0}$ and $Q_1$ are proportional to $|\e|^{-a}$ (the last one for $\gotm=1$),
while $Q_0^*$ and $\gotK_2$ are bounded by an $\e$-independent constant for all $a\in[0,1)$
so that $a_1^*$ is also bounded independently of $\e$. 

The last bound is also discussed as in the proof of Lemma \ref{lem:5.3}, and is proved to hold with
\begin{equation} \label{eq:6.7}
a_2^* := 24 \ga^{-1} \gotm \Phi \gotD_{1}^{\gotm+1} 
\left( 1 + \frac{\gotK_0 Q_2 N \gotB_1 (6\Phi)^2 \gotR_0^2}{\gotm \Phi \gotD_{1}} \right) ,
\end{equation}
with $\gotK_0$ as in \eqref{eq:6.5} and $Q_2$ as in the proof of Lemma \ref{lem:5.3}. \qed

\vspace{.3cm}

A result analogous to Lemma \ref{lem:5.5} follows, with only difference being that the further condition to be required on $\e$
involves the parameters $a_0^*$, $a_1^*$ and $a_2^*$.

\begin{lemma} \label{lem:6.5}
Assume Hypotheses \ref{hyp:1} and \ref{hyp:2}, and let $\be$ be defined as in \eqref{eq:5.1}.
Then there exist constants $\al_0$ and $\bar\e_0$ such that, defining
$\om_0$ and $\be_0$ as in Lemma \ref{lem:5.5},
for all $|\e|\le \bar\e_0$ there exists a counterterm $\tilde\al(\e,\al_0)$, differentiable in both $\e$ and $\al_0$,
and a function $H(\psi;\e)$, continuous in $\e$ and analytic in $\psi$, such that
$\al_0 + \e \tilde \al(\e,\al_0)=\al$, $\be_0 \ge \be/2$ and
$H(\om_0 t;\e)$ solves \eqref{eq:2.8}.
\end{lemma}

\noindent\emph{Proof}. 
One reasons as in the proof of Lemma \ref{lem:5.5}, by defining
\begin{equation} \label{eq:6.8}
\bar \e_0 = \min\{\e_6,\tilde\e_8\} , \qquad \tilde\e_8 := \min\left\{  \sqrt{\frac{1}{2a_1^*}} , \sqrt{\frac{\eta_0}{2(a_0^*+a_2^*)}} \right\} ,
\end{equation}
with $a_0^*$, $a_1^*$ and $a_2^*$ defined as in Lemma \ref{lem:6.3}.
\qed

\vspace{.3cm}

The following result, analogous to Lemma \ref{lem:5.6}, completes the proof of Theorem \ref{thm:2}.

\begin{lemma} \label{lem:6.6}
Assume Hypotheses \ref{hyp:1} and \ref{hyp:2}, and let $\be$ be defined as in \eqref{eq:5.1}.
Let $\bar\e_0$ be defined as in Lemma \ref{lem:6.5}.
Then for all $|\e| \le \bar\e_0$ there is a multi-periodic solution $\theta(t)$ to \eqref{eq:2.1},
with $f$ given by \eqref{eq:2.2}, 
of the form $\theta(t)=\theta_0+\al_0 t + \Theta(\om_0 t;\e)$, with $\theta_0\in\TTT^2$ and $\om_0=(\al_0,1)$, such that 
\begin{equation} \nonumber
\lim_{\e\to0} \al_0 = \al , \qquad
\lim_{\e\to0} \Theta(\psi;\e) = 0 .
\end{equation}
\end{lemma}

\noindent\emph{Proof}. 
As for Lemma \ref{lem:5.6}. \qed

\zerarcounters 
\section{Conclusions}
\label{sec:7} 

We have considered a class of systems on the one-dimensional  torus $\TTT$, described by \eqref{eq:1.2},
and have shown that, as long as $\ga_0$ is sublinear in $\e$ and $f$ is a trigonometric polynomial,
if $\al$ satisfies a non-resonance condition up to some finite order,
depending on both the degree of the polynomial and the size of $\ga_0$,
then for $|\e|$ small enough there exists a multi-periodic solution (cf.~Theorem \ref{thm:2}).
The frequency vector of the multi-periodic solution is $(\al_0,1)$, 
where $1$ is the fixed intrinsic frequency of $f$ in the variable $t$ and $\al_0$ is such that
\begin{equation} \label{eq:7.1}
\al_0+\e \tilde \al(\e,\al_0)=\al ,
\end{equation}
for a suitable function $\tilde \al(\e,\al_0)$. The function $\tilde \al(\e,\al_0)$ is smooth in both $\e$ and $\al_0$,
as long as $\e$ is small and $\al_0$ is close to $\al$, so that \eqref{eq:7.1} can be solved by the implicit function theorem.

The multi-periodic solution has two frequencies $\al_0$ and $1$;  it is quasi-periodic if $\al_0$
is incommensurate with 1 and periodic otherwise. In the latter case, the period is expected to be large,
so that, up to very large time scales, the solution, for any practical purpose, appears to be quasi-periodic
even when it is periodic; for this reason, it is said to be pseudo-synchronous.

The maximum size $\e_0$ allowed for $\e$ depends on the degree $N$ of the trigonometric polynomial,
on the size of the dissipation and on the parameter $\al$. For fixed values of $\ga_0$, one may write
$\ga_0=|\e|^{\gotm -a}\bar\ga$, for some $a\in[0,1)$ and $\bar\ga$ independent of $\e$,
and require a non-resonance condition of order $4\gotm N$
on the frequency vector $\om=(\al,1)$:
then $\e_0$ goes to $0$ when either $N$ or $\gotm$ goes to infinity.
The fact that $\e_0$ goes to $0$ with $\gotm$ going to infinity is unavoidable:
when the dissipation goes to zero, the system becomes a quasi-integrable
Hamiltonian system and the resonant tori are destroyed.

Thus, the result leaves several open problems,
which would require further investigation:
\vspace{-.2cm}
\begin{enumerate}
\itemsep0em
\item in principle, the assumption that $f$ is a trigonometric polynomial might be removable,
possibly by requiring additional hypotheses on $\al$;
\item the quasi-periodic solutions proved to exist are expected to be local attractors, on both physical and mathematical grounds,
but this should be studied explicitly;
\item in order to apply the results to real situations, one should use for the parameters $\e$, $\ga$ and $\al$ in \eqref{eq:1.2}
the values arising from astronomical data, and check whether they are allowed for the theorems to hold.
\item while the quasi-periodic solutions investigated here bifurcate from an unperturbed maximal KAM torus,
it would be interesting to extend the results to the Langrangean tori appearing inside the resonant gaps.
\vspace{-.2cm}
\end{enumerate}

With respect to the existing literature, the main advantage of Theorem \ref{thm:2} is that
the external parameter does not need to be tuned, but can be set equal to all except finitely many values,
and that the frequencies are required to satisfy very mild non-resonance conditions.
On the other hand, the main drawback is that only trigonometric polynomials can be considered.
So, a natural question, according to item 1, is whether one can extend the result so as to include
analytic functions with arbitrarily many harmonics, and, if the answer is positive,
whether one must require stronger non-resonance conditions on $\al$ than those allowed by Hypothesis \ref{hyp:2}.

Already in the case discussed in this paper,
attractiveness of the solution is not obvious (see item 2),
since one cannot rely on the normal form construction used in
refs.~\cite{SL,M}, where the frequency vector $\om_0$ satisfies a Diophantine condition.
Linearisation around the multi-periodic solution
leads to Hill's equation with multi-periodic potential and no Diophantine assumption on the frequency vector.

As far as item 3 is concerned, we note that even in the case of the periodic attractors comparison with the astronomical data
is quite non-trivial \cite{ABC,BC1}.
Moreover, as said in the introduction, the attractors found in problems of celestial mechanics essentially are either periodic orbits or
Lagrangian tori, so that a direct application to real situations is far from obvious. The case of KAM tori could be
seen as a first, preliminary step before attacking the (harder) case of Lagrangian tori (see item 4).
Note that the Lagrangian tori which appear inside the gaps among the surviving tori are rather irrelevant
in the conservative case, since they cover in aggregate a negligible fraction of the phase space
(see refs.~\cite{MNT,BC2} for recent results); in the dissipative case, instead, they can play a
fundamental role, since the results in the literature suggest that, when a more realistic
model for the tidal torque is considered, a Lagrangian torus survives as the dominant attractor \cite{NFME,BDG4}.

\vspace{1cm}

\noindent \textbf{Ackowledgments.} We thank Luigi Chierchia and Jessica Massetti for useful discussions.


\end{document}